\documentclass[12pt]{amsart}

\usepackage{bm,latexsym,amsfonts,amsmath,amssymb,amsthm,url,graphics,psfrag,amscd,stmaryrd, mathrsfs}
\usepackage[english]{babel}
\usepackage[T1]{fontenc}
\usepackage{float}
\usepackage{graphics}
\usepackage{graphicx}
\usepackage{color}
\usepackage[colorinlistoftodos]{todonotes}  
\usepackage{hyperref}
\usepackage{comment}

\newcommand{\R}{\mathbb{R}}

\newcommand{\N}{\mathbb{N}}
\newcommand{\E}{\mathbb{E}}
\newcommand{\Z}{\mathbb{Z}}

\newcommand{\pp}{\mathbb{P}}

\newcommand{\kZ}{\mathcal{Z}}

\newcommand{\bx}{\boldsymbol{x}}

\newcommand{\by}{\boldsymbol{y}}

\newcommand{\bB}{\boldsymbol{B}}
\newcommand{\bs}{\boldsymbol{\sigma}}
\newcommand{\lin}{\left[\kern-0.15em\left[}
\newcommand{\rin} {\right]\kern-0.15em\right]}
\newcommand{\linf}{[\kern-0.15em [}
\newcommand{\rinf} {]\kern-0.15em ]}
\newcommand{\ilin}{\left]\kern-0.15em\left]}
\newcommand{\irin} {\right[\kern-0.15em\right[}

\newcommand{\bpr}{\beta^\prime}

\newcommand{\TTP}{{\sf TTP}}
\newcommand{\TP}{{\sf TP}}

\newcommand{\bsigma}{\boldsymbol{\sigma}}

\newtheorem{theorem}{Theorem}[section]

\newtheorem{definition}[theorem]{Definition}
\newtheorem{lemma}[theorem]{Lemma}

\newtheorem{proposition}[theorem]{Proposition}

\newtheorem{remark}[theorem]{Remark}
\newtheorem{cond}[theorem]{Condition}

\newcommand{\be}{\begin{equation}}
\newcommand{\ee}{\end{equation}}
\newcommand{\beal}{\begin{align}}
\newcommand{\eal}{\end{align}}
\newcommand{\nn}{\nonumber}

\newcommand{\prob}{\mathbb P}
\newcommand{\expec}{\mathbb E}

\newcommand{\calF}{\mathscr{F}}

\newcommand{\sss}{\scriptscriptstyle}

\numberwithin{equation}{section}
\newcommand{\vep}{\varepsilon}
\newcommand{\e}{{\rm e}}

\newcommand{\rem}[1]{}

\newcommand{\NM}[1]{\todo[inline,color=cyan]{#1}}
\newcommand{\RvdH}[1]{\todo[inline, color=magenta]{{\rm Remco: #1}}}
\newcommand{\Gib}[1]{\todo[inline, color=yellow]{Claudio: #1}}

\def\1{{\mathchoice {1\mskip-4mu\mathrm l}      
{1\mskip-4mu\mathrm l}
{1\mskip-4.5mu\mathrm l} {1\mskip-5mu\mathrm l}}}
\newcommand{\indic}[1]{\1_{\{#1\}}}
\newcommand{\indicwo}[1]{\1_{#1}}

\newcommand{\GRGw}{\mathrm{GRG}_n(\boldsymbol{w})}

\newcommand{\eqn}[1]{\begin{equation}#1\end{equation}}
\newcommand{\eqan}[1]{\begin{align}#1\end{align}}

\newcommand{\bfwit}{\boldsymbol{w}}

\newcommand{\wpn}{W^\prime_n}

\newcommand{\s}{\sigma}
\newcommand{\omqn}{[q]^n}
\newcommand{\omq}{[q]}
\newcommand{\bpinv}{\frac{1}{\beta'}}

\newcommand{\invisible}[1]{}
\newcommand{\ulW}{1}
\newcommand{\olW}{w}

\DeclareSymbolFont{extraup}{U}{zavm}{m}{n}
\DeclareMathSymbol{\varheart}{\mathalpha}{extraup}{86}
\DeclareMathSymbol{\vardiamond}{\mathalpha}{extraup}{87}
\newcommand{\ensymboldefinition}{$\blacktriangleleft$}

\newcommand{\Nel}[1]{\todo[inline, color=orange]{\color{black}Neeladri: #1 }}

\usepackage{fancyhdr}
\begin{document}

\title{Annealed Potts models on\\rank-1 inhomogeneous random graphs}
\author{Cristian Giardin\`a}
\address{Cristian Giardin\`a, University of Modena and Reggio Emilia, via G. Campi 213/b, 41125 Modena, Italy, cristian.giardina@unimore.it}
\author{Claudio Giberti}
\address{Claudio Giberti, University of Modena and Reggio Emilia, via G. Amendola 2, Pad. Morselli, 42122 Reggio E., Italy, claudio.giberti@unimore.it}
\author{Remco van der Hofstad}
\address{Remco van der Hofstad, Eindhoven University of Technology, P.O. Box 513, 5600 MB Eindhoven, The Netherlands, r.w.v.d.hofstad@tue.nl}
\author{\\Guido Janssen}
\address{Guido Janssen, Eindhoven University of Technology, P.O. Box 513, 5600 MB Eindhoven, The Netherlands, a.j.e.m.janssen@tue.nl}
\author{Neeladri Maitra}
\address{Neeladri Maitra, Department of Mathematics, University of Illinois at Urbana-Champaign, 69.2.1 Computing Applications Building, 605 East Springfield Avenue, 61820, Champaign, Illinois, USA, nmaitra@illinois.edu}

\begin{abstract}
In this paper, we study the annealed ferromagnetic $q$-state Potts model on sparse rank-1 random graphs, where vertices are equipped with a vertex weight, and the probability of an edge is proportional to the product of the vertex weights. 
In an annealed system, we take the average on both numerator and denominator of the ratio defining the Boltzmann-Gibbs measure of the Potts model. 

We show that the thermodynamic limit of the pressure per particle exists for rather general vertex weights. In the infinite-variance weight case, we show that the critical temperature equals infinity. For finite-variance weights, we show that, under a rather general condition, the phase transition is {\em first order} for all $q\geq 3$. However, we cannot generally show that the discontinuity of the order parameter is {\em unique}. We prove this uniqueness under a reasonable condition that holds for various distributions, including uniform, gamma, log-normal, Rayleigh and Pareto distributions. Further, we show that the first-order phase transition {\em persists} even for some small positive external field.

In the rather relevant case of Pareto distributions with power-law exponent $\tau$, remarkably, the phase transition is first order when $\tau\geq 4$, but not necessarily when the weights have an infinite third-moment, i.e., when $\tau\in(3,4)$. More precisely, the phase transition is second order for $\tau\in (3,\tau(q)]$,  while it is first order when $\tau>\tau(q)$, where we give an explicit equation that $\tau(q)$ solves.

\end{abstract}

\maketitle

\setcounter{tocdepth}{1} 
\tableofcontents






\section{Introduction and main results}
\label{sec-intro-results}
\subsection{Introduction}
\label{sec-intro}

The ferromagnetic nearest-neighbor $q$-state Potts model \cite{potts1952some} has been extensively studied {\em on regular lattices} (such as $\Z^d$) as a prototype model of ferromagnetic phase transitions. The emergence of the phase transition in the thermodynamic limit is captured by the order parameter, i.e., the spontaneous magnetization, which changes from zero to a non-zero value as the inverse temperature $\beta$ crosses a critical value $\beta_c$.
The existence and {\em nature} of the phase transition for the Potts model depends sensitively on the dimension $d$ of the lattice and the number $q$ of colors. It is expected that the spontaneous magnetization of the Potts model with $d\ge 2$ and with $q\ge q_c(d)$ has a discontinuity (associated with latent heat) and therefore the model has a first-order phase transition.
Recent years have witnessed
spectacular signs of progress in the rigorous understanding of the phase transition of the Potts model, thanks to the use of the random cluster representation and percolation theory.
We refer to \cite{baxter2016exactly} and \cite{wu1982potts} for reviews of the Potts model on regular lattices, and to  \cite{duminil2017lectures} for a review of recent developments.
\smallskip

In this paper, we address the properties of the {\em ferromagnetic Potts model on sparse random graphs}. In this setting, the behavior of the Potts model is much less understood, even in the physics literature \cite{dorogovtsev2008critical}. In particular, the addition of the randomness of the spatial structures, which we model by inhomogeneous random graphs, poses additional new challenges compared to the usual set-up of the regular lattice. The first issue to deal with is the choice of how to average over the random graphs. Traditionally, one distinguishes between the {\em quenched} and the {\em annealed system}. In a quenched system, the graph is frozen when observed from the point of view of the spins, so that the quenched probability measure is the average of the random Boltzmann-Gibbs measure. In an annealed system, instead, the graph is treated on the same footing as the spins. Thus, in an annealed system, we take the average on both the numerator and denominator of the ratio defining the Boltzmann-Gibbs measure; see \eqref{annealed-measure} for a precise definition of the annealed probability measure. The quenched measure is appropriate to describe physical systems with sparse impurities (disordered ferromagnet with frozen couplings among the spins). The annealed measure, instead, is used when the Potts model is viewed as an opinion model for a population (here the couplings among individuals, representing their acquaintances, are rapidly changing). For spin models on random graphs, which are often mean-field systems, the folklore is that the quenched state is close to that of spin models on (unimodular) {\em random trees}, whereas the annealed state should be given by an {\em inhomogeneous mean-field theory}.
However, both in the quenched and annealed setting, rigorous results about the Potts model on random graphs are rather limited .
\smallskip

More precisely, a clear picture has emerged so far only for the Ising case $(q=2)$, while results are scarce for $q\ge 3$.
In the Ising case, for graph sequences converging locally to a random tree in the sense of Benjamini–Schramm \cite{benjamini2011recurrence}, the quenched pressure has been shown to coincide with the replica symmetric Bethe pressure, which is expressed in terms of solutions of 
a random recursion relation \cite{dembo2010ising,dembo2013factor, dommers2010ising}.
The annealed pressure instead is given
by the pressure of an inhomogeneous mean-field model,
whose precise form depends on the particular
random graph model. In particular, for rank-1 or generalized random graphs, the annealed pressure coincides with that of an inhomogeneous Curie-Weiss model \cite{giardina2015annealed, dommers2016ising}; for the random regular graph with degrees all equal \cite{can2019annealed}, or more generally for the configuration model with a prescribed degree sequence \cite{can2022annealed}, the annealed pressure has been obtained by a combinatorial argument, using a mapping to a random matching problem. Interestingly,
the critical inverse temperatures of quenched and annealed Ising ferromagnets do coincide when the random graph model
has a fixed number of edges, while the annealed critical inverse temperature is strictly lower than the quenched one when the number of edges in the graph is allowed to fluctuate. In other words, in the annealed state, the graph rearranges itself by increasing the number of edges (if it is allowed to do so), which in turn favors the alignment of spins. 
\smallskip

For the Potts model with $q\ge 3$, the coincidence between the quenched pressure and the replica symmetric formula has been proved only for
uniformly sparse graph sequences converging locally to the $r$-regular tree (first with the restriction of $r$ being an even number
\cite{dembo2014replica} and very recently, without any restriction, for a sequence of $r$-regular graphs  \cite{bencs2023random}). For other graph sequences,
the only available results are limited to the uniqueness regime \cite{dembo2013factor}.
As far as we know, the problem
of computing the quenched pressure of the Potts model on a locally-tree-like graph sequence remains among the most important open problems for statistical mechanics models on random graphs. 
As for the annealed pressure, this problem has not been the subject of rigorous mathematical studies so far.
\medskip

\paragraph{\bf Main results of this paper.}
In this paper, we start such an analysis by studying the ferromagnetic Potts model on rank-1 inhomogeneous random graphs. A natural expectation, that we prove rigorously, is that the model can be mapped to an inhomogeneous Curie-Weiss model with $q$ colors. As a consequence of this mapping, we can prove a number of surprising results, as described in more detail in Section \ref{sec-main-results}. In particular, we give an almost complete characterization of the thermodynamic limit and phase transition of the annealed $q$-states Potts on the rank-1 inhomogeneous random graph, one of the most popular random graph models for complex networks. We identify the first-order phase transition via the so-called ``equal area construction'' \cite{presutti2015equilibrium}. The first-order phase transition persists in the presence of a small enough magnetic field.
For the important case of Pareto weights in a zero external field, we show that 
a switch from first-order to second-order phase transition  takes place when a certain moment of the vertex weights becomes infinite. Let us now describe the precise model, as well as our main results.

\subsection{Annealed Potts model on rank-1 inhomogeneous random graphs}
\label{sec-model-intro}
In this section, we introduce the random graph that we will work on, as well as the annealed Potts model on it.
\medskip

\paragraph{\bf The rank-1 inhomogenenous or generalized  random graph.}
Let $G_n=(V_n, E_n)$ be a graph with vertex set $V_n= \{1,\ldots, n \}=:[n]$ and  random edge set $E_n$. Each vertex $i$ is equipped with a weight $w_i >0$. With $I_{i,j}$ denoting the indicator that the edge between vertex $i$ and $j$ is present, the variables $(I_{i,j})_{1\leq i<j\leq n}$ are independent, with
    \eqn{
    \prob(I_{ij}=1)=\frac{w_iw_j}{\ell_n+w_iw_j},
    \qquad\text{where}
    \qquad
    \ell_n=\sum_{i\in [n]}w_i
    }
denotes the total weight.



We assume that the positive weight sequence $(w_i)_{i\in [n] }$ satisfies a {\em regularity condition}. Let $F_n$ denote the empirical weight distribution, i.e.,
    \eqn{
    \label{Fn-def}
    F_n(x)=\frac{1}{n}\sum_{i\in [n]}\indic{w_i\leq x},
    }
and let $W_n$ have distribution function $F_n$. Equivalently, $W_n$ is the weight of a vertex chosen uniformly at random from the vertex set $[n]$. Our main assumption on the weights is as follows:

\begin{cond}[Weight regularity] 
\label{cond-WR-GRG}
There exists a random variable $W$ such that, as $n\rightarrow\infty$,
\begin{itemize}
\item[$(a)$] $W_n \stackrel{\mathcal D}{\longrightarrow} W$, 
where $\,\stackrel{\mathcal D}{\longrightarrow}\,$ denotes convergence in distribution;
\item[$(b)$] There exist $c>0$ and $\tau>2$ such that, for all $w\geq 1$,
    \eqn{
    \label{main-assumption-Fn}
    1-F_n(w)\leq \frac{c}{w^{\tau-1}},
    }
\end{itemize}
 Further, we assume that $\mathbb{E}[W]>0$. 
\end{cond}

\begin{remark}[Convergence of average weight]
\label{rem-average-weight}
{\rm
It is easy to see that Condition \ref{cond-WR-GRG}(b) implies that $\mathbb{E}[W_n] = \frac{1}{n}\sum_{i\in[n]} w_i\longrightarrow \mathbb{E}[W]< \infty$.}
\hfill \ensymboldefinition
\end{remark}
\medskip

\paragraph{\bf The annealed Potts model on the rank-1 random graph.}
To each $i\in [n]$, we assign a Potts spin $\sigma_i \in \{1,2,\ldots, q \}=:[q]$, with $q\in \N$. The Hamiltonian on a finite graph with $n$ vertices is 
\be
\label{eq:ham-potts}
H_n(\bs)=\beta \sum_{1\le i<j\le n}  I_{i,j}  \indic{\s_i=\s_j} + \sum_{i\in [n]} \sum_{k\in [q]} B_k \indic{\s_i=k}.
\ee
Here $\bs= (\sigma_i)_{i\in[n]}$ denotes the collection of all spin variables, $\beta\ge 0$ is the inverse temperature, and $\bB=(B_k)_{k \in[q]}$ is the vector of external magnetic fields, where $B_k \ge 0$ for all $k \in  [q]$.

Denoting the expectation with respect to the law of the random graph by $\E[\cdot]$, the finite graph {\em annealed} partition function is $\mathbb{E}[Z_n(\beta, \bB)]$, where the random partition function is given by
\be
    Z_n(\beta, \bB)=\sum_{\bs \in \omqn} \exp \left [ \beta \sum_{i<j}  I_{i,j} \indic{\s_i=\s_j} + \sum_{i\in [n]} \sum_{k\in [q]} B_k \indic{\s_i=k}   \right ].
\ee
This allows us to define the annealed probability measure of the Potts model on a graph $G_n=(V_n,E_n)$ as
\be
\label{annealed-measure}
\mu_{n;\beta,\bB}(\bs) = \frac{\E[\exp [H_n(\bs)]]}{{\E [ Z_n(\beta, \bB)]}}.
\ee
To study the annealed probability measure in the thermodynamic limit, it is useful to introduce the finite-volume annealed pressure per particle
    \eqn{
    \label{finite-pressure}
    \varphi_n(\beta,\bB) =\frac{1}{n}\log{\expec[Z_n(\beta,\bB)]}.
    }

\subsection{Main results}
\label{sec-main-results}
In this section, we state our main results. This section is organised as follows. We start by discussing the thermodynamic limit of the pressure per particle, in which a variational problem appears. We then consider the solution of this variational problem, and its relation to the phase transition. After this, we state a central condition in this paper, which we call the {\em zero-crossing condition}. It allows us to give a detailed description of the first-order phase transition for $q\geq 3$. We continue by giving examples of weight distributions for which this condition holds, and one example where it turns out to be false. We then zoom in into the special case of Pareto weights, a case that has attracted substantial attention in the physics literature. Finally, we state a result on the existence of a first-order phase transition for general weight distributions.
\medskip

\paragraph{\bf Thermodynamic limit of the pressure per particle.}
Before stating our results, we introduce the notation
\be
\label{betaprime}
\beta'=\e^\beta-1,
\ee
and we define the function
    \eqn{
    \label{F-pressure}
    F_{\beta,\bB}( \boldsymbol{y}) = \expec \left [ \log \left ( \e^{\frac{\beta'}{\expec[W]} y_1   W  +B_1}+\cdots + \e^{ \frac{\beta'}{\expec[W]} y_q  W +B_q} \right )\right ],
    }
where  $\boldsymbol{y}=(y_1,\dots,y_q) \in \mathbb{R}^q$.
The following theorem shows that the infinite volume limit of the {\em annealed pressure} per particle exists: 

\begin{theorem}[Existence of annealed pressure]
\label{thm-press-particle}
Fix $\beta>0$ and $q\in \mathbb{N}$. Consider the annealed Potts model on the generalized random graph with weights satisfying Condition \ref{cond-WR-GRG}. Then,
    \eqn{
    \lim_{n\rightarrow \infty} \varphi_n(\beta,\bB)
    =: \varphi(\beta,\bB),
    }
where
    \eqn{
    \label{eq:laplace}
    \varphi(\beta,\bB)
    =\sup_{\boldsymbol{y} \in \mathbb{R}^q} \Big[F_{\beta,\bB}( \boldsymbol{y} )-\frac{\beta'}{2 \mathbb{E}[W]} \sum_{k=1}^q y_k^2 \Big].
    }
\end{theorem}

Theorem \ref{thm-press-particle} writes the pressure per particle as a variational problem. Any maximum point $\by^\star(\beta,\bB)=(y_1^\star(\beta,\bB),\ldots , y_q^\star(\beta,\bB))$ satisfies the {\em stationarity condition} 
\be\label{eq:stat0}
 \expec \left [ \frac{ \e^{y_k  {\frac{\beta'}{\expec[W]}}  W +B_k} W }{ \e^{y_1  {\frac{\beta'}{\expec[W]}}  W+B_1}+\cdots + \e^{y_q  {\frac{\beta'}{ \expec[W]}}  W +B_q}}\right ] = y_k,\qquad k\in[q].
\ee

\begin{remark}[Properties of optimizers]\label{rem-optim-positive}{\rm 
By the stationarity condition and since $\expec[W] > 0$, we have that  $y_k^\star(\beta,\bB)> 0$ for all $k\in [q]$. Moreover, $\sum_{k\in[q]}y_k^\star(\beta,\bB) = \expec[W]$.}\hfill\ensymboldefinition
\end{remark}

We next investigate the structure of the optimizers in more detail.
\medskip

\paragraph{\bf The structure of the optimizers.}
From now on, unless otherwise stated, we fix $\bB  
=(B, 0,\ldots,0)$.
The next theorem shows that 
the optimizer of \eqref{eq:laplace} 
has the form $\boldsymbol{y}=(y_k(s))_{k \in [q]}$, where
    \be\label{eq:soln_form}
    y_k (s) = 
    \begin{cases}
    \frac{\expec[W]}{q}(1+(q-1)s)  & \text{if } k =1, \\
    &\\
    \frac{\expec[W]}{q} (1-s) & \text{otherwise,}
    \end{cases}
    \ee
for a suitable function $s=s(\beta, B)$ taking values in $[0,1]$. 
Of course, it requires a proof that the optimal solution really is of this form.
We sometimes slightly abuse notation
and write $y_k^\star(\beta,B)$ instead of $y_k^\star(s(\beta,B))$. 


\invisible{\begin{theorem}[Annealed pressure with zero magnetic field]
\label{thm-soln_form}
Let the external magnetic fields satisfy $\bB=\boldsymbol{0}$, and fix an inverse temperature $\beta>0$. Then any permutation of the vector $\boldsymbol{y}(s)$ defined above in \eqref{eq:soln_form} solves the optimization problem \eqref{eq:laplace}, for a suitable function $s=s(\beta,0) \in [0,1]$. Consequently, 
    \be
    \label{var-problem-psi}
    \varphi(\beta,0) = \sup_{s\in[0,1]} p_{\beta}(s),
    \ee
where
\eqan{\label{eq:pressure_evaluated_at_soln-a}
p_{\beta}(s) &=\mathbb{E}\left[ \log\left(\e^{\beta'W s} + q-1 \right) \right]- \frac{\beta' \expec[W]} {2 q}\left[ (q-1)s^2 
+2s-1\right].
}
The optimal solution $s(\beta,0)$ of \eqref{var-problem-psi} satisfies the {\em stationarity condition}
    \eqn{
    \label{FOPT-1-B1}
    \mathbb{E}\left[ \frac{W}{\expec[W]} \frac{\e^{\beta' W s}-1}{\e^{\beta' Ws}+q-1}\right] = s.
    }
\end{theorem}

We next extend the above theorem to positive
external field:}
\begin{theorem}[Annealed pressure]
\label{thm-soln_form-positive}
Fix an inverse temperature $\beta>0$ and let the external magnetic fields satisfy $\bB=(B,0, \ldots, 0)$ for some $B\geq 0$.
\begin{enumerate} 
\item[(a)] If $B=0$ then any optimizer $\by^\star(\beta,0)$ of the optimization problem \eqref{eq:laplace} is a permutation of the vector $\boldsymbol{y}(s)$ defined in \eqref{eq:soln_form} for a suitable function $s=s(\beta,0) \in [0,1]$;
\item[(b)] If $B>0$ then the optimizer of \eqref{eq:laplace} is given by $\boldsymbol{y}(s)$ in \eqref{eq:soln_form} for a suitable function $s=s(\beta,B) \in [0,1]$.
 
\end{enumerate}
Consequently, 
    \be
    \label{var-problem-psi-B>0}
    \varphi(\beta,B) = \sup_{s\in[0,1]} p_{\beta,B}(s),
    \ee
where
\eqan{\label{eq:pressure_evaluated_at_soln-a-B>0}
p_{\beta,B}(s) &=\mathbb{E}\left[ \log\left(\e^{\beta'W s+B} + q-1 \right) \right]- \frac{\beta' \expec[W]} {2 q}\left[ (q-1)s^2 
+2s-1\right].
}
Any optimal solution $s(\beta,B)$ of \eqref{var-problem-psi-B>0} satisfies the {\em stationarity condition}
    \eqn{
    \label{FOPT-1-B1-B>0}
    \mathbb{E}\left[ \frac{W}{\expec[W]} \frac{\e^{\beta' W s+B}-1}{\e^{\beta' Ws+B}+q-1}\right] = s.
    }
\end{theorem}


\invisible{Theorem \ref{thm-soln_form-positive} allows us to define the {\em critical value} $\beta_c$ as
    \eqn{
    \label{def-beta-c}
    \beta_c=\inf \{\beta>0\colon \text{ any optimizer $s(\beta, 0)$ of \eqref{var-problem-psi-B>0} satisfies } s(\beta, 0)>0\}.
    }
\begin{remark}[Uniqueness critical value]
\label{rem-critical-value}
{\rm
We have not yet shown that $s(\beta, 0)$ is unique, nor that $s(\beta_1, 0)>0$ implies that $s(\beta_2, 0)>0$ for all $\beta_2> \beta_1$. We will return to this uniqueness question later on.}
\hfill \ensymboldefinition
\end{remark}}

\invisible{Indeed, if $s(\beta)>0$ then, in the thermodynamic limit, the {\em instantaneous} fraction of spins in state $1$ defined as
\be
\label{x1-beta}
x_1(\beta)=\lim_{B\searrow 0} x_1(\beta,B)
\ee
satisfies $x_1(\beta)>1/q$, 
meaning that the system is spontaneously magnetized. 
}

We first note that the system is always magnetized when the second moment of $W$ is infinite:

\begin{theorem}[Spontaneous magnetization for infinite-variance weights]
\label{thm-spont-magnetization}
Assume that $B=0$ and $\expec[W^2]=\infty$ and fix $\beta>0$. Then any optimizer $s(\beta, 0)$ of \eqref{var-problem-psi-B>0} satisfies $s(\beta, 0)>0$.
\end{theorem}

\invisible{We next investigate the critical value for our annealed Potts model:

\begin{theorem}[Critical value]
\label{thm-crit-value}
 Then, there exists 
Let 
Then, the critical value of the annealed Potts model equals $\beta_c$, and
\end{theorem}

\RvdH{Do we know that $\beta_c>0$?}
\RvdH{Are we implicitly assuming that $\expec[W^2]<\infty$ in Theorem \ref{thm-crit-value}?}}

Theorem \ref{thm-spont-magnetization} shows that there is no phase transition when $\expec[W^2]=\infty$. We next move to the setting where $\expec[W^2]<\infty$, in which case we need to make an additional assumption to identify the nature of the phase transition, as we explain next.
\medskip

\paragraph{\bf Condition implying a first-order phase transition:  zero-crossing.}

To investigate the nature of the phase transition when $\expec[W^2]<\infty$, we 
introduce, for $t\ge 0$, 
    \eqn{
    \label{def-Fcal}
    \mathscr{F}_B(t)=\expec \left[\frac{W}{\expec[W]}\frac{\e^{tW+B}-1}{\e^{tW+B}+q-1} \right].
    }
To see that this function is relevant, we note that
    \eqn{
    \label{stationary-F-def}
    \frac{d}{ds}p_{\beta,B}(s)=\frac{q-1}{q} \beta' \expec[W]\Big[\mathscr{F}_B(s\beta')-s\Big],
    }
so that the stationarity equation \eqref{FOPT-1-B1-B>0} translates into
    $
    \mathscr{F}_B(s\beta')-s=0,
    $
or, when denoting $t=s \beta'$,
    \eqn{
    \label{stationary-condition-F}
    \mathscr{F}_B(t)=t/\beta'.
    }
In the formulation in \eqref{stationary-F-def}, the stationarity condition is translated into a fixed-point equation for $\mathscr{F}_B$, which does not depend on $\beta'$, but $\beta'$ appears as the inverse slope of the line on the r.h.s.\ of \eqref{stationary-condition-F}, which turns out to be quite convenient.

We can compute
    \be\label{eq:second_derivative_gen}
    \frac{d^2}{dt^2}\mathscr{F}_B(t)
    =q\e^{-B} \expec\left[\frac{W^3}{\expec[W]}\frac{\e^{tW}(\e^{-B}(q-1)-\e^{tW})}{(\e^{tW}+\e^{-B}(q-1))^3}\right],
    \ee
where we observe that, apart from a trivial factor $q\e^{-B}$, this is the same function as for $B=0$, but with $q-1$ replaced by $\e^{-B}(q-1)$. 
\invisible{\color{red}
    \eqan{
    \label{def-Fcal-again}
    \mathscr{F}_B(t)&=\expec \left[\frac{W}{\expec[W]}\frac{\e^{tW+B}-1}{\e^{tW+B}+q-1} \right]\nonumber\\
    &=1-q\expec \left[\frac{W}{\expec[W]}\frac{1}{\e^{tW+B}+q-1} \right]\nonumber\\
    &=1-q\e^{-B}\expec \left[\frac{W}{\expec[W]}\frac{1}{\e^{tW}+\e^{-B}(q-1)} \right].
    }
Thus, even $\mathscr{F}_B(t)$ depends on $B$ only through $\e^{-B}(q-1)$.
\color{black}}
Note also that the function $t\mapsto\mathscr{F}_B(t)$ is concave  when $q\ge 3$ and $B >\log(q-1)$. 
\smallskip

The following condition states a fundamental property of $\mathscr{F}_B$ that will play a central role throughout this paper:
\invisible{
\begin{cond}[{Steep} unique zero-crossing condition]
\label{cond-zer-cross}
We say that $\mathscr{F}_B$ satisfies the {\em steep unique zero-crossing condition} if
\begin{itemize}
\item[(a)] the {\em unique zero-crossing condition} holds, in that the function $t\mapsto \frac{d^2}{dt^2}\mathscr{F}_B(t)$ is first positive and then negative for $t\in [0,\infty)$;
\item[(b)] the {\em steepness condition holds}, in that the unique inflection point $t_*(B)>0$  of $t\mapsto \mathscr{F}_B(t)$ satisfies
\be
\label{eq:steep}
\frac{d}{dt}\mathscr{F}_B(t_*(B)) > \frac{\mathscr{F}_B(t_*(B))}{t_*(B)}.
\ee
\end{itemize}
\color{black}
\end{cond}}

\begin{cond}[Crossing conditions at the inflection points]
\label{cond-zer-cross}
We say that $t\mapsto\mathscr{F}_B(t)$ satisfies
\begin{itemize}
\item[(a)] the {\em unique zero-crossing condition} if the function $t\mapsto \frac{d^2}{dt^2}\mathscr{F}_B(t)$ is first positive and then negative for $t\in [0,\infty)$, i.e., there is a unique inflection point $t_*(B)>0$  of $t\mapsto \mathscr{F}_B(t)$;
\item[(b)] the  {\em steep zero-crossing condition} at an inflection point $t_*(B)>0$ of 
$t\mapsto \mathscr{F}_B(t)$
if
\be
\label{eq:steep}
\frac{d}{dt}\mathscr{F}_B(t_*(B)) > \frac{\mathscr{F}_B(t_*(B))}{t_*(B)}.
\ee
\end{itemize}
If both conditions (a) and (b) hold, we say that $t\mapsto \mathscr{F}_B(t)$ satisfies the
{\em zero-crossing condition}.
\end{cond}

\begin{remark}[Steepness follows from uniqueness for small $B$]
\label{steep-remark}
{\rm We remark that if $B=0$ then the steepness Condition \ref{cond-zer-cross}(b) follows from the unique-zero crossing Condition \ref{cond-zer-cross} (a).
Furthermore, a continuity argument guarantees that there exists $\delta>0$ such that  Condition \ref{cond-zer-cross}(b) for $0<B<\delta$ follows from Condition \ref{cond-zer-cross}(b) with $B=0$.}\hfill\ensymboldefinition
\end{remark}

\invisible{From now on we fix $\bB = (B, 0,\ldots,0)$ with $B\ge 0$ and denote $\boldsymbol{y}(\beta,B)=(y_k (\beta,B))_{1\leq k \leq q}$ the optimizer 
of the variational problem in \eqref{eq:laplace}
with such a choice of the external magnetic fields. 
As we shall prove below, $y_k (\beta,B)$ 
has the form
    \be\label{eq:soln_form}
    y_k (\beta,B) = 
    \begin{cases}
    \frac{\expec[W]}{q}(1+(q-1)s(\beta,B))  & \text{if } k =1, \\
    &\\
    \frac{\expec[W]}{q} (1-s(\beta,B)) & \text{otherwise,}
    \end{cases}
    \ee
where the function $s(\beta,B)$ takes values in $[0,1]$ and is identified in Theorem  \ref{thm-soln_form}.

We also define the average proportion (in the thermodynamic limit) of spins that are
in the spin state $k\in [q]$ as 
\be
x_k(\beta,B) =\lim_{n\to\infty}
\sum_{\bs\in [q]^n}  \left( \frac{1}{n}\sum_{i=1}^n \indic{\sigma_i=k} \right) \mu_{n;\beta,B}(\bs).
\ee
where $\mu_{n;\beta,B}$ is the annealed measure \eqref{annealed-measure} with
inverse temperature $\beta$ and external field $B \ge 0$
acting on the first spin component.

\RvdH{Should we restrict Theorem \ref{thm-soln_form} to $B=0$? If so, then how do we break the symmetry? In the proof, we order the $y_i$ by their size, making $y_1$ the largest. How do we define the phase transition when we only use $B=0$?}
\NM{Neeladri: I think a \emph{uniform} $y_i$ becomes big (instead of $y_1$). At least that's what seems to happen in the homogeneous case in \cite{Ellis_Wang_90}.}
\RvdH{Yes, that is suggestive. But can we prove it?}
\begin{theorem}[Pressure with single positive external field]
\label{thm-soln_form}
Let the external magnetic fields satisfy $\mathbf{B}=(B, 0, \ldots, 0)$ with $B\ge 0$, and fix an inverse temperature $\beta>0$. Then the solution to the optimization problem \eqref{eq:laplace} has the form \eqref{eq:soln_form}, for some $s=s(\beta,B) \in [0,1]$. Consequently, 
    \be
    \label{var-problem-psi}
    \varphi(\beta,B) = \sup_{s\in[0,1]} p_{\beta, B}(s),
    \ee
where
\eqan{\label{eq:pressure_evaluated_at_soln-a}
p_{\beta, B}(s) &=\mathbb{E}\left[ \log\left(\e^{\beta'W s+B} + q-1 \right) \right]- \frac{\beta' \expec[W]} {2 q}\left[ (q-1)s^2 
+2s-1\right].
}
The optimal solution $s(\beta,B)$ of \eqref{var-problem-psi} satisfies the {\em stationarity condition}
    \eqn{
    \label{FOPT-1-B1}
    \mathbb{E}\left[ \frac{W}{\expec[W]} \frac{\e^{\beta' W s+B}-1}{\e^{\beta' Ws+B}+q-1}\right] = s.
    }
The average proportion of 1-spins equals 
\eqn{
\label{x1-def}
    x_1(\beta,B)=
    \mathbb{E}\left[\frac{\e^{\beta' W s(\beta,B) }}{\e^{\beta' W s(\beta,B)} + q-1 }\right]
}
while the proportion of $k$-spins for $k\ge2$ equals
    \eqn{
    \label{xk-beta}
    x_k(\beta,B)=     \mathbb{E}\left[\frac{1}{\e^{\beta' W s(\beta,B)} + q-1 }\right]
    }
Further, $B\mapsto s(\beta, B)$ is decreasing, so that 
    \eqn{
    s(\beta)=\lim_{B\searrow 0} s(\beta,B)
    }
is well defined.
\end{theorem}
}

\invisible{We start by drawing the analogy to the Ising model:

\begin{remark}[The Ising critical value]
\label{rem-Ising-case}
For $q=2$,
    \be
    \mathbb{E}\left[ \frac{W}{\expec[W]} \frac{\e^{\beta' W s}-1}{\e^{\beta' W s}+1}\right] = s.
    \ee
and linearizing around $s=0$ we find
\be
\beta_c' =  2 \frac{\expec[W]}{\expec[W^2]}
\ee 
\end{remark}

The situation should be different when $q \ge 3$ as there should be a jump in the 
optimal solution $y^{\star}(s)$, corresponding to a first order phase transition.
The value of the jump is located by imposing the condition $p(s) = p(0)$
which reads
\be
\mathbb{E}\left[\log\left(\frac{\e^{\beta' W s}+q-1}{q}\right) \frac{1}{\expec[W]}\right] = \frac{\beta'}{q} \left(s + \frac12 (q-1) s^2\right).
\ee 
Thus, the task is to solve the system of equations \eqref{FOPT-1-B1-B>0}--\eqref{crit-cond-B>0}.

\begin{remark}
In the homogeneous set up $W=1$ one recover the equations
\be
\frac{e^{\beta'  s}-1}{e^{\beta'  s}+q-1}= s
\ee
\be
\log\left(\frac{e^{\beta'  s}+q-1}{q}\right) = \frac{\beta'}{q} \left(s + \frac12 (q-1) s^2\right)
\ee 
whose solutions are
\be
s(\beta_c)(q) = \frac{q-2}{q-1}
\ee
\be
\beta'_c(q) = 2 \frac{q-1}{q-2} \log (q-1)
\ee
Note that $\lim_{q\to 2} s(\beta_c)(q) = 0$ and $\lim_{q\to 2} \beta'_c(q) = 2$.
\end{remark}}

\invisible{Indeed, assume that 
\be
\label{eq:steep-2}
\frac{d}{dt}\mathscr{F}_0(t_*(0)) > \frac{\mathscr{F}_0(t_*(0))}{t_*(0)}.
\ee
Then, for $0<B<\delta$,
\be
\label{eq:steep-3}
\frac{d}{dt}\mathscr{F}_B(t_*(B)) > \frac{\mathscr{F}_B(t_*(B))}{t_*(B)},
\ee
because $\mathscr{F}_B(t)$ and $\frac{d}{dt}\mathscr{F}_B(t)$ are continuous functions both in $t$ and $B$ and furthermore, $t_*(B)$ is also continuous in $B$ and different from 0.
These properties of $t_*(B)$ follow from the implicit function theorem.
 Indeed, $t_*(0)$ solves $\calF_0''(t_*(0))=0$ and we know that $t_*(0)>0$. Therefore $t_*(B)$ satisfying  $\calF_B''(t_*(B))=0$ is also positive
for small $B$. The application of the implicit function theorem is possible
because $\frac{d^3\calF_0}{dt^3} (t_*(0)) \neq 0$ as
\be
    \frac{d^3}{dt^3}\mathscr{F}_0(t)=q\expec\left[\frac{W^4}{\expec[W]}\frac{\e^{tW}\Big(\e^{2tW}-4\e^{tW}(q-1)+(q-1)^2\Big)}{(\e^{tW+B}+q-1)^4}\right] >0.
\ee
\RvdH{Should we not move this to later? Bit long here. And why is the above integral positive?}
}

\paragraph{\bf The nature of the phase transition under the zero-crossing condition.} 

We have the following two main results about the classification of stationary points and the nature of the phase transition:

\begin{theorem}[Stationary solutions]
\label{thm-nature-pt-a}
Let $\expec[W^2]<\infty$ and fix $q\geq 3$. Assume that the external field is such that $0\le B < \log(q-1)$. The solutions to the stationarity condition in \eqref{FOPT-1-B1-B>0} satisfy the following:
\begin{itemize}
    \item[(a)] Fix $B=0$, and assume that  Condition \ref{cond-zer-cross}(a) holds. Then, there are $\beta_0\in (0,\infty)$ and $\beta_1\in (\beta_0,\infty)$ such that only the trivial solution $s_1(\beta,0)=0$ exists for $\beta\in (0,\beta_0)$, there are three stationary solutions $s_1(\beta,0)=0<s_2(\beta,0)<s_3(\beta,0)$ for $\beta\in (\beta_0,\beta_1),$ and two stationary solutions $s_1(\beta,0)=0<s_3(\beta,0)$ for $\beta>\beta_1.$
    
    \item[(b)] Fix $B\in(0,\log(q-1))$, and assume that  Condition \ref{cond-zer-cross} holds. Then, there are $\beta_0(B)\in (0,\infty)$ and $\beta_1(B)\in (\beta_0(B),\infty)$ such that there is a unique solution $s_1(\beta,B)>0$ for $\beta\in (0,\beta_0(B))$, there are three stationary solutions $0<s_1(\beta,B)<s_2(\beta,B)<s_3(\beta,B)$ for $\beta\in (\beta_0(B),\beta_1(B)),$ and a unique stationary solution $s_3(\beta,B)>0$ for $\beta>\beta_1(B).$
\end{itemize}
\end{theorem}

\begin{remark}[Conditions on $B$]
\label{rem-cond-B}
{\rm The restriction to $0\le B < \log(q-1)$ can be understood from \eqref{def-Fcal}, which writes $\mathscr{F}_B(t)=\expec[h_B(t;W)]$, where $h_B(t;W)=1-\frac{W}{\expec[W]}\frac{q\e^{-B}}{\e^{tW}+\e^{-B}(q-1)}$. Note that $t\mapsto h_B(t;w)$ is concave for $B\geq \log(q-1)$ and all $w\geq 0$, so that also $t\mapsto \mathscr{F}_B(t)$ is. As a result, the nature of solutions of the stationarity condition is rather different for $B\geq \log(q-1)$, in that there always is a {\em unique} stationary solution. When $0\le B < \log(q-1)$ the unique zero-crossing Condition \ref{cond-zer-cross} is needed to guarantee that there is a regime where there are three stationary solutions. When there does not exist a regime where there are three stationary solutions, there will not be a phase transition, as the optimizer depends continuously on $\beta$. As such, both restrictions on $B$ (i.e., $0\le B < \log(q-1)$ and the Condition \ref{cond-zer-cross} are necessary and sufficient for there to be three stationary solutions.}\hfill\ensymboldefinition
\end{remark}

We next investigate the consequences of the structure of the solutions to the stationarity condition on the nature of the phase transition. 
%

\begin{theorem}[Nature of the phase transition]
\label{thm-nature-pt-b}
Let $\expec[W^2]<\infty$ and fix $q\geq 3$. Assume that the external field is such that $0\le B < \log(q-1)$ and assume  Condition \ref{cond-zer-cross}. 
Then the optimal value $s^\star(\beta,B)$,  which is  the stationary solution that maximizes $p_{\beta,B}(s)$ in \eqref{eq:pressure_evaluated_at_soln-a-B>0}, is {\em unique}, except for one $\beta=\beta_c(B)$. More precisely,
there exists a {\em critical inverse temperature} $\beta_c(B)\in(0,\infty)$ such that
\be
s^\star(\beta,B) = 
\begin{cases}
s_1(\beta,B) & \text{for } \beta<\beta_c(B)\\
s_3(\beta,B) & \text{for } \beta>\beta_c(B),
\end{cases}
\ee
where $s_1(\beta,B)$ and $s_3(\beta,B)$ are
the stationary solutions identified in Theorem \ref{thm-nature-pt-a}.
Further, $(\beta_c(B),s_1(\beta_c(B),B))$ and $(\beta_c(B),s_3(\beta_c(B),B))$
are the unique solutions to the system of equations given by the stationarity condition \eqref{FOPT-1-B1-B>0}
and the {\em criticality condition}
    \eqn{
    \label{crit-cond-B>0}
    p_{\beta_c(B),B}( s_1(\beta_c(B),B))=p_{\beta_c(B),B}( s_3(\beta_c(B),B)).
    }
Finally, $\beta\mapsto s^\star(\beta,B)$ is increasing and $0\le s_1(\beta_c(B),B) < s_3(\beta_c(B),B) < \infty$, i.e., the phase transition is \rm{first order}.
\end{theorem}

\medskip

\paragraph{\bf The order parameter}
\invisible{and we let $\mathscr{Y}_{\beta,B}$ be the set of all permutations of the vector $\by^\star(\beta,B)$ in \eqref{eq:soln_form}. By symmetry, each of the solutions in $\mathscr{Y}_\beta$ is an optimizer of \eqref{eq:laplace}. Further, this optimizer is {\em unique} when $s(\beta)=0$, while it consists of $q$ choices when  $s(\beta)>0$. 
\smallskip}
We now discuss the order parameter associated to the phase transition. We define the fraction of spins that are
of color $k\in [q]$ as 
\be
\label{X-k-number-def}
X_k=X_k(\bs) =\frac{1}{n}\sum_{i=1}^n \indic{\sigma_i=k}.
\ee
\invisible{
where $\mu_{n,\beta,0}$ is the annealed measure \eqref{annealed-measure} with
inverse temperature $\beta$ and external field $B=0$
acting on the first spin component. Further, define the vector $\bx^\star(\beta)$ by
    \eqn{
    \label{x1-def}
    x_1(\beta)=
    \mathbb{E}
    \left[\frac{\e^{\beta' W s(\beta) }}{\e^{\beta' W s(\beta)} + q-1 }\right],
    }
}
The following theorem shows how the limiting behavior of color proportions can be extracted from the optimizer
$\by^\star(\beta,B)$ of 
\eqref{eq:laplace}.



\begin{theorem}[Limiting color proportions]
\label{thm-proportions}
Let $\expec[W^2]<\infty$ and fix $q\geq 3$.
Fix $0 < B < \log(q-1)$ and assume Condition \ref{cond-zer-cross}. 
\invisible{Then, for all $t\geq 0$,
    \eqn{
    \label{varphi-beta-t}
    \lim_{n\rightarrow \infty}\frac{1}{n}\log\expec_{\mu_{n;\beta,B}}[\e^{tY_1}]
    =\psi_{\beta,B}(t).
    }
    }
Then, for $\beta\neq \beta_c(B)$,
    \eqn{
    \label{der-x1-equation}
    \lim_{n\rightarrow \infty}\expec_{\mu_{n;\beta,B}}\big[X_1\big] = \frac{\partial \varphi(\beta,B)}{\partial B}\equiv x^\star_1(\beta,B),
    } 
    where
    \be
    \label{x1-equation}
    x^\star_1(\beta,B) = \expec \Bigg[ \frac
    { \e^{y^\star_1(\beta,B)  {\frac{\beta'}{\expec[W]}}  W +B}  }
    { \e^{y^\star_1(\beta,B)  {\frac{\beta'}{\expec[W]}}  W+B}
     +\e^{y^\star_2(\beta,B)  {\frac{\beta'}{\expec[W]}}  W  }
     +\cdots 
     +\e^{y^\star_q(\beta,B)  {\frac{\beta'}{\expec[W]}}  W }
     }\Bigg].
\ee
\end{theorem}
Clearly, due to the symmetry of the other colors,  for $k=2,\ldots,q$,
\be
\lim_{n\rightarrow \infty}\expec_{\mu_{n;\beta,B}}\big[X_k\big] = \frac{1- x^\star_1(\beta,B)}{q-1}.
\ee
Therefore, when the assumptions of Theorem \ref{thm-proportions} hold, the fraction of spins having color 1 
has a jump at $\beta= \beta_c(B)$. This justifies calling $x^\star_1(\beta,B)$ the order parameter when
$0 < B < \log(q-1)$. 

When $B=0$ we know that
any permutation of \eqref{eq:soln_form} is an optimizer
and therefore, for any $k\in[q]$ and any $\beta\in[0,\infty)$, 
\be
\lim_{n\rightarrow \infty}\expec_{\mu_{n;\beta,0}}\big[X_k\big] = \frac{1}{q}.
\ee
In this case, the order parameter can be obtained by considering the $\lim_{B\searrow 0}x^\star_1(\beta,B)$. Namely, by \eqref{x1-equation},  when $\beta<\beta_c(0)$,
\be
\lim_{B\searrow 0}x^\star_1(\beta,B) = 
1/q,
\ee
while, when $\beta>\beta_c(0)$,
\begin{eqnarray}
   \lim_{B\searrow 0}x^\star_1(\beta,B) 
   &=& 
   \expec \Bigg[ 
   \frac{\e^{\frac{\beta' W}{q}(1+(q-1)s^\star(\beta,0))}}
    {\e^{\frac{\beta' W}{q}(1+(q-1)s^\star(\beta,0))} + (q-1)\e^{\frac{\beta' W}{q}(1-s^\star(\beta,0))}}
    \Bigg] \nonumber\\
&=& \expec \Big[ 
   \frac{\e^{\beta' W s^\star(\beta,0)}}
    {\e^{\beta' W s^\star(\beta,0)}+q-1}
    \Big] > \frac{1}{q}.
\end{eqnarray}

\begin{remark}[Weighted color proportions]
{\em As will become  clear from the proof of Theorem \ref{thm-proportions}, the optimizer $\by^\star(\beta,B)$
is the average of the weighted proportions of colors. More precisely,
defining
\be
\label{Y-k-number-def}
Y_k=Y_k(\bs) =\frac{1}{n \expec[W_n]}\sum_{i=1}^n w_i\indic{\sigma_i=k}.
\ee
we have
 \eqn{
    \label{psi-mGF}
    \lim_{n\rightarrow \infty}  \expec_{\mu_{n;\beta,B}}[Y_k]=\frac{y^\star_k(\beta,B)}{\expec[W]}.
    }
    }
    \hfill\ensymboldefinition 
\end{remark}

\medskip

\paragraph{\bf The critical value with zero field.}
We next specialize to the important case of $B=0$.
In this case, since $s_1(\beta,0)=0$, the criticality condition \eqref{crit-cond-B>0} takes on the form
    \eqn{
    \frac{1}{\expec[W]}\expec\left[\log\left(\frac{\e^{\beta' Ws}+q-1}{q}\right)\right]
    =\frac{\beta'}{q}\Big(s+\frac{1}{2}(q-1)s^2\Big).
    }
As a consequence, denoting the critical value $\beta_c(0)$ by $\beta_c$, we have the following characterization of the critical inverse temperature: 
\begin{theorem}[Critical value]
\label{thm-critical-value}
Let $\expec[W^2]<\infty$ and fix $q\geq 3$. Assume that $B=0$ and that  Condition \ref{cond-zer-cross}(a) holds. Let 
    \eqn{
    \label{K-function}
    \mathscr{K}(t)=
    \frac{1}{\expec[W]}\expec\left[\log\left(\frac{\e^{tW}+q-1}{q}\right)\right]
    -\frac{q-1}{2q} t\mathscr{F}_0(t)-\frac{t}{q}.
    }
There is a {\em unique} $t_c>0$ such that $\mathscr{K}(t_c)=0$. In terms of this $t_c$, the critical value is given by
    \eqn{
    \label{beta-s-crit}
        \beta_c=\log\left(1+\frac{t_c}{\mathscr{F}_0(t_c)}\right),
        }
        and, furthermore,
    \eqn{
    s(\beta_c,0)=\mathscr{F}_0(t_c).
    }
\end{theorem}

\begin{remark}
[Critical value for the Erd\H{o}s-R\'enyi random graph]
\label{remark-ER}
{\rm
The choice of weights $w_i= \frac{n\lambda}{n-\lambda}$ corresponds to the sparse Erd\H{o}s-R\'enyi random graph
where edges are present with equal probability equal to $\lambda/n$. In this case the distribution of $W$ is a Dirac measure centered 
at $\lambda$ and therefore the $\mathscr{K}$ function 
reads
   \eqn{
    \mathscr{K}^{\sss\rm{ER}}(t)=
    \frac{1}{\lambda}\log\left(\frac{\e^{t\lambda }+q-1}{q}\right)
    -\frac{(q+1)}{2q} t + \frac{1}{2} \frac{(q-1)t}{\e^{t\lambda}+q-1}.    
    }
The solution of 
$\mathscr{K}^{\sss\rm{ER}}(t_c)=0$ is 
$t_c= \frac{2}{\lambda}\log(q-1)$
and thus
    \be
    {\beta'}_c^{\sss\rm{ER}} = 
    \e^{\beta_c^{\sss\rm{ER}}}-1
    =\frac{2}{\lambda} \frac{q-1}{q-2} \log (q-1),
    \ee
    \be
    s({\beta}_c^{\sss\rm{ER}},0) = \frac{q-2}{q-1}.
    \ee
We retrieve the solution of the homogenous Curie-Weiss-Potts model \cite{wu1982potts}.}\hfill\ensymboldefinition
\end{remark}

\invisible{\noindent
{\color{red} {\bf Cristian:}
Can we improve  Theorem \ref{thm-proportions}
 now that we know the form of the optimizer?
If we define 
\eqn{
\varphi(\beta,B,u) = \lim_{n\to\infty} \frac{1}{n}\log\mathbb{E}\left[\sum_{\bs \in \omqn} \exp \left [ \beta \sum_{i<j}  I_{i,j} \indic{\s_i=\s_j} + B \sum_{i\in [n]}  \indic{\s_i=1} + \frac{u}{\mathbb{E}[W]} \sum_{i\in[n]} w_i\indic{\s_i=1}    \right ]\right]
}
then we have
   \eqn{
   \label{eq-pippo}
    \varphi(\beta,B,u)
    =\sup_{\boldsymbol{y} \in \mathbb{R}^q} \Big[{F}_{\beta,B}(u; \boldsymbol{y} )-\frac{\beta'}{2 \mathbb{E}[W]} \sum_{k=1}^q y_k^2 \Big]
    }
    with
        \eqn{
    {F}_{\beta,B}(u; \boldsymbol{y})  = \expec \left [ \log \left ( \e^{\frac{W}{\expec[W]} (\beta'y_1+u)    +B} 
    + \e^{ \frac{W}{\expec[W]} \beta' y_2} + \cdots + \e^{ \frac{W}{\expec[W]} \beta'y_q} \right )\right ].
    }
    If we call $y^*(\beta,B,u)$ the optimizer of \eqref{eq-pippo} we then have
    \eqn{
    \lim_{n\rightarrow \infty}\frac{1}{n}\expec_{\mu_{n;\beta,B}}[Y_1] = \frac{\partial \varphi(\beta,B,u)}{\partial u}|_{u=0} = y_1^*(\beta,B,0)=y^*_1(\beta,B)
    }
       where $y(\beta,B)$ is the solution of \eqref{eq:stat0}.
    Furthermore we have
\eqn{
    \lim_{n\rightarrow \infty}\frac{1}{n}\expec_{\mu_{n;\beta,B}}\left[\sum_{i}\indic{\s_i=1}\right] = \frac{\partial \varphi(\beta,B,u)}{\partial B}|_{u=0} = x_1^*(\beta,B,0)=x^*_1(\beta,B)
    } 
    where
    \be
    \label{x1-equation}
    x^*_1(\beta,B) = \expec \left [ \frac
    { \e^{y^*_1  {\frac{\beta'}{\expec[W]}}  W +B}  }
    { \e^{y^*_1  {\frac{\beta'}{\expec[W]}}  W+B}
     +\e^{y^*_2  {\frac{\beta'}{\expec[W]}}  W  }
     +\cdots 
     +\e^{y^*_q  {\frac{\beta'}{\expec[W]}}  W  }
     }\right ].
\ee
}

\RvdH{The above is correct, and follows from our proof. However, as far as I am aware, we have no way of proving that $x^*_1(\beta,0)>1/q$ precisely when $s(\beta, 0)>0.$ We we wish to add the unweighted spin proportions?}
\medskip

\noindent
{\color{red} {\bf Cristian:}
If $s(\beta,0)=0$ then from \eqref{x1-equation}
\be
   x^*_1(\beta,0) = \expec \left [ 
   \frac{\e^{\frac{\beta' W}{q}}}
    {q\e^{\frac{\beta' W}{q}}}
    \right ] = \frac{1}{q}
\ee
If $s(\beta,0)>0$ then from \eqref{x1-equation}
\begin{eqnarray}
   x^*_1(\beta,0) 
   &=& 
   \expec \left [ 
   \frac{\e^{\frac{\beta' W}{q}\big(1+(q-1)s(\beta,0)\big)}}
    {\e^{\frac{\beta' W}{q}\big(1+(q-1)s(\beta,0)\big)} + (q-1)\e^{\frac{\beta' W}{q}\big(1-s(\beta,0)\big)}}
    \right ] \nonumber\\
&=& \expec \left [ 
   \frac{\e^{\beta' W s(\beta,0)}}
    {\e^{\beta' W s(\beta,0)}+(q-1)}
    \right ] > \frac{1}{q}
\end{eqnarray}
}}


\paragraph{\bf When does the zero-crossing condition 
hold?}
While we cannot prove that Condition \ref{cond-zer-cross} holds in general, we can show that $t\mapsto \frac{d^2}{dt^2}\mathscr{F}_B(t)$ is first positive and then negative for $\e^{-B}(q-1)>1$ for a rather broad class of distributions. Therefore, for this class of distributions, Condition \ref{cond-zer-cross} holds for $B=0$ and for $B$ small enough (see Remark \ref{steep-remark}). As a consequence, the annealed Potts model on rank-1 inhomogeneous random graphs with these weight distributions have a first-order phase transition (see Theorem \ref{thm-nature-pt-b}). We write, for general $0\leq b<c\leq \infty$ and $p\in \mathbb{R}$, the probability density function as
	\eqn{
	\label{Guido(3)}
	f_{\sss W}(w)=C w^{p} \e^{-\varphi(w)} \indicwo{[b, c)}(w), 
	}
with $\varphi(w) \geq 0$ and $C \geq 0$ such that $\mathbb{E}[W]>0$ and for which \eqref{main-assumption-Fn} in Condition \ref{cond-WR-GRG} holds. 
\invisible{that there exists $D>0$ and $\tau>2$ such that
	\eqn{
	\label{main-assumption-Fn}
	1-\prob(W \leq w) \leq \frac{D}{w^{\tau-1}},\qquad w \geq 1.
	}
\RvdH{This actually follows from \eqref{main-assumption-Fn}}}

\invisible{where we call a random variable $W$ {\em monomial} when its density is given by
    \eqn{
    \label{def-monomial}
    f_W(w) = 
    \begin{cases}
        \frac{(p+1)w^p}{b^{p+1}-a^{p+1}}
        &\text{for }w\in[a,b],\\
        0 &\text{otherwise}.
    \end{cases}
    }
The class of monomial distributions in \eqref{def-monomial} is quite general. It includes uniform distributions, as well as power-law distributions for $b=\infty$ and $p<-3$.
\Nel{do we interpret the $b^{p+1}$ term above to be $0$ for $b=\infty$?}}

The following theorem gives conditions for which densities of the form \eqref{Guido(3)} satisfy  Condition \ref{cond-zer-cross}(a):

\begin{theorem}[Special weight cases]
\label{thm-density-Guido}
Fix $\e^{-B}(q-1)>1$. The function $t\mapsto \frac{d^2}{dt^2}\mathscr{F}_B(t)$ is first positive and then negative for densities of the form \eqref{Guido(3)} that satisfy \eqref{main-assumption-Fn}, and for which
\begin{itemize}
    \item[$\rhd$] $\varphi \in C^{2}([b, c])$ and $w\varphi^{\prime}(w)$ is non-negative and non-decreasing in $w \in[b, c]$ for $0<b<c<\infty$;
    \item[$\rhd$] $\varphi \in C^{2}([b, \infty))$ and  $w\varphi^{\prime}(w)$ is non-negative and non-decreasing in $w\in [b, \infty)$ with $\lim _{w \rightarrow \infty} w\varphi^{\prime}(w)=\infty$ for $b>0, c=\infty;$
    \item[$\rhd$] $p>-1$, $\varphi \in C^{2}([0, c])$ and $w\varphi^{\prime}(w)$ is non-negative and non-decreasing in $w \in[0, c]$ for $b=0, c<\infty$;
    \item[$\rhd$] $p>-1$, $\varphi \in C^{2}([b, c])$ and $w\varphi^{\prime}(w)$ is non-negative and non-decreasing in $w \in[0, \infty)$ with $\lim _{w \rightarrow \infty} w\varphi^{\prime}(w)=\infty$ for $b=0, c=\infty$.
\end{itemize}
\end{theorem}

\begin{remark}[Conditions]
\label{rem-conditions}
{\rm The conditions above can be interpreted as follows. When $b=0$, the density $f_{\sss W}$ needs to be integrable in 0, which explains why we need $p>-1$. When $c=\infty$, then we need the density to be integrable in infinity as well. For this, the condition $\lim _{w \rightarrow \infty} w\varphi^{\prime}(w)=\infty$ is highly relevant, but stronger than necessary. The precise form of the conditions is dictated by our proof.}\\
{\rm Note that the form of the density in \eqref{Guido(3)} is preserved by replacing $W$ by $W^\alpha$ for any $\alpha>0$. For example, after this transformation, $\varphi(w)$ is replaced by $\varphi(w^{1/\alpha})$, so that also $\lim_{w \rightarrow \infty} w\varphi^{\prime}(w)=\infty$ is preserved. Further, $p$ is replaced by $(p+1)/\alpha -1$, so that $p>-1$ implies that $(p+1)/\alpha -1>-1$.}
\hfill \ensymboldefinition
\end{remark}

\begin{remark}[Examples]
\label{rem-examples}
{\rm There are many examples for which the conditions in Theorem \ref{thm-density-Guido} hold. Indeed, the conditions hold for general uniform, gamma, Rayleigh, Weibull, log-normal, and Pareto random variables amongst others, as well as positive powers of them.}
\hfill \ensymboldefinition
\end{remark}

\invisible{We believe that the function $t\mapsto \frac{d^2}{dt^2}\mathscr{F}(t)$ is first positive and then negative for {\em all} distributions $W$, but lack a proof for this, except for the case where $W$ has a monomial distribution and when $W$ is constant. The latter corresponds to the Potts-Curie-Weiss model.
\Gib{This comment is obsolete. Now we know that statement is not true in general, but we have a large set of cases in which it works.}

\Nel{Should we write a new theorem at this point dedicated to the special Pareto case? We could first state the $\tau>4$ case, then the exotic $\tau \in (3,4)$ case.}}

We next consider cases where $W$ has a small compact support. The case where $W$ has a single atom, or a singleton set as support, corresponds to the annealed Erd\H{o}s-R\'enyi random graph, for which the function $t\mapsto \frac{d^2}{dt^2}\mathscr{F}_B(t)$ is indeed first positive and then negative. However, for random variables $W$ with small support, the picture becomes more subtle:


\begin{theorem}[Small support]
\label{thm-two-atoms}
Fix $q\geq 3$ and $0\leq B<\log(q-1)$. Let $W$ have a distribution that has support contained inside a compact set $[w',w]$, with $w'>0$ and 
\begin{equation}\label{eq:small_support_condition}
    \frac{w}{w'}\leq 1+\frac{\log{(2+\sqrt{3})}}{\log{(q-1)}-B}.
\end{equation}
Then, $t\mapsto \frac{d^2}{dt^2}\mathscr{F}_B(t)$ is first positive and then negative.
\end{theorem}
\begin{remark}[Counterexample with two atoms]
\label{counterex}
{\rm In Appendix \ref{appA}, we give an example with $B=0$ and $q=7$, and with a choice of $W$ having exactly two atoms that  are sufficiently far apart, and for which $t\mapsto \frac{d^2}{dt^2}\mathscr{F}_0(t)$ crosses zero more than once.}\hfill \ensymboldefinition
\end{remark}

\paragraph{\bf The Pareto case: smoothed phase transition.}
We finally consider the Pareto case with parameter $\tau>2$, for which
	\eqn{
	\label{Pareto}
	f_{\sss W}(w)=(\tau-1) w^{-\tau}\, \indicwo{[1, \infty)}(w).
	}
Such distributions have attracted tremendous attention, since they give rise to random graphs with power-law degrees. Theorem \ref{thm-spont-magnetization} shows that there is instantaneous order when $\tau\in (2,3]$, since then $\expec[W^2]=\infty$. Theorem \ref{thm-density-Guido} shows that the unique zero-crossing condition of $t\mapsto \frac{d^2}{dt^2}\mathscr{F}_B(t)$ holds for $\tau \ge 4$, which we state separately here:
\begin{theorem}[Pareto weights with $\tau \ge 4$]
\label{thm-Pareto-tau>4}
Fix $q\geq 3$ and $\e^{-B}(q-1)>1$. The function $t\mapsto \frac{d^2}{dt^2}\mathscr{F}_B(t)$ is first positive and then negative for the Pareto density in \eqref{Pareto} when $\tau \ge 4$. Consequently, for $\tau \ge 4$, there exists $\delta>0$ such that when $0\le B\le \delta$ there is a first-order phase transition.
\end{theorem}

\begin{remark}
[Critical value Pareto]{\rm
Using Theorem \ref{thm-critical-value} and applying integration by parts, the  critical values $(s^{\text{Par}}_c,\beta^{\text{Par}}_c)$ 
can be obtained from the root $t=t_c$
of $\mathscr{K}^{\text{Par}}(t)=0$, where
   \eqan{
    \mathscr{K}^{\sss\rm{Par}}(t)&=
    \frac{\tau-2}{\tau-1}\log\left(\frac{\e^{t }+q-1}{q}\right)
    -\left(\frac{q+1}{2q}-\frac{1}{\tau-1}\right) t\nonumber\\
    &\qquad+ \frac{(\tau-3)(q-1)t}{2(\tau-1)} \int_{1}^{\infty}\frac{(\tau-2)w^{-\tau+1}}{\e^{t w}+q-1}dw.    
    }
    Note the similarity between $\mathscr{K}^{\sss\rm{Par}}(t)$ in this remark and $\mathscr{K}^{\sss\rm{ER}}(t)$ in Remark \ref{remark-ER}, which is in fact exactly recovered in the limit limit $\tau\to\infty$. In the Pareto case, giving a closed formula for $t_c$ is impossible. An efficient iterative procedure for computing $t_c$ with arbitrary precision is discussed in Appendix \ref{appC}.
}\hfill \ensymboldefinition
\end{remark}

We next investigate the setting where $\tau\in (3,4]$ and $B=0,$ which turns out to be rather surprising:

\begin{theorem}[Pareto weights with $\tau\in \protect{(3,4]}$]
\label{thm-Pareto}
Fix $B=0$. Then for every $q\geq 3$, there exists a $\tau(q)\in (3,4]$ such that the function $t\mapsto \frac{d^2}{dt^2}\mathscr{F}_0(t)$ is first positive and then negative for the Pareto density in \eqref{Pareto} when $\tau\in(\tau(q),4]$, while $\frac{d^2}{dt^2}\mathscr{F}_0(t)$ is negative on $(0,\infty)$ for $\tau\in (3,\tau(q)].$ 
Consequently, for $\tau\in(\tau(q),4]$, the phase transition is {\em first order}, while for $\tau\in (3,\tau(q)],$ the phase transition is {\em second order} with critical value $\beta_c=\log(1+q/\nu)$ with $\nu=\expec[W^2]/\expec[W]$. 
Moreover, the critical value $\tau(q)$ satisfies
\be
\int_0^\infty x^{3-\tau(q)}\frac{(q-1)\e^x-\e^{2x}}{(\e^x+q-1)^3} dx=0.
\ee
Finally, as $q\rightarrow \infty,$ $\tau(q)= 3+\frac{\log{q}}{q}(1+o(1)).$
\end{theorem}



\paragraph{\bf A general approach to first-order phase transitions}
We finally give a general condition for the occurrence of a first-order phase transition, which avoids the zero-crossing condition:

\begin{theorem}[General condition for first-order phase transition]
\label{thm-general-first-order}
Let $q\geq 3$, $B=0$ and assume that $\expec[W^3]<\infty$. Then, there exists a critical value $\beta_c(0)$ such that $\beta\mapsto s^\star(\beta,0)$ makes a jump at  $\beta_c(0)$, where $s^\star(\beta,0)$ is the optimal value in \ref{var-problem-psi-B>0}. Thus the annealed Potts model has a first-order phase transition at $\beta_c(0)$. Further, $s^\star(\beta,0)=0$ for $\beta\in (0,\beta_c(0))$ and $s^\star(\beta,0)>0$ for all $\beta > \beta_c(0)$.
\end{theorem}
\begin{remark}[Extension to $B>0$]
{\em
By a continuity argument, we expect the same result to hold for $B>0$ sufficiently small,
the difference being that $s^\star(\beta,B)>0$ for $\beta\in (0,\beta_c(B))$. This is further discussed in Section \ref{sxec-ext-B>0-general-first-order-phase transition}. }
\end{remark}
\begin{remark}[Possible non-uniqueness jump]
\label{rem-general-first-order-PT}
{\rm Theorem \ref{thm-general-first-order} shows that, for $q\geq 3$ and $\expec[W^3]<\infty$, the occurrence of first-order phase transitions is {\em general.} We will see that this follows from the fact that, when $\expec[W^3]<\infty$, $\mathscr{F}_0(t)$ is convex for small $t$, and well-defined for $t=0$. Theorem \ref{thm-general-first-order} shows that generally $\beta\mapsto s(\beta, 0)$ jumps at $\beta=\beta_c(0)$. However, we cannot prove that this discontinuity is {\em unique}, since we cannot rule out multiple jumps. Multiple jumps would be rather surprising though.}
\hfill \ensymboldefinition
\end{remark}

\subsection{Discussion and open problems}
\label{sec-disc}
In this section, we discuss our results and state some open problems.
\medskip

\paragraph{\bf Unique jump.} In the case where the zero-crossing condition 
does not hold, we know that the phase transition is first order for $q\geq 3$ and $0\leq B<\log(q-1)$ when $\expec[W^3]<\infty$, but we cannot prove that the jump at $\beta_c$ is {\em unique}. It is of great interest to investigate this further.
\medskip

\paragraph{\bf Identification of the critical value.} In the case where the phase transition is first order, we can give an analytic criterion for the critical value $\beta_c$, but we cannot compute it explicitly except in the homogeneous case. Is there, aside from the homogeneous case, another example where one can give an {\em explicit} formula for the critical value?
\medskip

\paragraph{\bf First-order phase transition for non-zero magnetic field.}
When the zero-crossing condition holds, the first-order phase transition persists in the presence of a magnetic field that is small enough. In the homogeneous set-up, the case of non-zero external field was considered in \cite{biskup2006mean,blanchard2008thermodynamic}, and the first-order transition remains on a {\em critical line}, which is computed explicitly. Remarkably, 
the maximum value of the magnetic field for which a phase transition exists (i.e., $B_{\text{max}}=\log(q-1)$) is the same for all rank-1 inhomogeneous random graphs and coincides with the value that is also found in the homogeneous case.

\medskip
\paragraph{\bf Critical exponents for Pareto weights with second-order phase transitions.} Critical exponents give a wealth of information about the nature of the phase transition.
when this transition is second order. 
For Pareto weights and $\tau\in(3,\tau(q)]$, it would be of interest to identify the critical exponents for the magnetization, susceptibility, internal energy and specific heat. This may be possible by linearising the fixed-point equation around the critical value, as is common for Ising models.
\medskip

\paragraph{\bf Smoothed transition.} 
The smoothing from a first-order phase transition to a second-order one proved in Theorem \ref{thm-Pareto} for the annealed Potts model on rank-1  random graphs with Pareto distributions having exponent $3<\tau\le\tau(q)$ has been predicted in \cite{igloi2002first}. We make rigorous sense of this prediction and see that it is a phenomenon occurring in the annealed setting. It is an open question whether the same phenomenon also occurs in the quenched setting. 

The smoothing on the scale-free networks in Theorem \ref{thm-Pareto} is reminiscent of the smoothing proved in \cite{gobron2007first} on lattices. There, for the Potts model with Kac potential in two spatial dimensions and $q=3$, it is proved that the system undergoes a change from first- to second-order phase transition by changing the finite range of the interaction. Comparing the two results, we see that while the smoothing on two-dimensional lattices occurs
only for $q=3$, in the rank-1 random graph setting with Pareto weights the smoothing happens for all $q\ge 3$.
\medskip

\paragraph{\bf Quenched setting.} In this paper, we have investigated the {\em annealed} Potts model, finding that the phase transition is often first order. The same is expected for the quenched setting, but so far, this has only been proved for the random regular graph \cite{ bencs2023random,dembo2013factor,dembo2014replica}. It would be of great interest, but probably quite hard, to compute the quenched pressure per particle, and use it to identify the nature of the phase transition for the quenched setting. In the mean-field setting, often the quenched and annealed phase transition behave similarly (even though the critical values may differ). Is this also true for the rank-1 random graph?
\medskip

\paragraph{\bf Organisation of this paper.}
This paper is organised as follows. In Section \ref{sec-potts-grg}, we identify the pressure per particle and prove 
Theorem \ref{thm-press-particle}. %
In Section \ref{sec-soln_form}, we analyse the variational problem in \eqref{var-problem-psi-B>0} and prove  Theorem \ref{thm-soln_form-positive}. 
In Section \ref{sec-tau(2,3)}, we prove that, when $\expec[W^2]=\infty$, the spontaneous magnetization is positive for all $\beta>0$, and thus obtain Theorem \ref{thm-spont-magnetization}. %
In Section \ref{sec-nature-pt}, we assume that $\expec[W^2]<\infty$ and that $t\mapsto \frac{d^2}{dt^2}\mathscr{F}_B(t)$ is first positive and then negative. This allows us to establish the existence of a first-order phase transition for $q\ge 3$, and prove Theorems \ref{thm-nature-pt-a} and \ref{thm-nature-pt-b}. 
In Section \ref{sec-proportions}, we study the order parameter and prove Theorem \ref{thm-proportions}.
In Section \ref{sec-unique-zero-crossing-or-not}, we identify the critical temperature in zero field and investigate the special examples in \eqref{Guido(3)}, \eqref{eq:small_support_condition} and the Pareto case, and prove Theorems \ref{thm-critical-value}, \ref{thm-density-Guido}, \ref{thm-Pareto} and \ref{thm-two-atoms}.
%
%
%
In Section \ref{sec-FO-general}, we give a general argument for the existence of a first-order phase transition and prove Theorem \ref{thm-general-first-order}. %
Finally, there are three appendices.
In Appendix \ref{appA} we show that the zero-crossing condition fails for a specific choice of $W$ with two atoms.
In Appendix \ref{appB} we prove the asymptotics of $\tau(q)$ stated in Theorem \ref{thm-Pareto}.
We close with Appendix \ref{appC} 
where we discuss the computation of the critical temperature via Theorem \ref{thm-critical-value}.

\section{Existence of the partition function: Proof of Theorem \ref{thm-press-particle}}
\label{sec-potts-grg}
\invisible{\RvdH{TO DO 2: Clean up the proof of this section. Remco will do this, Neeladri will do the second round.\\
Remco has done his round of polishing.}
\NM{The proof reads okay to me. Should I add a small paragraph with Varadhan's lemma to explain the proof strategy a little at a heuristic level?}}

\subsection{Setting the stage}
Averaging over the randomness and recalling that in the rank-1 inhomogeneous random graph the occupation variables $I_{i,j}$ are independent Bernoulli variables with success probability $p_{ij},$
\begin{align}
\E [ Z_n(\beta, \bB)]= &\E \left ( \sum_{\boldsymbol{\sigma} \in \omqn} \exp \left [ \beta \sum_{i<j}  I_{i,j}  \indic{\s_i=\s_j} + \sum_{i\in [n]} \sum_{k\in [q]} B_k \indic{\s_i=k}   \right ]  \right )\\\nonumber
&=   \sum_{\boldsymbol{\sigma} \in \omqn} \e^ {\sum_{i\in [n]} \sum_{k\in [q]} B_k \indic{\s_i=k} }  \prod_{i<j} \E \left [  \e^{\beta   I_{i,j}  \indic{\s_i=\s_j}} \right ]\\\nonumber
&=  \sum_{\boldsymbol{\sigma} \in \omqn} \e^ {\sum_{i\in [n]} \sum_{k\in [q]} B_k \indic{\s_i=k} }  \prod_{i<j} \left (  \e^{\beta  \indic{\s_i=\s_j}} p_{ij}+ (1-p_{ij} )\right),
\end{align}
where $p_{ij}=\pp(I_{ij}=1)$. We introduce
$$
\beta_{ij}=\log \left(1+ (\e^\beta-1) p_{ij}  \right),
$$
and note that $\beta_{ij}$ satisfies
$$
 \e^{\beta  \indic{\s_i=\s_j}} p_{ij}+ (1-p_{ij} )= \e^{\beta_{ij}  \indic{\s_i=\s_j}}
$$
for any $\s_i$ and $\s_j$. Therefore, we can write
\be
\label{eq:z-ann}
\E [ Z_n(\beta, \bB)]= G_n(\beta)  \sum_{\boldsymbol{\sigma} \in \omqn} 
\e^{ \frac 1 2 \sum_{i,j\in [n]} \beta_{ij}  \indic{\s_i=\s_j} }
\e^ {\sum_{i\in [n]} \sum_{k\in [q]} B_k \indic{\s_i=k},} 
\ee
where we have used $\beta_{ij}=\beta_{ji}$, and we have introduced self-loops in the Hamiltonian that give rise to the factor 
    \eqn{
    G_n(\beta)=\prod_{i \in [n]} \e^{-\beta_{ii}/2}.
    }
Apart from this factor, \eqref{eq:z-ann}
is the partition function of an inhomogeneous Potts model on the complete graph with Hamiltonian, for $\bsigma\in \omqn,$
\be
\label{eq:in-potts}
\mathcal{H}_n(\boldsymbol{\sigma})= \frac 1 2 \sum_{i,j\in [n]}  \beta_{i,j}  \indic{\s_i=\s_j} + \sum_{i\in [n]} \sum_{k\in [q]} B_k \indic{\s_i=k}.
\ee
We expand $\beta_{ij}$ for small $p_{ij}$ as
$$
\beta_{ij}=(\e^\beta-1) p_{ij} + O(p_{ij}^2),
$$
and rewrite \eqref{eq:in-potts} as
\be
\label{eq:in-potts2}
\frac12 (\e^\beta-1) \sum_{i,j\in [n]}  p_{i,j}  \indic{\s_i=\s_j} + \sum_{i\in [n]} \sum_{k\in [q]} B_k \indic{\s_i=k}+ O ( \sum_{i,j\in [n]} p_{ij}^2 ).
\ee
Now we make use of the explicit structure of the $\GRGw$ with {\em deterministic weights}  $\bfwit=(w_i)_{i\in [n]}$ by setting
$$
p_{ij}=\frac{w_i w_j}{\ell_n + w_i w_j},
$$
where $\ell_n= \sum_{i\in [n]} w_i$ and assume the  weight regularity Condition \ref{cond-WR-GRG}. Under this condition, we will show below that the big-O term in \eqref{eq:in-potts2} is $o(n)$. Writing
$$
p_{ij}=\frac{w_i w_j}{\ell_n}+O\left (\left (\frac{w_i w_j}{\ell_n}\right)^2\wedge 1\right ),
$$
we obtain
$$
 \sum_{i,j\in [n]}  p_{i,j}  \indic{\s_i=\s_j} =  \sum_{i,j\in [n]} \frac{w_i w_j}{\ell_n} \indic{\s_i=\s_j}+O \left (\sum_{i,j\in [n]}  \left (\frac{w_i w_j}{\ell_n}\wedge 1\right)^2 \right).
$$
We next investigate the error term, which we can bound from above as
    \eqn{
    \sum_{i,j\in [n]}  \left (\frac{w_i w_j}{\ell_n}\right)^2
    \leq \frac{C}{n^2} \sum_{i,j\in [n]}  \left (w_i w_j\wedge n\right)^2
    =C \expec\Big[(W_1^{\sss(n)}W_2^{\sss(n)}\wedge n)^2\Big],
    }
where $W_1^{\sss(n)}$ and $W_2^{\sss(n)}$ are two i.i.d.\ random variables with distribution function as in \eqref{Fn-def}. Recall the power-law tail assumption in \eqref{main-assumption-Fn}, and let $\bar{W}_1$ and $\bar{W}_2$ be two i.i.d.\ random variables with complementary distribution function, for $w>c^{1/(\tau-1)},$
    \eqn{
    \prob(\bar{W}>w)=\frac{c}{w^{\tau-1}}.
    }
Then, \eqref{main-assumption-Fn} implies that $W_1^{\sss(n)}$ and $W_2^{\sss(n)}$ are stochastically bounded from above by $\bar{W}_1$ and $\bar{W}_2$, so that
    \eqn{
    \sum_{i,j\in [n]}  \left (\frac{w_i w_j}{\ell_n}\right)^2
    \leq C \expec\Big[(W_1^{\sss(n)}W_2^{\sss(n)}\wedge n)^2\Big]\leq
    C\expec\Big[(\bar{W}_1\bar{W}_2\wedge n)^2\Big].
    }
Since $\bar{W}_1$ and $\bar{W}_2$ have power-law complementary distribution function, it is not hard to show that there exists a constant $C>1$ such that
    \eqn{
    \prob(\bar{W}_1\bar{W}_2>x)\leq \frac{C\log x}{x^{\tau-1}}.
    }
This implies that, uniformly in $a>1$, for some $C>1$, since $\tau\in(2,3)$,
    \eqn{
    \expec\Big[(\bar{W}_1\bar{W}_2\wedge a)^2\Big]
    \leq C a^{3-\tau}\log{a}.
    }
We conclude that
    \eqn{
    \sum_{i,j\in [n]}  \left (\frac{w_i w_j}{\ell_n}\right)^2
    \leq C\expec\Big[(\bar{W}_1\bar{W}_2\wedge n)^2\Big]
    \leq C n^{3-\tau}\log{n}.
    }
Since $\tau>2$, this is $o(n)$.

All in all, the inhomogeneous Potts Hamiltonian on the complete graph defined in \eqref{eq:in-potts} is 
\be
\label{eq:h-in-potts}
\mathcal{H}_n(\boldsymbol{\sigma})=\frac{\e^\beta-1}{2} \sum_{i,j\in [n]} \frac{w_i w_j}{\ell_n} \indic{\s_i=\s_j} + \sum_{i\in [n]} \sum_{k\in \omq} B_k \indic{\s_i=k} +o(n).
\ee
Since 
$$
 \indic{\s_i=\s_j} = \sum_{k \in \omq} \indic{\s_i=k}  \indic{\s_j=k},
$$
\eqref{eq:h-in-potts} can be rewritten as 
\be
\mathcal{H}(\boldsymbol{\sigma})=\frac{\e^\beta-1}{2\ell_n}\sum_{k\in \omq} \Big( \sum_{i\in [n]} w_i  \indic{\s_i=k}\Big)^2+  \sum_{k\in \omq} B_k  \sum_{i\in [n]}\indic{\s_i=k} +o(n).
\ee
Let $U_n \in [n]$ be a uniformly chosen vertex in $\GRGw$ and $W_n=w_{U_n}$ its weight, then  $\expec[W_n]= \ell_n/n$. Recalling the notation \eqref{betaprime}, we write
\be
\mathcal{H}(\boldsymbol{\sigma})=\frac{\beta'}{2n\,   \expec[W_n]}\sum_{k\in \omq} \Big( \sum_{i\in [n]} w_i  \indic{\s_i=k}\Big)^2+  \sum_{k\in \omq} B_k  \sum_{i\in [n]}\indic{\s_i=k} +o(n).
\ee
Finally, we rewrite  \eqref{eq:z-ann} as
\eqan{\label{eq:z-ann-2}
&\E [ Z_n(\beta, \bB)]\nn\\
&= 
G_n(\beta) \e^{o(n)}  \sum_{\boldsymbol{\sigma} \in \omqn} 
\exp \left [   \frac{\beta'}{2n\,   \expec[W_n]}\sum_{k\in \omq} \Big( \sum_{i\in [n]} w_i  \indic{\s_i=k}\Big)^2  \right ]
\exp \Big[  \sum_{k\in \omq} B_k  \sum_{i\in [n]}\indic{\s_i=k}   \Big] 
.
}

\invisible{Our aim is to compute the annealed pressure
$$
\psi(\beta,\bB) = \lim_{n\to\infty} \frac{1}{n} \log \mathbb{E}[Z_n(\beta,\bB)].
$$}

\subsection{Applying Hubbard-Stratonovich}
\label{sec-Hub-Strat}
We introduce the random vector $\kZ=(\kZ_1,\ldots,\kZ_q)$ with $\kZ_k$ independent standard Gaussian. Using the identity
$$\expec[\e^{t_1 \kZ_1+\cdots  + t_q \kZ_q}]=\e^{\frac{1}{2 }(t_1^2+\cdots + t_q^2)},$$ 
valid for all $t_k \in \R$, we compute the quadratic term in \eqref{eq:z-ann-2} as
$$
 \exp \left [   \frac{\beta'}{2n\,   \expec[W_n]}\sum_{k\in \omq} \Big( \sum_{i\in [n]} w_i  \indic{\s_i=k}\Big)^2  \right ]= \expec_\kZ \left [  \exp \left ( \sqrt{\frac{\beta'}{n \expec[W_n]}}  \sum_{k \in  \omq}\sum_{i\in [n]} w_i  \indic{\s_i=k} \kZ_k \right)\right ],
$$
where $\expec_\kZ [\cdot]$ denotes the expectation with respect to the normal vector $\kZ$.
Thus, exchanging the sums over $k \in  \omq$ and $i\in [n]$ we can write $\eqref{eq:z-ann-2}$ as follows
\begin{align}
\nonumber
&\E [ Z_n(\beta, \bB)]\\
&=G_n(\beta) \e^{o(n)} \hspace{-0.2cm}  \sum_{\boldsymbol{\sigma} \in \omqn} 
\expec_\kZ \left [  \exp \left ( \sqrt{\frac{\beta'}{n \expec[W_n]}} \sum_{i\in [n]}   \sum_{k \in  \omq} w_i  \indic{\s_i=k} \kZ_k \right)\right ]
\exp \left [   \sum_{i\in [n]} \sum_{k\in \omq} B_k \indic{\s_i=k}   \right ] 
\\\nonumber
& = G_n(\beta) \e^{o(n)}  \expec_\kZ \left [ \sum_{\sigma_1\in \omq}\cdots \sum_{\sigma_n \in \omq} \prod_{i \in [n]} \exp \left (  
\sqrt{\frac{\beta'}{n \expec[W_n]}}  \sum_{k \in  \omq} w_i  \indic{\s_i=k} \kZ_k +\sum_{k\in \omq} B_k \indic{\s_i=k} 
    \right )\right  ]\\\nonumber
&= G_n(\beta) \e^{o(n)}   \expec_\kZ \left [  \prod_{i \in [n]}   \sum_{\sigma_i\in \omq}\ \exp \left (  \sum_{k\in \omq} \left \{  \sqrt{\frac{\beta'}{n \expec[W_n]}}  w_i \kZ_k +B_k \right \} \indic{\s_i=k}\right ) \right  ]\\\nonumber
&= G_n(\beta) \e^{o(n)}    \expec_\kZ \left [  \prod_{i \in [n]} \sum_{k\in \omq} \exp \left (  \sqrt{\frac{\beta'}{n \expec[W_n]}}  w_i  \kZ_k + B_k \right )  \right  ].
\end{align}
In the last step we have used the fact that for each choice of $\sigma_i \in \omq $ only one term of the sum in the exponential is different from zero. Introducing the function
\eqn{
\label{f-n-w-z-def}
f_{n;\beta,\bB}(w;z_1,\ldots z_q)=\e^{ \sqrt{\frac{\beta'}{n \expec[W_n]}}  w z_1 +B_1}+\cdots + \e^{\sqrt{\frac{\beta'}{n \expec[W_n]}}  wz_q +B_q},
}
we can write the annealed partition function as
\begin{align}
\E [ Z_n(\beta, \bB)]&=G_n(\beta) \e^{o(n)} \expec_\kZ  \left [ \prod_{i\in [n]}  f_{n;\beta,\bB}(w_i;\kZ_1,\ldots, \kZ_q)   \right ]\\\nonumber
& = G_n(\beta) \e^{o(n)} \expec_\kZ  \left [ \exp \sum_{i \in [n]} \log f_{n;\beta,\bB}(w_i; \kZ_1,\ldots, \kZ_q) \right ]\\\nonumber
&= G_n(\beta) \e^{o(n)} \expec_\kZ  \left [   \e^{n \expec_{W_n} \left [ \log f_{n;\beta,\bB}(W_n;\kZ_1,\ldots,\kZ_q)\right ]}   \right ] \\\nonumber
&= G_n(\beta) \e^{o(n)} \expec_\kZ  \left [ \e^{n \tilde{F}_{n;\beta,\bB}(\kZ_1,\ldots,\kZ_q )} \right ],
\end{align}
where
\begin{align}
\tilde{F}_{n;\beta,\bB}(z_1,\ldots,z_q )&= \expec_{W_n} \left [ \log f_{n;\beta,\bB}(W_n;z_1,\ldots,z_q)\right ]\\\nonumber
 &=\expec_{W_n} \left [ \log \left ( \e^{z_1  \sqrt{\frac{\beta'}{n \expec[W_n]}}  W_n +B_1}+\cdots + \e^{z_q  \sqrt{\frac{\beta'}{n \expec[W_n]}}  W_n +B_q} \right )\right ],
 \end{align}
 and the expectation in $\expec_{W_n}$ is w.r.t.\ $W_n$. Then, applying the change of variable $z_k=y_k \sqrt{\frac{n \beta'}{\mathbb{E}[W_n]}}$, we arrive at
 \begin{align}\label{eq:partf}
  &\E [ Z_n(\beta, \bB)]=G_n(\beta) \e^{o(n)} {(2 \pi)^{-\frac q 2}} \int_{\R^q} dz_1\cdots d z_q \exp \left [n \tilde{F}_{n;\beta,\bB}(z_1,\ldots,z_q )-\frac 1 2  \sum_{i=1}^q z_i^2 \right ]\\\nonumber
  &= G_n(\beta) \e^{o(n)}  \left (\frac{n \beta' }{2 \pi \expec[W_n]} \right)^{\frac q 2} \int_{\R^q} dy_1\cdots d y_q  \exp \Bigg\{n  \left [  {F}_{n;\beta,\bB}(\by)   -\frac{\beta'}{2 \mathbb{E}[W_n]}  \sum_{i=1}^q y_i^2 \right ]\Bigg\},
 \end{align}
where
\be\label{eq:FN}
 {F}_{n;\beta,\bB}(\by) = \expec_{W_n} \left [ \log \left ( \e^{\beta' y_1   \frac{W_n}{\expec[W_n]}  +B_1}+\cdots + \e^{ \beta' y_q  \frac{W_n}{ \expec[W_n]}  +B_q} \right )\right ].
\ee

\subsection{Annealed pressure: conclusion of the proof of Theorem \ref{thm-press-particle} }
We compute the thermodynamic limit of the pressure
\begin{align}\label{eq:limpress}
\varphi(\beta,\bB)&= \lim_{n\to \infty} \frac 1 n \log \E [ Z_n(\beta, \bB)]= 
 \lim_{n\to \infty} \frac 1 n \log G_n(\beta) \\\nonumber
&\quad+ \lim_{n\to \infty} \frac 1 n \log  \int_{\R^q} dy_1\cdots d y_q  \exp \Bigg\{n  \left [  {F}_{n;\beta,\bB}(\by)   -\frac{\beta'}{2 \mathbb{E}[W_n]} \sum_{i=1}^q y_i^2 \right ]\Bigg\}.
\end{align}
The limit of the integral can be evaluated by Varadhan's Lemma \cite{varadhan1966asymptotic}.
\invisible{For the time being,  we assume that the function $f_n$ can be substituted with its limit
\be
\label{eq:effe}
 {F}(\bB; \by) = \expec_{W} \left [ \log \left ( \e^{\beta' y_1   \frac{W}{\expec[W]}  +B_1}+\cdots + \e^{ \beta' y_q  \frac{W}{ \expec[W]}  +B_q} \right )\right ],\ee
$W$ being the limiting weight. 
Recalling that
\be\label{eq:FN2}
 {F}_n(\bB; \by) = \expec_{W_n} \left [ \log \left ( \e^{\beta' y_1   \frac{W_n}{\expec[W_n]}  +B_1}+\cdots + \e^{ \beta' y_q  \frac{W_n}{ \expec[W_n]}  +B_q} \right )\right ]
\ee
and 
\be
\label{eq:effe2}
 {F}(\bB; \by) = \expec_{W} \left [ \log \left ( \e^{\beta' y_1   \frac{W}{\expec[W]}  +B_1}+\cdots + \e^{ \beta' y_q  \frac{W}{ \expec[W]}  +B_q} \right )\right ].\ee
 }
 Recalling the notation introduced in \eqref{F-pressure}, i.e.,
 \be
\label{eq:effe2}
 F_{\beta,\bB}( \by) = \expec_{W} \left [ \log \left ( \e^{\beta' y_1   \frac{W}{\expec[W]}  +B_1}+\cdots + \e^{ \beta' y_q  \frac{W}{ \expec[W]}  +B_q} \right )\right ],\ee
 we first prove the following result:
 \begin{lemma}[Pointwise convergence $F_{n,\beta}$]
 \label{lemma:convfn}
     Let $M>0$. Then, for all $\by\in \R^q$ with $\|\by\|_\infty < M$,
     \be
     \lim_{n\to\infty} {F}_{n;\beta,\bB}(\by)=F_{\beta,\bB}( \by ).
     \ee
 \end{lemma}

 
\begin{proof}
For the sake of notation we introduce the normalized sequence $\wpn=W_n/\expec[W_n]$, and observe that, by Condition \ref{cond-WR-GRG}, $\wpn \stackrel{\mathcal D}{\longrightarrow} W^\prime$, where $W^\prime=\frac{W}{\expec[W]}$. 
Since $\expec[\wpn] \to \expec[W^\prime]$, we have (see \cite[Theorem 3.6]{billingsley2013convergence}) that the variables $\wpn$ are uniformly integrable.\\
Now we consider the variables
\be
X_n= h(\wpn)\equiv \log \Big(\sum_{i=1}^q \e^{\bpr \wpn y_i +B_i}\Big).
\ee
By hypothesis $|y_i|\le M$,
so that
$$
|X_n|\le \log q + ||B||_\infty + M\, \bpr\, \wpn.
$$
From the uniform integrability of $\wpn$ we have that the r.h.s.\ of the previous inequality is uniformly integrable and the same is also true for $X_n=h(\wpn)$, as $h$ is a continuous function. The convergence in distribution of $\wpn$ and the uniform integrability of $h(\wpn)$ imply (by \cite[Theorem A.8.6]{ellis2007entropy}) the convergence of 
$\expec[h(\wpn)]$ to $\expec[h(W)]$, that is, $ {F}_{n;\beta,\bB}(\by)\to  F_{\beta,\bB}( \by)$ for $\|\by\|_\infty < M$.
 \end{proof}
The proof of Theorem \ref{thm-press-particle}, i.e., that the limit  
\eqref{eq:limpress}
is given by \eqref{eq:laplace},  can be obtained by applying  Varadhan's Lemma \cite{varadhan1966asymptotic} 
to compute the limit of $\frac {1}{a_n} \log \int_\Omega \exp\{a_n f_n(\by)\} P_n(d\by)$. Here $(P_n)_n$ is a sequence of probability measures on a regular topological space $\Omega$ that satisfies a large deviation property with rate function $\by\mapsto I(\by)$ and speed $a_n \to \infty$. On $\Omega$, the functions $f_n(\by)$ defined in \eqref{f-n-w-z-def} are required to satisfy a set of conditions weaker than the convergence stated in Lemma \ref{lemma:convfn}, plus the condition
\be\label{eq:weakbound}
\lim_{L\to \infty} \limsup_{n\to \infty} \frac 1 a_n \log \int_{f_n(\by)\ge L} \exp\{a_n f_n(\by)\} P_n(d\by) = -\infty
\ee
(which is less restrictive than uniform boundedness from above).
\invisible{The proof of \eqref{eq:laplace} can be obtained by applying a generalization of Varadhan's Lemma \cite{varadhan1966asymptotic} that gives the limit of the sequence $\frac {1}{a_N} \log \int_\Omega \exp\{a_N f_n(y)\} dP_n$, where $(P_n)_N$ is a sequence of probability measures on a regular topological space $\Omega$ that satisfies a large deviation property for $a_N \to \infty$ with rate function $I(y)$, while $f_n(y)$ are functions on $\Omega$ that fulfill a set of conditions  weaker than the convergence established in Lemma \ref{lemma:convfn}, plus the following one (which is less restrictive than the uniform boundedness from above of the family $(f_n)_{n\geq 1}$):
\be\label{eq:weakbound}
\lim_{L\to \infty} \limsup_{n\to \infty} \frac 1 a_n \log \int_{f_n(\by)\ge L} \exp\{a_n f_n(\by)\} dP_n = -\infty.
\ee}
In our case, see \eqref{eq:limpress}, 
$P_n$ is the $q$-fold product of normal distributions on $\Omega=\R^q$, with variance $\expec[W_n]/\bpr n$, i.e., $P_n(d\by)\propto \e^{-n \bpr\frac{|\by|^2}{2 \expec[W_n] }} dy_1\cdots dy_q$. The sequence $(P_n)_n$ satisfies a large deviation principle with speed $a_n=n$ and large deviation function 
\begin{equation}\label{eq:ratenormal}
  I(\by)=\frac 1 2 \frac{\bpr}{\expec[W]} |\by|^2,
\end{equation}
(where here and in the following $|\cdot|$ denotes the $L^2$ norm on $\R^q$).
This fact can be seen by observing that $P_n$ is the distribution of a random vector $Q_n\in \R^q$ with i.i.d.\ normal components, whose logarithmic moment generating function 
$$\log \phi_n(t)=\frac{1}{2}\frac{\expec[W_n]}{\bpr\, n} |t|^2$$ 
satisfies the hypotheses of the G\"artner-Ellis Theorem (see \cite[Section 2.3]{Dembo1998LargeDT}), i.e.,
$$ \lim_{n\to \infty} \frac 1 n \log 
 \phi_n(nt)=\frac{1}{2}\frac{\expec[W]}{\bpr} |t|^2.$$ 
 As a consequence, the Legendre transform of $\frac{1}{2}\frac{\expec[W]}{\bpr} |t|^2$, given in \eqref{eq:ratenormal}, is the large deviation function $\by\mapsto I(\by)$ of the sequence $(P_n)_n$.
 \begin{lemma}[Large values of $f_n$ do not contribute to exponential moment]
 \label{lemma:weakbound}
 The sequence of functions \eqref{eq:FN} satisfies the condition \eqref{eq:weakbound}.
 \end{lemma}
 \begin{proof}
For the sake of notation, we let $f_n(\by)$ be the function 
 defined in \eqref{eq:FN}. 
    Since there exist positive constants $A$ and $B$ such that $f_n(\by)\le A+ B |\by|$, it follows that
    \begin{align}\label{eq:ineq1}
        &\int_{f_n(\by)\ge L} \exp\{n f_n(\by)\} dP_n(\by)  \le
        \e^{A\,n}\int_{|\by|\ge (L-A)/B} \exp\{  n\, B |\by|\} P_n(d\by).
    \end{align}

\invisible{  \begin{align}\label{eq:ineq1}
        &\int_{f_n(\by)\ge L} \exp\{n f_n(\by)\} dP_n(\by) \le
        \int_{f_n(\by)\ge L} \exp\{  n(A+ B |\by|)\} P_n(d\by)\nonumber \\
        & \le
        \e^{A\,n}\int_{|\by|\ge (L-A)/B} \exp\{  n\, B |\by|\} dP_n(y).
    \end{align}
}    
    Therefore,
    \begin{align}\label{eq:ineq2}
       & \frac 1 n \log \int_{f_n(\by)\ge L} \exp\{n f_n(\by)\} dP_n (\by)\le A+  \frac 1 n \log \int_{|\by|\ge (L-A)/B} \exp\{  n\, B |\by|\} dP_n(\by)
        \\
        & = A+ \frac 1 n \log \int_{|\by|\ge L'} \exp \left \{  n\,|\by| \left (B-\frac{\bpr |\by|}{2 \expec[W_n]} \right)\right \} dy_1\cdots dy_q\nonumber \\
        &\le 
        A+ \frac 1 n \log \int_{|\by|\ge L'} \exp \left \{  n\,|\by| \left (B-\frac{\bpr L'}{2 \expec[W_n]} \right)\right \} dy_1\cdots dy_q, \nonumber
    \end{align}
    with $L'=(L-A)/B$. For any fixed $c>0$, taking  $L'$ and $n$ sufficiently large, the inequality  $(B-\frac{\bpr L'}{2 \expec[W_n]}) <-c$ is satisfied. Then,
    \eqan{
       &\frac 1 n \log \int_{|\by|\ge L'} \exp \left \{  n\,|\by| \left (B-\frac{\bpr L'}{2 \expec[W_n]} \right)\right \} dy_1\cdots dy_q\nn\\
       &\qquad\le
       \frac 1 n \log \int_{|\by|\ge L'} \exp \left \{  -c\, n\,|\by| \right \} dy_1\cdots dy_q.
    }
    In order to estimate the last integral we introduce polar coordinates in $\R^q$ with radial variable $r=|\by|$.  Assuming $L'$ to be sufficiently large, we can write the integral as 
    \begin{align}   
    \int_{|\by|\ge L'} \exp \left \{  -c\, n\,|\by| \right \} dy_1\cdots dy_p 
    &= a_{q-1} \int_{L'}^\infty \e^{-c\,n\,r} r^{q-1} d\,r \le a_{q-1} \int_{L'}^\infty \e^{(-c\,n\,+1)r} d\,r\nonumber \\
    & = a_{q-1} \frac{\e^{(-c\,n+1)L'}}{c\,n-1},\nonumber
     \end{align}
     where $a_{q-1}$ is the surface of the unit sphere in $\R^{q-1}$. Finally, putting together the previous inequalities, we get
     \begin{equation}
         \frac 1 n \log \int_{f_n(\by)\ge L} \exp\{n f_n(\by)\} dP_n (\by)
         \le A- cL' + \frac 1 n \log a_{q-1} -\frac 1 n \log{(c\,n-1)} + \frac 1 n L',
     \end{equation}
     which shows that
     $$
     \lim_{L\to \infty} \limsup_{n\to \infty} \frac 1 n \log \int_{f_n(\by)\ge L} \exp\{n f_n(\by)\} P_n (d\by) =-\infty.
     $$
 \end{proof}
 \begin{proposition}[Convergence of annealed pressure per particle]
     The thermodynamic limit \eqref{eq:limpress} of the finite-volume annealed pressure per particle \eqref{finite-pressure} satisfies \eqref{eq:laplace}.
 \end{proposition}
  \begin{proof}
  Varadhan's Lemma, whose hypotheses are satisfied by Lemmas  \ref{lemma:convfn} and \ref{lemma:weakbound}, implies that the second limit in the r.h.s.\ of \eqref{eq:limpress} is $\sup_{\by\in \R^q} [F_{\beta,\bB}( \by)-I(\by)]$, with $F_{\beta,\bB}(\by)$ and $I(\by)$ given in \eqref{eq:effe2} and \eqref{eq:ratenormal}. We further observe that the $\sup$ is in fact a $\max$. Indeed, from the upper bound $F_{\beta,\bB}( \by)-I(\by)\le q+|B|+\bpr \expec[W] |\by|-\frac{\bpr}{2\expec[W]}|\by|^2$, we get that  $F_{\beta,\bB}( \by)-I(\by)$ is strictly negative outside the compact set $|\by|\le M$, for $M>0$ sufficiently large. Inside this compact set, the function takes its global positive maximum value. 

  The first limit in the r.h.s.\ of \eqref{eq:limpress}   involves the function $G_n(\beta)=\prod_{i \in [n]} \e^{-\beta_{i,i}/2}$, where $\beta_{i,i}=\log(1+\frac{\bpr w_i^2}{\ell_n + w_i^2})$
  with $\ell_n=\sum_{i\in [n]} w_i$. Thus,
  \begin{align*}
     \left | \frac 1 n \log G_n(\beta) \right | & = \frac{1}{2n} \sum_{i=1}^n \log\Big(1+\frac{ \bpr w_i^2}{\ell_n + w_i^2}\Big) =\frac{1}{2}\expec\Big[\log\Big(1+\frac{ \bpr W_n^2}{\ell_n + W_n^2}\Big)\Big].
     \end{align*}
We next use Condition \ref{cond-WR-GRG}(b). Let $a\in (0,\tfrac{1}{2}),$ and split the expectation depending on whether $W_n\leq n^a$ or $W_n>n^a$. Since $x \to x/(\ell_n+x)$ is increasing, we can bound the first contribution by
    \eqn{
    \expec\Big[\log\Big(1+\frac{ \bpr W_n^2}{\ell_n + W_n^2}\Big)\indic{W_n\leq n^a}\Big]
    \leq \log\Big(1+\bpr n^{2a-1}\Big)=o(1),
    }
since, by the weight regularity condition, $\ell_n/n=\expec[W_n] \to \expec[W] \in(0,\infty)$. Further, using that $\frac{ \bpr W_n^2}{\ell_n + W_n^2}\leq \bpr$, we obtain that
    \eqn{
    \expec\Big[\log\Big(1+\frac{ \bpr W_n^2}{\ell_n + W_n^2}\Big)\indic{W_n> n^a}\Big]
    \leq \log\big(1+\bpr\big)\prob(W_n>n^a)=o(1),
    }
by Condition \ref{cond-WR-GRG}(b).
Thus, we conclude that
  $$
  \lim_{n \to \infty}  \frac 1 n \log G_n(\beta)=0,
  $$
\invisible{
We further observe that if the convergence of the second moments is lacking, the sequence $(\frac 1 n \log G_n(\beta))_n$ has still a converging sub sequence. Indeed, the sequence is bounded because 
  $\log(1+\frac{ \bpr w_i^2}{\ell_n + w_i^2}) \le \frac{ \bpr w_i^2}{\ell_n + w_i^2} \le \bpr$.}
and finally that
\begin{align}\label{eq:laplace-2}
\varphi(\beta,\bB) &= \lim_{n\to \infty} \frac 1 n \log  \int_{\R^q} dy_1\cdots d y_q  \exp \Bigg\{n  \left [  {F}_{n;\beta,\bB}(\by)    -\frac{\beta'}{2 \mathbb{E}[W_n]}  \sum_{i=1}^q y_i^2 \right ]\Bigg\}\nonumber\\
& = \max_{y_1,\ldots,  y_q} \left [  {F_\beta}(\bB; \by)    -\frac{\beta'}{2 \mathbb{E}[W]}  \sum_{i=1}^q y_i^2 \right ] .
\end{align}
This completes the proof of Theorem \ref{thm-press-particle}.
\end{proof}
\invisible{
\RvdH{TO DO 3: How do we know that this solution is {\em unique}? It probably need not be, right? If it is not, then we may have a problem with the analysis of the proportion of spins $X_k$ below.}
and the annealed pressure  is
\be\label{eq:pressure}
\varphi(\beta,\bB)=  \expec_{W} \left [ \log \left ( \e^{y_1^\star(\beta,\bB)  \frac{\beta'}{\expec[W]}  W +B_1}+\cdots +  \e^{y_q^\star(\beta,\bB)  \frac{\beta'}{\expec[W]}  W +B_q} \right )\right ]-\frac{\beta'}{2 \expec[W]} \sum_{i=1}^q y_i^\star(\beta,\bB)^2.
\ee
\RvdH{The part below is fishy, since we do not even know that the maximizer $\by^\star(\beta,\bB)$ is unique. How to define $x_k$ in \eqref{x-k-form} when $\by^\star(\beta,\bB)$ is not unique? We may have to abandon this part.}
The limit $x_k(\beta,\bB)$ of the average proportion $X_k$ of spin with value $k\in \omq$, i.e., 
\be
X_k(\beta,\bB) = \sum_{\sigma\in [q]^n} \frac{\frac{1}{n}\sum_{i=1}^n \indic{\sigma_i=k} \E[\exp[{\mathcal H(\boldsymbol{\sigma})}]] }{{\E [ Z_n(\beta, \bB)]}},
\ee
can be obtained by computing the derivative
\be
x_k(\beta,\bB) =\frac{\partial}{\partial B_k} \varphi(\beta,\bB).
\ee
This gives
 \begin{align}
 x_k(\beta,\bB)&= \expec_W \left [ \frac{\sum_{j\in[q]} \e^{y_j^\star(\beta,\bB)  \frac{\beta'}{\expec[W]} W + B_j} \frac{\beta'}{\expec[W]} W  \frac{\partial y_j^\star(\beta,\bB)}{\partial B_k} +  \e^{y_k^\star(\beta,\bB)  \frac{\beta'}{\expec[W]} W + B_k }}{ \sum_{j\in[q]} \e^{y_j^\star(\beta,\bB)  \frac{\beta'}{\expec[W]}  W +B_j}} \right ]\\\nonumber
&\qquad -\frac{\beta'}{ \expec[W]}\sum_{j \in \omq} y_j^\star(\beta,\bB)  \frac{\partial y_j^\star(\beta,\bB)}{\partial B_k} \\
&=
 \frac{\beta'}{ \expec[W]}\ \expec_W \left [ \sum_{j\in[q]} \left (\frac{ \e^{y_j^\star(\beta,\bB)   W + B_j}  W }{ \sum_{i\in \omq} \e^{y_i^\star(\beta,\bB)  \frac{\beta'}{\expec[W]}  W +B_i}} - y_j^\star(\beta,\bB) \right) \frac{\partial y_j^\star(\beta,\bB)}{\partial B_k}\right ]\nn\\\nonumber
&\qquad+\expec_W \left [  \frac{\e^{y_k^\star(\beta,\bB)  \frac{\beta'}{\expec[W]} W + B_k }}{ \sum_{j\in[q]} \e^{y_j^\star(\beta,\bB)  \frac{\beta'}{\expec[W]}  W +B_j}} \right ].\nn
 \end{align}
 Since the $y_j^\star(\beta,\bB)$'s  satisfy the stationarity conditions \eqref{eq:stat0}, the first term in the r.h.s.\ of the previous equation vanishes and we conclude that
 \be
 \label{x-k-form}
 x_k(\beta,\bB)=\expec_W \left [  \frac{\e^{y_k^\star(\beta,\bB)  \frac{\beta'}{\expec[W]} W + B_k }}{ \sum_{j\in[q]} \e^{y_j^\star(\beta,\bB)  \frac{\beta'}{\expec[W]} W +B_j}} \right ].
 \ee}

\section{Form of the optimizer: Proof of Theorem \ref{thm-soln_form-positive}}
\label{sec-soln_form}


In this section, we prove that the optimizer of \eqref{eq:laplace} is of the form \eqref{eq:soln_form}, as stated in Theorem \ref{thm-soln_form-positive}. 
First we 
need a useful technical lemma: 
\begin{lemma}[A useful convexity property]
\label{lem:g_convexity}
    Let $X$ be any random variable with $\E[X]< \infty$. Consider any finite collection $\{y_i\colon  i \in [L]\}$ of non-negative reals, and let $y_{\max}=\max\{y_i
    \colon  i \in [L]\}$. Then, the function 
    $$g(x):=\frac{1}{x}\E\left[\frac{X\e^{x X}}{\sum_{j=1}^L \e^{y_j X}} \right]$$ 
    is a convex continuous function of $x$ on $(0,y_{\max}]$. In particular, either $g$ is strictly decreasing on $(0,y_{\max})$, or $g$ has a unique global minimum at some $0<x_*<y_{\max}$, and $g$ is strictly decreasing on $(0,x_*)$, while strictly increasing on $(x_*,y_{\max})$. Furthermore, $\lim_{x \to 0}g(x)=\lim_{x \to \infty}g(x)=\infty$. 
\end{lemma}
\begin{proof}
    The second statement of the lemma clearly follows from the first statement, and the limiting statement $\lim_{x \to 0}g(x)=\lim_{x \to \infty}g(x)=\infty$ is a routine check, so we only verify the first statement of the lemma.
    Note that $x\mapsto g(x)$ is a well-defined continuous function on $(0,y_{\max}]$. To verify convexity, for any $x,y \in (0,y_{\max}]$ and $\theta \in [0,1]$, we let $h_X(x):=\e^{x X}/x$, and note that
    \begin{align*}
        g(\theta x + (1-\theta) y)&=\E\left[\frac{X}{\sum_{j=1}^L \e^{y_j X}} h_{X}(\theta x+(1-\theta)y)\right]\\&\leq \E\left[\frac{X}{\sum_{j=1}^L \e^{y_j X}} (\theta h_{X}(x)+(1-\theta)h_X(y))\right]\\&=\theta \E\left[\frac{Xh_X(x)}{\sum_{j=1}^L \e^{y_j X}} \right]+(1-\theta)\E\left[\frac{Xh_X(y)}{\sum_{j=1}^L \e^{y_j X}} \right]\\&= \theta g(x)+(1-\theta)g(y),
    \end{align*}
    where to obtain the second inequality above we use the fact that for any fixed $X$, the function $h_X(x)$ is a convex function of $x$. This verifies convexity of $x\mapsto g(x)$. 
\end{proof}
\begin{remark}[Restriction of domain to $(0,y_{\max}\protect{]}$]
    {\rm Note that unless $X$ has finite exponential moments in some neighborhood of $0$, the function $g$ might not be defined for $x=r$ for some $r>y_{\max}$, which is why we restrict to the domain $(0,y_{\max}]$.  }\hfill\ensymboldefinition
\end{remark}

Let us now consider the form of the optimizer of \eqref{eq:laplace} when $B$ is strictly positive:

\begin{proposition}[General solution structure for a field with one non-zero coordinate]
\label{prop-form-solution-B>0}
    Consider a field of the form $\bB=(B,0,\dots,0)$ for arbitrary $B>0$. Then any optimizer of the variational problem \eqref{eq:laplace} has the form \eqref{eq:soln_form}.
\end{proposition}
\begin{proof} 

Let us define the function
\begin{equation}\label{eq:g_y_B}
        g_{\by,\bB}(x)=\frac{1}{x}\E\left[\frac{W\e^{x\frac{\beta'W}{\E[W]}}}{\sum_{j=1}^q\e^{y_j(\beta,\bB)\frac{\beta'W}{\E[W]}+B_j}}\right],
    \end{equation}
    where $(y_1(\beta,\bB),\ldots,y_q(\beta,\bB))$ is an optimizer of \eqref{eq:laplace-2}.
Note that the stationarity conditions \eqref{eq:stat0} read
\begin{align*}
    \e^{B}g_{\by,\bB}(y_1(\beta,\bB))=g_{\by,\bB}(y_i(\beta,\bB)),
\end{align*}
whenever $i \in \{2,\dots,q\}$. It is straightforward to check that, since $y_i(\beta,\bB)$ are non-negative (see Remark \ref{rem-optim-positive}),  the function $g_{\by,\bB}$ satisfies the properties of Lemma \ref{lem:g_convexity}, and so the conclusions of Lemma \ref{lem:g_convexity} are true for it. Let us define $y_{\max}=y_{\max}(\beta,\bB):=\max\{y_i(\beta,\bB)\colon  i \in [q]\}$.

We now do a case-by-case analysis:
\smallskip

\textbf{Case $\min_{x \in (0,y_{\max}]}g_{\by,\bB}(x)=g(y_{\max})$.} In this case, using the convexity and continuity properties of $g_{\by,\bB}$, we note that the function $g_{\by,\bB}$ is a continuous strictly decreasing function on $(0,y_{\max}]$. Consequently, $y_i(\beta,\bB)=y_j(\beta,\bB)$ holds whenever $i,j \in \{2,\dots, q\}$, using the stationarity conditions (i.e., $ \e^{B}g_{\by,\bB}(y_i(\beta,\bB))=g_{\by,\bB}(y_j(\beta,\bB)))$ and injectivity. On the other hand, letting $r=g_{\by,\bB}(y_2(\beta,\bB))$, we  again use the stationarity conditions and $B>0$ to note that $g_{\by,\bB}(y_1(\beta,\bB))=\e^{-B}r<g_{\by,\bB}(y_i(\beta,\bB))$, which shows that $y_1(\beta,\bB)>y_i(\beta,\bB)$ for any $i \in \{2,\dots, q\}$ as $g_{\by,\bB}$ is decreasing on $(0,y_{\max}]$. In particular, the conclusion of the proposition is true in this case.
\smallskip

\textbf{Case $\min_{x \in (0,y_{\max}]}g_{\by,\bB}(x)\neq g(y_{\max})$.}
In this case, the convex continuous function $g_{\by,\bB}$ has a global minimum $x_*(\by,\bB)$ on $ (0,y_{\max}]$, $g_{\by,\bB}$ is strictly decreasing on $(0,x_*(\beta,\bB))$ and strictly increasing on $[x_*(\beta,\bB),y_{\max}]$. Recalling $r=g_{\by,\bB}(y_2(\beta,\bB))$, this implies that $y_1(\beta,\bB)\in \{z_0,z_1\}$, and $y_2(\beta,\bB),\dots,y_q(\beta,\bB) \in \{x_0,x_1\}$, where $z_0 \in (0,x_*(\by,\bB))$ and $z_1 \in (x_*(\by,\bB),y_{\max}]$ satisfy $g_{\by,\bB}(z_0)=g_{\by,\bB}(z_1)=\e^{-B}r$, while $x_0 \in (0,x_*(\by,\bB))$ and $x_1 \in (x_*(\by,\bB),y_{\max}]$ satisfy $g_{\by,\bB}(x_0)=g_{\by,\bB}(x_1)=r$. Observe that $x_0<z_0<z_1<x_1$, since $\e^{-B}r<r$ and due to the fact that $g_{\by,\bB}$ is strictly decreasing on $(0,x_*(\by,\bB))$ while strictly increasing on $(x_*(\by,\bB),y_{\max}]$.

We claim that in this case $y_1(\beta,\bB)=\max\{y_i(\beta,\bB)\colon i \in [q]\}=y_{\max}$. Before proving this claim, let us finish the proof of the proposition subject to it. Note that since both $z_0,z_1<x_1$, none of the $y_i(\beta,\bB)$ for $i\in \{2,\dots,q\}$ can take the value $x_1$, as otherwise we would have $y_i(\beta,\bB)>y_1(\beta,\bB)=\max\{y_i(\beta,\bB)\colon i \in [q]\}$, which leads to a contradiction. Consequently, $y_2(\beta,\bB)=y_3(\beta,\bB)=\dots=y_q(\beta,\bB)=x_0$, and $\by(\beta,\bB)$ has the form \eqref{eq:soln_form}.
\smallskip

Finally, we check the claim that $y_1(\beta,\bB)=\max\{y_i(\beta,\bB)\colon i \in [q]\}$. We will assume otherwise, and show that this leads to a contradiction. The only way for $y_1(\beta,\bB)$ to fail to be the largest coordinate of $\by(\beta,\bB)$ is to have some $y_i(\beta,\bB)=x_1$ for some $i\in \{2,\dots,q\}$. Without loss of generality, let $y_2(\beta,\bB)=x_1$. Note that since $y_2(\beta,\bB)>y_1(\beta,\bB)$,
\eqan{\label{eq:switching}
F_{\beta,\bB}(y_1,y_2,y_3,\dots,y_q)&=\E_{W}\Big[\log\Big(\e^{\beta'y_1\frac{W}{\E[W]}+B}+\e^{\beta'y_2\frac{W}{\E[W]}}+\sum_{j=3}^q\e^{\beta'y_j\frac{W}{\E[W]}}\Big)\Big]\nn\\&<\E_{W}\Big[\log\Big(\e^{\beta'y_2\frac{W}{\E[W]}+B}+\e^{\beta'y_1\frac{W}{\E[W]}}+\sum_{j=3}^q\e^{\beta'y_j\frac{W}{\E[W]}}\Big)\Big]\\&=F_{\beta,\bB}(y_2,y_1,y_3,\dots,y_q)\nn,    
}
where we have abbreviated $y_i=y_i(\beta,\bB)$. The last display shows that $\by(\beta,\bB)$ is not a global supremum of $F_{\beta,\bB}(\by)-\frac{\beta'}{2\E[W]}\sum_{i=1}^q y_i^2$, which lead to a contradiction. Hence $y_1(\beta,\bB)=\max\{y_i(\beta,\bB):i \in [q]\}$, which proves the claim and hence the proposition in this case.
\end{proof}

We next consider the zero-field case, i.e., when one has a field of the form $\bB=(0,0,\dots,0)$. Interestingly, the proof of Proposition \ref{prop-form-solution-B>0} does not generalize to this case. The reason is that, when $B=0$, we cannot claim that when $y_1(\by,\bB)$ is not the largest coordinate, there has to be one \emph{even larger}, which is the basis of the \emph{switching} inequality \eqref{eq:switching}, where we switched the values of $y_1(\beta,\bB)$ and $y_2(\beta,\bB)$ to obtain a contradiction. Thus we have to provide a different proof for the $B=0$ case.
\smallskip

For this, we provide an extension of the result in Proposition \ref{prop-form-solution-B>0} to general non-negative external fields, under the restriction that they are small enough. This result is interesting in its own right, and also allows us to recover the solution structure for the zero-field case:

\begin{proposition}[Solution structure for small magnetic fields]\label{prop:small_field_sol}
    Consider the optimization problem \eqref{eq:laplace} for a general field of the form $\bB=(B_1,\dots, B_q)$, and let $\by(\beta,\bB):=(y_1(\beta,\bB),\dots,y_q(\beta,\bB))$ be any optimizer. Recall the function $g_{\by, \bB}(x)$ from \eqref{eq:g_y_B}, and let $y_{\max}=y_{\max}(\beta,\bB):=\max\{y_i(\beta,\bB)\colon  i \in [q]\}$. Define $x_*(\by,\bB):=\arg\min\{g_{\by,\bB}(x):x\in(0,y_{\max}]\}$. There exists $\varepsilon>0$, such that, under the assumption $0\leq \max_{i\in[q]}B_i < \varepsilon$, the number of $i$ such that $y_i(\beta,\bB)>x_*(\beta,\bB)$ is at most 1.
\end{proposition}
\begin{proof}
Clearly, the conclusions of Lemma \ref{lem:g_convexity} are true for $g_{\by,\bB}$. Further, note that if $x_*(\by,\bB)=y_{\max}$, then the statement of the proposition is true, so for the rest of this proof, we consider the case that $0<x_*(\by,\bB)<y_{\max}$. Recall from Lemma \ref{lem:g_convexity} that $g_{\by,\bB}$ is strictly decreasing on $(0,x_*(\by,\bB))$ and strictly increasing on $[x_*(\by,\bB),y_{\max}]$.

    Note that the stationarity conditions \eqref{eq:stat0} imply that for $i\neq j$ with $i,j \in [q]$,
    \begin{equation}\label{eq:stat_g}
        \e^{B_i}g_{\by,\bB}(y_i(\beta,\bB))=\e^{B_j}g_{\by,\bB}(y_j(\beta,\bB)). 
    \end{equation}
Observe that none of the $y_i(\beta,\bB)$ can be equal to zero; this is because otherwise all of them have to be equal to zero using \eqref{eq:stat_g} and the fact that $\lim_{x \to 0}g_{\by,\bB}(x)=\infty$. However, this violates \eqref{eq:stat0} which asserts that the sum $\sum_{i=1}^q y_i(\beta,B)=\E[W]>0$. 

Recall from Remark \ref{rem-optim-positive} that $\E[W]\geq y_i(\beta,B)>0$ for all $i$. This implies that if we let $\bB \to \boldsymbol{0}=(0,\dots,0) \in \R^q$, along some sequence $\bB_n$, then using compactness on the interval $[0,\E[W]]$ we can extract a further subsequence along which $\by(\beta, B)=(y_i(\beta,\bB))_{i \in [q]}$ converges coordinate-wise to a vector $\bx(\beta)=(x_1(\beta),\dots,x_q(\beta))$.
\invisible{This means that if we let $\bB \to \boldsymbol{0}$, where we recall that $\boldsymbol{0}=(0,\dots,0)$ is the zero vector in $\R^q$, along some sequence $\bB_n$, then using compactness on the interval $[0,\E[W]]$ we can extract a further subsequence along which $\by(\beta, B)=y_i(\beta,\bB)$ converges coordinate-wise to a vector $\bx(\beta)=}

From now on, unless mentioned otherwise, we work along this convergent subsequence. We suppress the notation of this subsequence, but whenever we write $\bB_n \to \boldsymbol{0}$, we think that it goes to zero along this convergent subsequence.

Observe that letting $\bB=\bB_n$ in the expression \eqref{eq:stat0}, and then $n \to \infty$, by dominated convergence, we get that the vector $(x_1(\beta),\dots,x_q(\beta))$ satisfies the analogous stationarity conditions
\begin{equation}\label{eq:stat_g_0}
    g_{\by,\mathbf{0}}(x_i(\beta))=g_{\by,\mathbf{0}}(x_j(\beta)),
\end{equation}
where by $g_{\by,\mathbf{0}}$ we simply mean the function $g_{\by,\bB}$ with $\bB=\mathbf{0}$. It is easy to check that the conclusions of Lemma \ref{lem:g_convexity} hold true for the function $x \mapsto g_{\by,\mathbf{0}}$, on the domain $(0,x_{\max}]$, where $x_{\max}=x_{\max}(\beta):=\max\{x_i(\beta)\colon  i \in [q]\}$ is the (subsequential) limit of $y_{\max}$ as $\bB \to \mathbf{0}$. In particular we note that none of the $x_i(\beta)$ can be equal to zero, as otherwise, all of them must equal zero using the stationarity conditions \eqref{eq:stat_g_0}, which would then violate $\sum_{i=1}^qx_i(\beta)=\E[W]$.

With all these observations in place, we are ready to go into the heart of the proof. Let us assume without loss of generality that $y_{\max}=y_{\max}(\beta,\bB_n)=y_1(\beta,\bB_n)$, and further assume for a contradiction that $y_1(\beta,\bB_n)>x_*(\by,\bB_n)$, and also some other coordinate, which without loss of generality we can assume to be $y_2(\beta,\bB_n)$, satisfies $y_2(\beta,\bB_n)>x_*(\by,\bB_n)$.

Then we define the function
\begin{align*}
  H_{\by,\bB}(x)&=\E\left[\log(\e^{x\frac{\beta'W}{\E[W]}+B_1}+\e^{(y_1+y_2-x)\frac{\beta'W}{\E[W]}+B_2}+\sum_{i=3}^q\e^{y_i\frac{\beta'W}{\E[W]}+B_i})\right]\\
  &\qquad -\frac{\beta'}{2\E[W]}\Big(x^2+(y_1+y_2-x)^2+\sum_{i=3}^q y_i^2\Big),\nn
\end{align*}
for $x\in (0,y_1+y_2)$, where we abbreviated $y_i=y_i(\beta,\bB)$ and $\bB=\bB_n$ in the above expression. Note that since $\by(\beta,\bB)$ is a maximizer of \eqref{eq:laplace}, $x=y_1$ is a global supremum of the function $H_{\by,\bB}(x)$. This gives the condition that $$\frac{d^2}{dx^2}H_{\by,\bB}(x)\big|_{x=y_1}\leq 0.$$

The second derivative of $H_{\by,\bB}(x)$ equals
\begin{align*}
    &\mathbb{E}\left[\left(\frac{\beta'W}{\E[W]}\right)^2\left(\frac{\e^{x\frac{\beta'W}{\E[W]}+B_1}+\e^{(y_1+y_2-x)\frac{\beta'W}{\E[W]}+B_2}}{\e^{x\frac{\beta'W}{\E[W]}+B_1}+\e^{(y_1+y_2-x)\frac{\beta'W}{\E[W]}+B_2}+\sum_{i=3}^q\e^{y_i\frac{\beta'W}{\E[W]}+B_i}}\right)\right]\\&-\mathbb{E}\left[\left(\frac{\beta'W}{\E[W]}\right)^2\left(\frac{(\e^{x\frac{\beta'W}{\E[W]}+B_1}-\e^{(y_1+y_2-x)\frac{\beta'W}{\E[W]}+B_2})^2}{(\e^{x\frac{\beta'W}{\E[W]}+B_1}+\e^{(y_1+y_2-x)\frac{\beta'W}{\E[W]}+B_2}+\sum_{i=3}^q\e^{y_i\frac{\beta'W}{\E[W]}+B_i})^2}\right)\right]-\frac{2\beta'}{\E[W]}.
\end{align*}
Consequently, the condition $\frac{d^2}{dx^2}H_{\by,\bB}(x)\big|_{x=y_1}\leq 0$ gives us that
\begin{align*}
    &\mathbb{E}\left[\frac{\beta'W^2}{\E[W]}\left(\frac{\e^{y_1\frac{\beta'W}{\E[W]}+B_1}+\e^{y_2\frac{\beta'W}{\E[W]}+B_2}}{\e^{y_1\frac{\beta'W}{\E[W]}+B_1}+\e^{y_2\frac{\beta'W}{\E[W]}+B_2}+\sum_{i=3}^q\e^{y_i\frac{\beta'W}{\E[W]}+B_i}}\right)\right]\\&-\mathbb{E}\left[\frac{\beta'W^2}{\E[W]}\left(\frac{(\e^{y_1\frac{\beta'W}{\E[W]}+B_1}-\e^{y_2\frac{\beta'W}{\E[W]}+B_2})^2}{(\e^{y_1\frac{\beta'W}{\E[W]}+B_1}+\e^{y_2\frac{\beta'W}{\E[W]}+B_2}+\sum_{i=3}^q\e^{y_i\frac{\beta'W}{\E[W]}+B_i})^2}\right)\right]-2\leq 0
\end{align*}
Now, we let $\bB\to 0$ along the convergent subsequence, to obtain that
\begin{align}
\label{H''-inequality}&\mathbb{E}\left[\frac{\beta'W^2}{\E[W]}\left(\frac{\e^{x_1\frac{\beta'W}{\E[W]}}+\e^{x_2\frac{\beta'W}{\E[W]}}}{\e^{x_1\frac{\beta'W}{\E[W]}}+\e^{x_2\frac{\beta'W}{\E[W]}}+\sum_{i=3}^q\e^{x_i\frac{\beta'W}{\E[W]}}}\right)\right]\nn\\&-\mathbb{E}\left[\frac{\beta'W^2}{\E[W]}\left(\frac{(\e^{x_1\frac{\beta'W}{\E[W]}}-\e^{x_2\frac{\beta'W}{\E[W]}})^2}{(\e^{x_1\frac{\beta'W}{\E[W]}}+\e^{x_2\frac{\beta'W}{\E[W]}}+\sum_{i=3}^q\e^{x_i\frac{\beta'W}{\E[W]}})^2}\right)\right]-2 \leq 0,
\end{align}
where as before, we have abbreviated $x_i(\beta)=x_i$ for $i\in[q]$.

Again by a compactness argument, and passing to a further (sub)subsequence, we have that $x_*(\by,\bB)=\text{argmin}\{g_{\by,\bB}(x):x\in(0,y_{\max}]\}$ converges to $x_*=x_*(\beta):=\text{argmin}\{g_{\by,\mathbf{0}}(x):x\in(0,x_{\max}]\}$. Furthermore, as (non-strict) inequalities are preserved after passing on to limits, we have that both $x_1(\beta),x_2(\beta)\geq x_*$. Note that then by injectivity of the function $g_{\by,\mathbf{0}}$ on $[x_*,x_{\max}]$, and using the stationary conditions \eqref{eq:stat_g_0} for $i=1$ and $j=2$, we must have that $x_1(\beta)=x_2(\beta)$. In particular, \eqref{H''-inequality} reduces to
\begin{align*}
\mathbb{E}\left[\frac{\beta'W^2}{\E[W]}\left(\frac{\e^{x_1\frac{\beta'W}{\E[W]}}}{\sum_{i=1}^q\e^{x_i\frac{\beta'W}{\E[W]}}}\right)\right]-1 \leq 0.
\end{align*}

Going back to the pre-limit, we see that for any given $\delta>0$, one can choose $\varepsilon>0$ such that whenever $\max_{i\in[q]}{B_i}<\varepsilon$ (along the convergent subsequence),
\begin{equation}\label{eq:sol_form_cont_1}
    \mathbb{E}\left[\frac{\beta'W^2}{\E[W]}\left(\frac{\e^{y_1(\beta,\bB)\frac{\beta'W}{\E[W]}}}{\sum_{i=1}^q\e^{y_i(\beta,\bB)\frac{\beta'W}{\E[W]}}}\right)\right]-1 < \delta.
\end{equation}
However, we note that as $y_1(\beta,\bB)>x_*(\by,\bB)$, the function $g_{\by,\bB}$ is strictly increasing at $y_1=y_1(\beta,\bB)$, and thus $\frac{d}{dx}g_{\by,\bB}(x)\big|_{x=y_1}>0$. Writing this out and using the stationarity conditions \eqref{eq:stat_g}, we arrive at
\begin{equation*}
    \mathbb{E}\left[\frac{\beta'W^2}{\E[W]}\left(\frac{\e^{y_1(\beta,\bB)\frac{\beta'W}{\E[W]}}}{\sum_{i=1}^q\e^{y_i(\beta,\bB)\frac{\beta'W}{\E[W]}}}\right)\right]-1>0.
\end{equation*}
In particular, the left-hand side in the last display above is strictly larger than $2\delta$, for some sufficiently small $\delta>0$, which contradicts \eqref{eq:sol_form_cont_1}. Hence, this leads to a contradiction with our assumption that there are two or more coordinates of $\by(\beta,\bB)$ that are larger than $x_*(\by,\bB)$, which finishes the proof.
\end{proof}

\begin{remark}[Recovering the solution structure for the one non-negative field case]\label{rem:general_sol_form_nnf}
{\rm 
    When all the field values $B_i=0$, the last result tells us that at most one coordinate of the optimizer $\by(\beta,\mathbf{0})=(y_1(\beta,\mathbf{0}),\dots,y_q(\beta,\mathbf{0}))$ can be larger than $\text{argmin}(g_{\by,\mathbf{0}})=x_*$. On the other hand, the value of the function $g_{\by,\mathbf{0}}$ at all these coordinates are equal, thanks to the stationarity conditions $g_{\by,\mathbf{0}}(y_i(\beta,\mathbf{0}))=g_{\by,\mathbf{0}}(y_j(\beta,\mathbf{0}))$ for any $i,j \in [q]$. Consequently, by injectivity of the function $g_{\by,\mathbf{0}}$ on the interval $(0,x_*)$, and the fact that at least $q-1$ many coordinates $y_i(\beta,\mathbf{0})$ are in $(0,x_*)$, we obtain that the optimizer $\by(\beta,\mathbf{0})$ has at least $q-1$ equal coordinates, i.e., it has the form \eqref{eq:soln_form}. Also, when the field is of the form $(B,0,\dots,0)$ for some $B>0$, Proposition \ref{prop-form-solution-B>0} tells us that the optimizer is of the form \eqref{eq:soln_form}. Combining, we get that the optimizer has the form \eqref{eq:soln_form} when the field is of the form $(B,0,\dots,0)$ for any $B\geq 0$.}\hfill\ensymboldefinition
\end{remark}

\medskip

Finally, we use the observation in Remark \ref{rem:general_sol_form_nnf} to complete the proof of Theorem \ref{thm-soln_form-positive}:
\medskip

\noindent
{\it Proof of Theorem \ref{thm-soln_form-positive}.}
By Remark \ref{rem:general_sol_form_nnf} the maximum in \eqref{eq:laplace} is indeed of the form \eqref{eq:soln_form}, and we are led to 
\be
\label{varphi-form-B=0}
\varphi(\beta,B) = \sup_{s} p_{\beta,B} (s)
\ee
where
\eqan{\label{eq:pressure_evaluated_at_soln}
p_{\beta,B}(s) &= \mathbb{E}\left[ \log\left(\e^{\frac{\beta'}{q}W(1+(q-1)s)+B} + (q-1) \e^{\frac{\beta'}{q}W(1-s)} \right) \right] - \frac{\beta' \expec[W]} {2 q}\left[ (q-1)s^2 
+ 1\right]\nn\\
&=\mathbb{E}\left[ \log\left(\e^{\frac{\beta'}{q}W(q-1)s+B} + (q-1) \e^{-\frac{\beta'}{q}Ws} \right) \right]- \frac{\beta' \expec[W]} {2 q}\left[ (q-1)s^2 
- 1\right]\nn\\
&=\mathbb{E}\left[ \log\left(\e^{\beta'sW+B} + (q-1) \right) \right]-\frac{s\beta' \expec[W]} {q}- \frac{\beta' \expec[W]} {2 q}\left[ (q-1)s^2 
-1\right]\nn\\
&=\mathbb{E}\left[ \log\left(\e^{\beta'sW+B} + (q-1) \right) \right]- \frac{\beta' \expec[W]} {2 q}\left[ (q-1)s^2 
+2s-1\right].
}
The optimal value of $s$ is determined from $p_{\beta,B}'(s) =0$, for which we compute
    \eqan{
    \label{p-beta-B-derivative}
    p_{\beta,B}'(s)&=
    \mathbb{E}\left[\beta' W \frac{\e^{\beta' W s+B}}{\e^{\beta' W s+B}+q-1}\right]
    - \frac{\beta' \expec[W]} { q}\left[ (q-1)s 
+1\right]\nn\\
    &=\frac{q-1}{q}\mathbb{E}\left[\beta' W \frac{\e^{\beta' W s+B}-1}{\e^{\beta' W s}+q-1}\right]
    - \frac{\beta' \expec[W](q-1)} { q}s\nn\\
    &=\frac{q-1}{q}\beta'\expec[W]\big[\mathscr{F}_B(\beta's)-s\big].
    }
Therefore, $s(\beta,B)$ satisfies the stationarity condition in \eqref{FOPT-1-B1-B>0}.
This completes the proof of Theorem \ref{thm-soln_form-positive}.
\qed

\section{Spontaneous magnetization for infinite-variance weights: Proof of Theorem \ref{thm-spont-magnetization}}
\label{sec-tau(2,3)}




Let $s(\beta,0)$ be any optimizer of \eqref{var-problem-psi-B>0}. To prove Theorem \ref{thm-spont-magnetization}, we have to show that when $\E[W^2]=\infty$, for any $\beta>0$, any argmax $s(\beta,0)$ of the function $p_{\beta,0}(s),$ for $s \in [0,1]$, is positive. By \eqref{p-beta-B-derivative},
\begin{align*}
    p_{\beta, 0}'(s)=\frac{q-1}{q}\beta'\expec[W]\mathscr{G}_{\beta}(s),
\end{align*}
where
\begin{align*}
    \mathscr{G}_{\beta}(s)=\expec\Big[\frac{W}{\expec[W]}\frac{\e^{\beta'Ws}-1}{\e^{\beta'Ws}+q-1}
    \Big]-s. 
\end{align*}

Note that the signs of $p_{\beta,0}''(s)$ and $\mathscr{G}_{\beta}'(s)$ are the same at any value of $s$. We further compute
\begin{align*}
    \mathscr{G}_{\beta}'(s)=-1+\frac{q\beta'}{\expec[W]}\mathbb{E}\left[W^2\left(\e^{(\beta'Ws)/2}+\frac{q-1}{\e^{(\beta'Ws)/2}} \right)^{-2}\right].
\end{align*}
Note that $\mathscr{G}_{\beta}'(s) \to \infty$ as $s \searrow 0$ for any $\beta>0$, since $\expec[W^2]=\infty$. Furthermore, $\mathscr{G}_{\beta}'(s)< \infty$ for any fixed $s>0$. Hence, by continuity, $\mathscr{G}_{\beta}'(s)$ and therefore $p_{\beta,0}''(s)$, is strictly positive in an interval of the form $(0, \varepsilon)$ for some $\varepsilon>0$. Note that this in particular implies that no maximizer of the function $p_{\beta,0}(s)$ in $[0,1]$ can lie in $(0,\varepsilon)$. As a consequence, for any $\beta>0$, and for any optimizer $s(\beta,0)$, we have $s(\beta,0)>\varepsilon>0$. \qed

\section{On the nature of the phase transition: Proof of Theorems \ref{thm-nature-pt-a}--\ref{thm-nature-pt-b}}
\label{sec-nature-pt}

In this section, we prove Theorems \ref{thm-nature-pt-a} and \ref{thm-nature-pt-b}. Thus we assume throughout this section that $\mathbb{E}[W^2] < \infty$, $q\ge 3$, $0\le B < \log(q-1)$ and that Condition \ref{cond-zer-cross} holds. For the proof, we proceed in three steps: in the first step in Section \ref{sec:stationary-points}, we analyze the solutions of the stationarity condition \eqref{FOPT-1-B1-B>0}, which
we recall here for reader's convenience
to read
\begin{align}
\label{eq:stat-cond-B0}
  \expec \left [\frac{W}{\expec[W]}\frac{\e^{\beta'Ws+B}-1}{\e^{\beta'Ws+B}+q-1}
    \right ]=s. 
\end{align}
In the second step in Section \ref{sec:optimal-solution}, we identify which of these solutions is optimal for the variational problem \eqref{var-problem-psi-B>0}-\eqref{eq:pressure_evaluated_at_soln-a-B>0}, i.e.,
\eqan{\label{eq:pressure_evaluated_at_soln-a-B=0}
    \varphi(\beta,B) = \sup_{s\in[0,1]} p_{\beta,B}(s),
    }
    where
\eqan{\label{eq:pressure-p-beta-zero}    
p_{\beta,B}(s) =\mathbb{E}\left[ \log\left(\e^{\beta'W s+B} + q-1 \right) \right]- \frac{\beta' \expec[W]} {2 q}\left[ (q-1)s^2 
+2s-1\right].
}
We show the existence of a critical value $0<\beta_c(B)<\infty$ such that when $\beta<\beta_c(B)$ the optimal solution is $s_1(\beta,B)$ (the smallest solution of  \eqref{eq:stat-cond-B0}), while for $\beta>\beta_c(B)$ the optimal solution is $s_3(\beta,B)$ (the largest solution of \eqref{eq:stat-cond-B0}). Furthermore the change of the optimal solution occurs via a discontinuity at $\beta_c(B)$. In the last step in Section \ref{sec:conclusion-proof}, we conclude and prove the statements
of Theorems \ref{thm-nature-pt-a} and \ref{thm-nature-pt-b}.

\subsection{Stationary points: Graphical analysis}
\label{sec:stationary-points}
It is convenient to change variables, and define $t = \beta' s$. Then,  \eqref{eq:stat-cond-B0} for the $s$ variable 
is rewritten as an equation for the $t$ variable as
\be
\label{eq:stat-cond-t}
{\mathscr{F}}_B(t) = \frac{t}{\beta'},
\ee
where ${\mathscr{F}}_B$ is the function defined \eqref{def-Fcal}, which does not depend on $\beta$.
We can then use a graphical approach to study the intersections between the function ${\mathscr{F}}_B$ and the lines with slope $1/\beta'$ passing through the origin. For the benefit of the reader, we include in Figure \ref{fig:second-deriv} a plot of the function
${\mathscr{F}}_0$ and in Figure \ref{fig:second-deriv-b} a plot of the function
${\mathscr{F}}_B$ for a small positive $B$. 

Clearly, we have ${\mathscr{F}_B(0)}=\frac{\e^{B}-1}{\e^{B}+q-1}$ and $\lim_{t\to\infty}{\mathscr{F}_B}(t) = 1$. 
Further,
\be
\label{eq:derivative-F}
\frac{d}{dt}\mathscr{F}_B(t)=\expec\left[\frac{W^2}{\expec[W]}\frac{q\e^{tW+B}}{(\e^{tW+B}+q-1)^2}\right]>0,
\ee
which means that $\mathscr{F}_B$ is increasing. Furthermore,
    \be\label{eq:second_derivative_gen-copy}
    \frac{d^2}{dt^2}\mathscr{F}_B(t)=q\expec\left[\frac{W^3}{\expec[W]}\frac{\e^{tW+B}(q-1-\e^{tW+B})}{(\e^{tW+B}+q-1)^3}\right],
    \ee
and Condition \ref{cond-zer-cross}(a) implies that $\mathscr{F}_B$ has a unique inflection point  $t_*(B) \in (0,\infty)$. Therefore, 
$\mathscr{F}_B(t)$ is strictly convex for $t\in[0,t_*(B))$ and strictly concave for $t\in(t_*(B),\infty)$. The assumption in Theorems \ref{thm-nature-pt-a}  and \ref{thm-nature-pt-b} that $0\le B< \log(q-1)$ is required so that $\frac{d^2}{dt^2}\mathscr{F}_B(0) >0$.
 Moreover, by Condition \ref{cond-zer-cross}(b), the inflection point is  steep, meaning that at this point the slope of $\mathscr{F}_B$ is larger
    than $\mathscr{F}_B(t_*(B))/t_*(B)$.
\smallskip

With these preliminaries, we show that the function
${\mathscr{F}_B}$ has two tangent lines passing through the origin:

\begin{proposition}[Tangent lines]\label{prop:tangents}
The equation $\mathscr{F}_B(t)=[\frac{d}{dt}\mathscr{F}_B(t)]t$ has two solutions,  $0\le t_a(B) < t_*(B) < t_b(B)$. In the following, we abbreviate $t_a = t_a(B)$, $t_b = t_b(B)$ and $t_*=t_*(B)$. As a consequence, there are two tangent lines to the curve $\mathscr{F}_B$ which also pass through the origin. We denote these tangent lines as
$G_a(t)=(\frac{d}{dt}\mathscr{F}_B(t_a))t$, which is the tangent at $(t_a,\mathscr{F}_B(t_a))$, and $G_b(t)=(\frac{d}{dt}\mathscr{F}_B(t_b))t$, which is the tangent at $(t_b,\mathscr{F}_B(t_b))$.     
\end{proposition}

\begin{proof} Let us consider the intercept $q_B(t)=\calF_B(t)-t\frac{d}{dt}\mathscr{F}_B(t)$ of the tangent to $\calF_B$, and search for solutions of the equation $q_B(t)=0$, so that the tangent passes through the origin.
We compute the derivative of $q_B(t)$ and find
$$\frac{dq_B}{dt}(t) = -t \frac{d^2\mathscr{F}_B}{dt^2}(t).$$ 
Condition \ref{cond-zer-cross} then implies that $q_B(t)$ is strictly decreasing in $(0,t_*)$ and strictly increasing in $(t_*, \infty)$. We observe that $q_B(0)={\mathscr{F}_B}(0) \in [0,\frac{(q-2)}{2(q-1)})$ when $0\le B <\log(q-1)$. The steepnees property of $t_*$ in Condition \ref{cond-zer-cross} implies $q(t_*)= \mathscr{F}_B(t_*) - t_* \frac{d{\mathscr{F}}_B}{dt}(t_*)<0$. Therefore, by continuity, there exists a unique $t_a \in [0,t_*)$ such that $q_B(t_a)=0$.
This solution gives the tangent $G_a(t)=(\frac{d}{dt}\mathscr{F}(t_a))t$.
Similarly, since $q_B(t_*)<0$ and $\lim_{t\to\infty}q_B(t)=1$,
by continuity there exists a unique $t_b \in (t_*,\infty)$ such that $q_B(t_b)=0$. This solution gives the tangent $G_b(t)=(\frac{d}{dt}\mathscr{F}(t_b))t$.
\end{proof}

We now investigate the solutions of the \eqref{eq:stat-cond-t} as we vary $\beta$. We split the analysis into three parts, corresponding to the three regions defined by the two tangent lines $G_a(t)$ and $G_b(t)$ with $t\ge0$. 

\begin{proposition}[Unique solution for straight lines above $ G_b$]
\label{prop:sol_lim_vec_below_CR}
    For any $0\le \beta' < \frac{1}{\frac{d\mathscr{F}_B}{dt}(t_b)}$, 
     \eqref{eq:stat-cond-t} has only one solution $0\le t_1(\beta,B) < t_a$. 
\end{proposition}

\begin{proof} Since $t\mapsto \calF_B(t)$ is
strictly concave for $t\in(t_*, \infty)$, on this interval, its graph lies below the line $G_b$ which is tangent at the point $(t_b, \calF(t_b))$ with $t_b>t_*$. Therefore, when $t\in(t_*, \infty)$, there are no solutions to \eqref{eq:stat-cond-t}.
On the other hand, since $0\le \mathscr{F}_B(0) <1$ and the function $\calF_B(t)$ is
increasing and strictly convex for $t\in[0, t_*)$ with $\frac{d\mathscr{F}_B}{dt}(t_*) > \frac{\mathscr{F}_B(t_*)}{t_*}$ by Condition \ref{cond-zer-cross}, then there exists  a unique non-negative intersection between the graph of $\calF_B$ and the line $t/\beta'$ with slope 
$\frac{1}{\beta'} >\frac{d\mathscr{F}_B}{dt}(t_b)$.
\end{proof}

When $\beta' = \frac{1}{\frac{d\mathscr{F}_B}{dt}(t_b)}$, 
\eqref{eq:stat-cond-t} becomes 
${\mathscr{F}}_B(t) = G_b(t)$ and we have 
two solutions, namely $t_1(\frac{1}{\frac{d\mathscr{F}_B}{dt}(t_b)},B)\in[0,t_a)$
and  $t_b$. 
If we continue to increase $\beta'$, then we enter a region with three solutions, as the following proposition shows:

\begin{proposition}[Three solutions for straight lines between  $G_a$ and  $G_b$]
\label{prop:soln_inside_h-bis}
For any $\beta'\in (\frac{1}{\frac{d\mathscr{F}_B}{dt}(t_b)}, 
\frac{1}{\frac{d\mathscr{F}_B}{dt}(t_a)})$, \eqref{eq:stat-cond-t}  has three solutions $t_1(\beta,B)<t_2(\beta,B)<t_3(\beta,B)$ that satisfy
$t_1\in[0,t_a)$, 
$t_2\in(t_a,t_b)$,
$t_3\in(t_b,\infty)$.
\end{proposition}

\begin{proof}
We define $\Delta (t)= t/\beta' -\calF_B(t)$ and
separately consider three sub-cases:

\smallskip
$\bullet$ {\bf Case $\frac{1}{\frac{d\mathscr{F}_B}{dt}(t_b)}< \beta' < \frac{t_*}{\calF_B(t_*)}$}:
in this case the line $t/\beta'$ is above the line $\frac{\calF_B(t_*)}{t_*}t$ for all $t\ge0$. As before, there exists a unique solution (call it $t_1(\beta,B)$) of \eqref{eq:stat-cond-t} for $t\in[0,t_*)$. Further, this solution must be smaller than the tangent point $t_a$, i.e., 
$0\le t_1(\beta,B) <t_a$. 
For $t>t_*$, we reason as follows.
Since $\frac{d\calF_B}{dt}(t)$ is decreasing for $t>t_*$, then we have that
$\frac{d \Delta}{dt}(t)
=\frac{1}{\beta'}-\frac{d\calF_B}{dt}(t)
<\frac{1}{\beta'}-\frac{d\calF_B}{dt}(t_b)
<0$ for all $t\in(t_*,t_b)$. As $\Delta(t_*)>0$ and $\Delta(t_b)<0$
we conclude that there exists a unique solution (call it $t_2(\beta,B)$)
of \eqref{eq:stat-cond-t} for $t\in(t_*,t_b)$.
Finally by concavity and boundedness of $\calF_B(t)$,
another unique solution (call it $t_3(\beta,B)$) must exist
for $t>t_b$.

\smallskip
$\bullet$ {\bf Case $\beta' = \frac{t_*}{\calF_B(t_*)}$}: obviously  the point $t_2(\beta,B)=t_*$ is a solution. As in the previous case, two additional solutions exist: $t_1(\beta,B)\in[0,t_a)$ and  $t_3(\beta,B)\in[t_b,\infty)$.

\smallskip
$\bullet$ {\bf Case $\frac{t_*} {\calF_B(t_*)} < \beta' < \frac{1}{\frac{d\mathscr{F}_B}{dt}(t_a)}$}: in this case we have
one unique solution $t_1(\beta,B)\in[0,t_a)$. 
Further, another unique solution of \eqref{eq:stat-cond-t} exists for $t\in (t_a,t_*)$. Indeed,
$\frac{d \Delta}{dt}(t)
=\frac{1}{\beta'}-\frac{d\calF_B}{dt}(t)
<\frac{1}{\beta'}-\frac{d\calF_B}{dt}(t_*)$ since $\frac{d\calF_B}{dt}(t)$ is increasing in this interval.
By Condition \ref{cond-zer-cross},
$\frac{d \Delta}{dt}(t)<0$. As $\Delta(t_a)>0$ and $\Delta(t_*)<0$,
we conclude that there exists a unique solution (call it $t_2(\beta,B)$)
of \eqref{eq:stat-cond-t} when $t\in(t_a,t_*)$.
Finally, a third unique solution $t_3(\beta,B)$ exists for 
$t\in (t_*,\infty)$.
\end{proof}

The solution $t_1(\beta,B)$ of Proposition \ref{prop:soln_inside_h-bis}
is increasing in $\beta'$, while the solution $t_2(\beta,B)$ is deccreasing in $\beta'$. When 
$\beta' = \frac{1}{\frac{d\mathscr{F}_B}{dt}(t_a)}$, 
\eqref{eq:stat-cond-t} becomes 
${\mathscr{F}}_B(t) = G_a(t)$ and we have 
two solutions, namely $t_1(\frac{1}{\frac{d\mathscr{F}_B}{dt}(t_a)},B)=t_2(\frac{1}{\frac{d\mathscr{F}_B}{dt}(t_a)},B)=t_a$
and  $t_3(\frac{1}{\frac{d\mathscr{F}_B}{dt}(t_a)},B)\in (t_*,\infty)$. 
If we continue increasing $\beta'$ beyond this point, we enter a region with two solutions when $B=0$ or one solution when $B>0$: 
\begin{proposition}[One/Two solutions for straight lines below or equal to $G_a$]\label{prop:soln_above_h}
    For any $\beta'\geq \frac{1}{\frac{d\mathscr{F}_B}{dt}(t_a)}$,
    \eqref{eq:stat-cond-t}  has  two solutions $0=t_1(\beta',0)<t_3(\beta',0) < \infty$ when $B=0$, and
    one solution $t_3(\beta,B)$ when $B>0$.
\end{proposition}
\begin{proof}
    Assume $B=0$. Then it is easy to see that $t_a=0$ and $t_1(\beta',0)=0$ is a solution. We need to prove the existence and uniqueness of the solution $t_3(\beta',0)$. If $\beta'\geq \frac{1}{\frac{d\mathscr{F}_0}{dt}(0)} = \frac{q\expec[W]}{\expec[W^2]}$ then the line $t/\beta'$ does not intersect ${\mathscr{F}_0}$ for $t\in(0,t_*)$, since its slope is less than the slope of $\mathscr{F}_0$ at $t=0$, and $\mathscr{F}_0$ is strictly convex in $(0,t_*)$.  For $t > t_*$ the line $t/\beta'$ is bound to intersect ${\mathscr{F}_0}$ exactly one time since
    ${\mathscr{F}_0}(t_*) > \frac{t_*}{\beta'}$, $\mathscr{F}_0$ is strictly concave and $\mathscr{F}_0(t)<1$. The proof of the case when $B>0$ is similar and left to the reader. 
\end{proof}

\subsection{Optimal solution}
\label{sec:optimal-solution}

After identifying all stationary points, the next step is to determine which is the optimal solution. Calling $s_i(\beta,B) = t_i(\beta,B)/\beta'$ for $i\in\{1,2,3\}$, where the $t_i(\beta,B)$ have been described in the previous section, the following lemma establishes the local nature of the stationary points of $p_{\beta,B}$:

\begin{lemma}[Nature of stationary points] 
\label{lemma:nature-stat-points-bis} 
As $\beta'$ varies in the interval $[0,\infty)$, the following holds:
\begin{enumerate}
\item[\text{(i)}] 
For any $\beta' < \frac{1}{\frac{d\mathscr{F}_B}{dt}(t_b)}$, the point  $s_1(\beta,B)=0$ is a maximum of $p_{\beta,B}$;
\item[\text{(ii)}] 
For $\beta'\in (\frac{1}{\frac{d\mathscr{F}_B}{dt}(t_b)}, 
\frac{1}{\frac{d\mathscr{F}_B}{dt}(t_a)})$,
among the points  $s_1(\beta,B)<s_2(\beta,B) <s_3(\beta,B)$, the first and the last are local maxima of $p_{\beta,B}$, while $s_2(\beta,B)$ is a minimum;
\item[\text{(iii)}] 
For $\beta' > \frac{1}{\frac{d\mathscr{F}_B}{dt}(t_a)}$, if $B=0$, the point $s_1(\beta',0)=0$ is a minimum of $p_{\beta',0}$ and $s_3(\beta',0)>0$ is a maximum. For $B>0$, the point  $s_3(\beta,B)>0$ is a maximum of $p_{\beta,B}$.
\end{enumerate}
\end{lemma}

\begin{proof}
By \eqref{p-beta-B-derivative},
\begin{eqnarray}
\frac{d^2p_{\beta,B}}{ds^2}(s) & = &
\beta'\frac{(q-1)}{q} \mathbb{E}[W]\left(
\beta'\frac{d{\mathscr{F}_B}}{ds}(\beta's) - 1
\right). 
\end{eqnarray}
In particular,
\begin{eqnarray}
\frac{d^2p_{\beta,B}}{ds^2}(s_1(\beta,B)) & = &
\beta'\frac{(q-1)}{q} \mathbb{E}[W]\left(
\beta'\frac{d{\mathscr{F}_B}}{dt}(t_1(\beta,B)) - 1
\right), 
\end{eqnarray}
which implies that $s_1(\beta,B)$ is a (local) maximum for $\beta' <\frac{1}{\frac{d{\mathscr{F}_B}}{dt}(t_1(\beta,B))}$, and a minimum otherwise. 
Because of convexity, the derivative of $\calF_B(t)$
for $t\in[0,t_a)$ is strictly smaller than $\frac{d{\mathscr{F}_B}}{dt}(t_a)$, which is the slope
of the tangent $G_a$. This proves that the point $s_1(\beta,B)$ is a maximum when $\beta' < \frac{1}{\frac{d\mathscr{F}_B}{dt}(t_a)}$ for all $0\le B<\log(q-1)$. If $\beta' > \frac{1}{\frac{d\mathscr{F}_B}{dt}(t_a)}$ then $s_1(\beta,B)=0$ exists only when $B=0$,
and is a minimum. 

As a consequence, since $s_1(\beta,B) < s_2(\beta,B) < s_3(\beta,B)$ by item (ii), it follows that
$s_2(\beta,B)$ is a minimum and $s_3(\beta,B)$ is a maximum. Similarly, since $s_1(\beta',0) < s_3(\beta',0)$ by item (iii), it follows that
$s_3(\beta',0)$ is a maximum.
\end{proof}

Lemma \ref{lemma:nature-stat-points-bis} tells us that
$s_1(\beta,B)$ is the best solution for any $\beta' < \frac{1}{\frac{d\mathscr{F}_B}{dt}(t_b)}$, whereas
$s_3(\beta,B)$ is best for $\beta'> \frac{1}{\frac{d\mathscr{F}_B}{dt}(t_a)}$. In the following proposition we establish the optimal value
in the region with three solutions:

\begin{proposition}[``Equal area'' construction]
\label{prop:optimal-value}
There exists $\beta'_c(B) \in \big(\frac{1}{\frac{d\mathscr{F}_B}{dt}(t_b)}, \frac{1}{\frac{d\mathscr{F}_B}{dt}(t_a)}\big)$ such that when $0\le\beta'<\beta'_c(B)$  the optimal value in \eqref{eq:pressure_evaluated_at_soln-a-B=0}
is $s_1(\beta,B)$ and when $\beta'>\beta'_c(B)$ the  optimal value is
$s_3(\beta,B)$.
\end{proposition}
\begin{proof}
Since we  have the results of Lemma \ref{lemma:nature-stat-points-bis}, we only need to 
compare $p_{\beta,B}(s_1(\beta,B))$ to $p_{\beta,B}(s_3(\beta,B))$ for 
$\beta' \in \big(\frac{1}{\frac{d\mathscr{F}_B}{dt}(t_b)}, \frac{1}{\frac{d\mathscr{F}_B}{dt}(t_a)}\big)$.
To this aim, we define 
\begin{align}
R_B(\beta')=\int_{t_1(\beta,B)}^{t_3(\beta,B)} \Big(\mathscr{F}_B(t)-\frac{t}{\beta'}\Big)dt,\label{eq:def_R}   
\end{align}
and observe that
\eqan{
&p_{\beta,B}(s_3(\beta,B)) - p_{\beta,B}(s_1(\beta,B))\nonumber\\
&\quad =
\beta' \frac{(q-1)}{q} \mathbb{E}[W]
\int_{s_1(\beta,B)}^{s_3(\beta,B)}  \left(\expec\Big[\frac{W}{\expec[W]}\frac{\e^{\beta'Ws'+B}-1}{\e^{\beta'Ws'+B}+q-1}
    \Big]-s'\right) ds'\nonumber\\
    &\quad=
\frac{(q-1)}{q} \mathbb{E}[W] 
R_B(\beta').
}
We now prove that the mapping
$\beta'\mapsto R_B(\beta')$ is monotonically increasing when $\beta'$ is in the compact set 
 $\big[\frac{1}{\frac{d\mathscr{F}_B}{dt}(t_b)}, \frac{1}{\frac{d\mathscr{F}_B}{dt}(t_a)}\big]$.
Furthermore, we prove that 
$R_B({\frac{1}{\frac{d\mathscr{F}_B}{dt}(t_b)}})<0$ 
and 
$R_B({\frac{1}{\frac{d\mathscr{F}_B}{dt}(t_a)}})>0$. Therefore, by the intermediate value theorem, there exists a unique 
$\beta'_c(B) \in 
(\frac{1}{\frac{d\mathscr{F}_B}{dt}(t_b)}, \frac{1}{\frac{d\mathscr{F}_B}{dt}(t_a)})$
such that $R_B(\beta_c')=0$
and the claim of the proposition follows.

To prove monotonicity, we compute
\begin{eqnarray}
\label{eq:derivative-R-bis}
\frac{d R_B}{d \beta'}(\beta')
&=&
\left (\mathscr{F}_B(t_3(\beta,B))-\frac{t_3(\beta,B)}{\beta'}\right ) \frac{d t_3}{d\beta'}(\beta,B)
\nonumber\\
&&\quad -
\left (\mathscr{F}_B(t_1(\beta,B))-\frac{t_1(\beta,B)}{\beta'}\right ) \frac{d t_1}{d\beta'}(\beta,B)
+
\int_{t_1(\beta,B)}^{t_3(\beta,B)} \frac{t}{{\beta'}^2}\, dt  
\nonumber\\
&=&
\frac{1}{2 {\beta'}^2}\Big(t_3(\beta,B)^2-t_1(\beta,B)^2\Big)>0,
\end{eqnarray}
where we observe that the two terms inside the round brackets are zero because $t_1$ and $t_3$ are stationary points,
and the differentiability of $t_1$ and $t_3$ is guaranteed 
by the implicit function theorem.
Thus \eqref{eq:derivative-R-bis} shows that $R_B(\beta')$ is strictly increasing for $\beta'\in 
\big[\frac{1}{\frac{d\mathscr{F}_B}{dt}(t_b)}, \frac{1}{\frac{d\mathscr{F}_B}{dt}(t_a)}\big]$.

To see that $R_B\big({\frac{1}{\frac{d\mathscr{F}_B}{dt}(t_b)}}\big)<0$, use that $t_3\big({\frac{1}{\frac{d\mathscr{F}_B}{dt}(t_b)}},B\big)=t_b$ and recall that the tangent line $G_b(t)=t\frac{d\mathscr{F}_B}{dt}(t_b)$ completely lies above the graph of $\mathscr{F}_B$ when $t\in[t_1(\beta,B),t_b]$, with the only points of intersection being $(t_b,\mathscr{F}_B(t_b))$ (and  also $(0,0)$ when $B=0$). Therefore $\mathscr{F}_B(t)-t\frac{d\mathscr{F}_B}{dt}(t_b) <0$ for all $t\in(t_1,t_b)$.

To see that 
$R_B\big({\frac{1}{\frac{d\mathscr{F}_B}{dt}(t_a)}}\big)>0$
the argument is similar. 
Use that $t_1\big({\frac{1}{\frac{d\mathscr{F}_B}{dt}(t_a)}},B\big)=t_a$, and recall that
the tangent line $G_a(t)=t\frac{d\mathscr{F}_B}{dt}(t_a)$ completely lies below the graph of $\mathscr{F}_B$. Therefore $\mathscr{F}_B(t)-t\frac{d\mathscr{F}_B}{dt}(t_a) >0$ for all $t\in(t_a,t_3)$.
\end{proof}
\begin{figure}
\centering
\includegraphics[width=10cm]{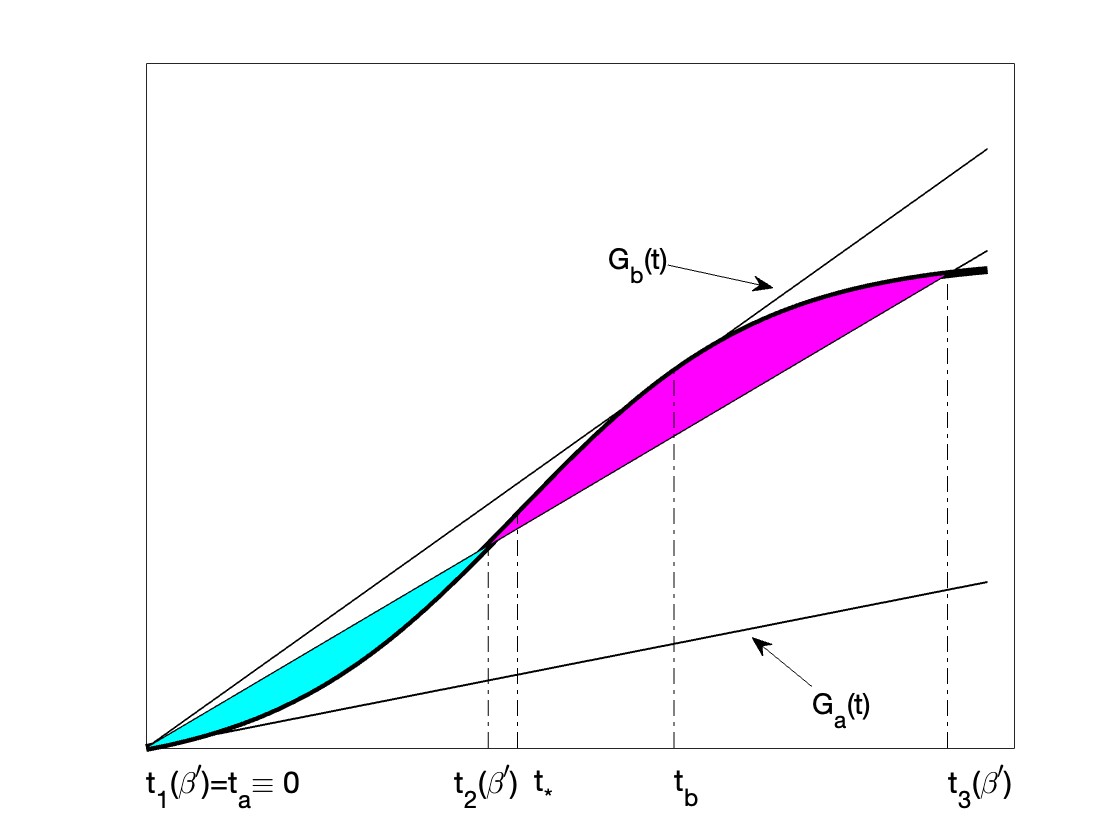}
\caption{The function $t\mapsto \mathscr{F}_0(t)$ for
  the sparse Erd\H{o}s–R\`enyi model with $\lambda=5$ and $q=21$.
  The two lines $G_a(t)$ tangent in $t_a=0$ and $G_b(t)$ tangent in $t_b$, as well as the inflection point $t_*$, are also shown. 
  For $\beta\in (\beta_0,\beta_1)$ the three solutions of \eqref{eq:stat-cond-t} are  represented.
  }
  \label{fig:second-deriv}
\end{figure}


\begin{figure}
\centering
\includegraphics[width=11cm]{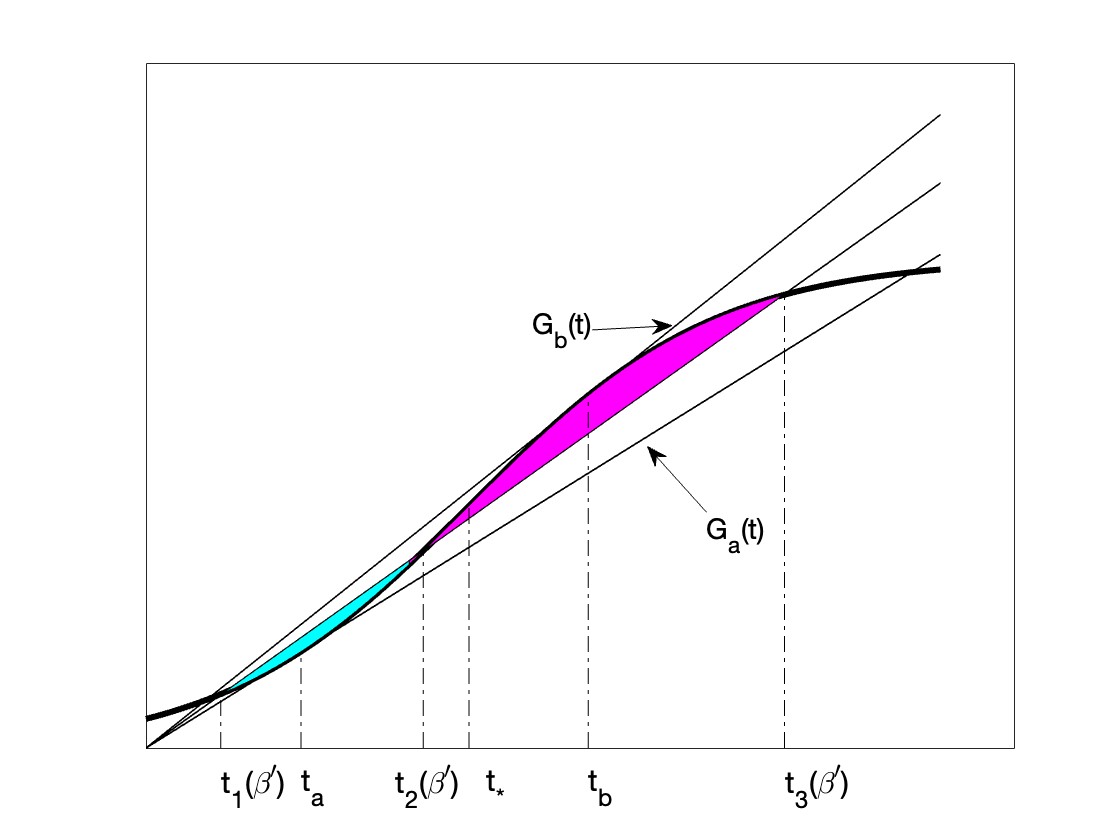}
\caption{The function $t\mapsto \mathscr{F}_B(t)$ for
  the sparse Erd\H{o}s–R\'enyi model with $\lambda=5$, $q=101$ and $B=2$.
  The two  lines $G_a(t)$ tangent in $t_a$ and $G_b(t)$ tangent in $t_b$, as well as the inflection point $t_*$, are also shown. 
  For $\beta\in (\beta_0,\beta_1)$ the three solutions of \eqref{eq:stat-cond-t} are  represented.
}
  \label{fig:second-deriv-b}
\end{figure}

%

\subsection{Conclusion of the proof of Theorems \ref{thm-nature-pt-a} and \ref{thm-nature-pt-b}}
\label{sec:conclusion-proof}

\invisible{{\color{red} Adapt the proof to B>0!}
First we prove that $\beta_c$ is well defined, i.e.
if $s(\beta_1)>0$ and $\beta_2>\beta_1$, then also $s(\beta_2)>0$. From Proposition \ref{prop:optimal-value}
it is enough to prove that $s_3(\beta')$ is increasing in $\beta'$, and then (see \eqref{betaprime}) in $\beta$. Since $$s_3(\beta')=\frac{t_3(\beta')}{\beta'}={\mathscr{F}(t_3(\beta'))}$$ we have
\be
\frac{ds_3}{d \beta'}(\beta')= 
\frac{d{\mathscr{F}}}{d t}(t_3(\beta'))
\frac{d{t_3}}{d \beta'}(\beta')
\ee
The derivative of $t_3(\beta')$ can be computed by applying the implicit function theorem:   
\be
\frac{d t_3}{d \beta'}(\beta')= - \frac{t_3(\beta')}{\beta'^2(\frac{d\calF}{dt}(t_3(\beta'))-\frac{1}{\beta'})}.
\ee
Then, since  $\frac{d{\mathscr{F}}}{d t}(t)>0$
and $\frac{d\calF}{dt}(t_3(\beta'))-\frac{1}{\beta'}<0$, we have \be
\frac{ds_3}{d \beta'}(\beta')>0.
\ee
This shows that $s(\beta')$ is strictly increasing.
The existence of a phase transition has been established
in Proposition \ref{prop:optimal-value}, together with the critical inverse temperature $\beta_c'$.
The equal area construction implies that 
$(s(\beta_c),\beta_c)$ are the solutions to the {\em stationarity condition} \eqref{FOPT-1-B1-B>0} and the {\em criticality condition} \eqref{crit-cond-B>0}.
Finally $s(\beta_c)=\lim_{\beta\searrow \beta_c} s(\beta)$ is clearly positive, implying the phase transition is first order.}

The existence of a phase transition has been established
in Proposition \ref{prop:optimal-value}, together with the critical inverse temperature $\beta_c'(B)$.
The equal area construction implies that the pairs
$(\beta_c(B),s_1(\beta_c(B),B))$  and $(\beta_c(B),s_3(\beta_c(B),B))$  are the solutions to the {\em stationarity condition} \eqref{FOPT-1-B1-B>0} and the {\em criticality condition} \eqref{crit-cond-B>0}.
Finally, the difference $s_3(\beta_c(B),B)-s_1(\beta_c(B),B))$ is clearly positive, implying that the phase transition is first order. To show the monotonicity of $\beta \mapsto s^\star(\beta,B)$
%
%
it is enough to prove that $s_1(\beta,B)$ and $s_3(\beta,B)$ are increasing in $\beta'$, and then (see \eqref{betaprime}) in $\beta$. We show this for  $s_3(\beta,B)$, the other case is analogous. Since $$s_3(\beta,B)=\frac{t_3(\beta,B)}{\beta'}={\mathscr{F}_B(t_3(\beta,B))},$$
\be
\frac{\partial s_3}{    \partial \beta'}(\beta,B)= 
\frac{d{\mathscr{F}_B}}{d t}(t_3(\beta,B))
\frac{\partial {t_3}}{\partial \beta'}(\beta,B).
\ee
The derivative of $t_3(\beta,B)$ w.r.t.\ $\beta'$ can be computed by implicit differentiation of the equation ${\mathscr{F}_B(t_3(\beta,B))}-t_3(\beta,B)/\beta'=0$, leading to   
\be
\frac{\partial t_3}{\partial \beta'}(\beta,B)= - \frac{t_3(\beta,B)}{\beta'^2(\frac{d\calF_B}{dt}(t_3(\beta,B))-\frac{1}{\beta'})}.
\ee
Then, since  $t_3(\beta,B)>0$
and $\frac{d\calF_B}{dt}(t_3(\beta,B))-\frac{1}{\beta'}<0$,
\be
\frac{\partial s_3}{\partial \beta'}(\beta,B)>0.
\ee
This shows that $\beta'\mapsto s_3(\beta,B)$ is strictly increasing.
\qed

\invisible{
TEXT WITH B=0

\begin{proposition}[Unique solution for straight lines above $ G_b$]
\label{prop:sol_lim_vec_below_CR}
    For any $0\le \beta' < \frac{1}{\frac{d\mathscr{F}_B}{dt}(t_b)}$, 
     \eqref{eq:stat-cond-t} has only the solution $t_1(\beta,B)\ge 0$. 
\end{proposition}
\begin{proof} Since the function $\calF(\cdot)$ is
strictly concave on $(t_*, +\infty)$, then on this interval its graph lies below the line $G_2(\cdot)$ which is tangent at the point $(t_0, \calF(t_0))$ with $t_0>t_*$. In particular, we have that $G_2(t_*)>\calF(t_*)$. By construction, the line $G_2(\cdot)$ is  above the line passing through the origin and the point $(t_*, \calF(t_*))$. This last line, in turn, is above the graph of $\calF(t)$ on the interval $(0,t_*)$ because of convexity. Hence $G_2(t)$ is above the graph of $\calF(t)$ for all $t>0$, except $t=t_0$ where it is tangent. Any line passing through the origin with slope $1/\beta'> \calF'(t_0)$, being above the line $G_2$, intersects the graph of $\calF(t)$ only at $t_1(\beta')=0$.  
\end{proof}

When $\beta' = \frac{1}{\mathscr{F}'(t_0)}$, 
\eqref{eq:stat-cond-t} becomes 
${\mathscr{F}}(t) = G_2(t)$ and we have 
two solutions: $t=0$ and $t=t_0$. 
If we continue to increase $\beta'$ we enter a region with three solutions, as the following proposition shows.

\begin{proposition}[Three solutions for straight lines between  $G_1$ and  $G_2$]
\label{prop:soln_inside_h}
For any $\beta'\in (\frac{1}{\frac{d\mathscr{F}}{dt}(t_0)}, \frac{q\expec[W]}{\expec[W^2]})$
the equation  \eqref{eq:stat-cond-t}  has three solutions: $0=t_1(\beta')<t_2(\beta')<t_3(\beta') < \infty$.  
\end{proposition}
\begin{proof}
The solution $t_1(\beta')=0$ exists for any $\beta'$.
In order to prove the existence of two additional solutions, we define $\Delta (t)= t/\beta' -\calF(t)$ and
consider separately three sub-cases.

\smallskip
$\bullet$ {\bf Case $\frac{1}{\frac{d\mathscr{F}}{dt}(t_0)}< \beta' < \frac{t_*}{\calF(t_*)}$}: in this case there are no solutions to $\Delta(t)=0$ in the open interval $(0,t_*)$. In fact, the line $t/\beta'$ is above the line $\frac{1}{t_*} \calF(t_*)\,t$ which, by strict convexity, is in turn above $\calF(t)$ on $(0,t_*)$. Now we observe that in the present case $\Delta(t_*)>0$, while $\Delta(t_0)=t/\beta'_0-\calF(t_0)=t/\beta'_0-G_2(t_0)<0$, since $t/\beta'<G_2(t)$ for all $t>0$. Then there exists a solution $t_2(\beta') \in (t_*, t_0)$. Since $\frac{d\calF}{dt}(t)$ is decreasing for $t>t_*$, then we have for all $t\in(t_*,t_0)$ that
$\frac{d \Delta}{dt}(t)
=\frac{1}{\beta'}-\frac{d\calF}{dt}(t)
<\frac{1}{\beta'}-\frac{d\calF}{dt}(t_0)
<0$. This implies that $t_2(\beta')$ is a unique solution in $(t_*,t_0)$. Since $\Delta(t_0)<0$ and $\Delta(t)\to +\infty$ for $t\to +\infty$, we conclude that there exists at least a solution in $(t_0,+\infty)$. This solution is unique because ${\mathscr{F}}(t)$ is strictly concave and bounded for $t\in(t_*,\infty)$.

\smallskip
$\bullet$ {\bf Case $\beta' = \frac{t_*}{\calF(t_*)}$}: obviously  the point $t_2(\beta')=t_*$ is a solution and no other solutions exist in the open interval $(0,t_*)$ since, being $\calF(t)$ strictly convex on the interval, its graph lies below $\frac{\mathscr{F}(t_*)}{t_*}t$. 
Further, we observe that $\frac{d \Delta}{dt}(t_*) = \frac{1}{\beta'} - \frac{d\calF}{dt}(t_*) =
\frac{1}{t_*} \calF(t_*) - \frac{d{\mathscr{F}}}{dt}(t_*) <0$. Then, being $\Delta(t_*)=0$ and 
$\lim_{t\to\infty}\Delta(t)\to +\infty$
we conclude that there exist a  solution $t_3(\beta')$ to $\Delta(t)=0$ in the open interval $(t_*, +\infty)$. This point is unique because $\calF(t)$ is a strictly concave and bounded for $t\in(t_*,\infty)$.

\smallskip
$\bullet$ {\bf Case $\frac{t_*} {\calF(t_*)} < \beta' < \frac{q\expec[W]}{\expec[W^2]}$}: in this case we have that $\Delta(t_*)<0$, while $\Delta(t)>0$ for $t>0$ sufficiently small, since $\Delta(0)=0$ and $\frac{d \Delta}{dt}(0)>0$. Thus there exists a point $0<t_2(\beta')<t_*$ such that 
$\Delta(t_2(\beta'))=0$. This point is unique because
$\calF(t)$ is a strictly convex and increasing function when $t\in(0,t_*)$.
Next we consider the solution on $(t_*,+\infty)$. Since $\Delta(t_*)<0$ and $\Delta(t)\to +\infty$ as $t\to +\infty$ we have that there exist a point $t_3(\beta')\in(t_*, +\infty)$ such that $\Delta(t_3)=0$.
This point is unique because $\calF(t)$ is a strictly concave, increasing and bounded function when $t\in(t_*,\infty)$.
\end{proof}

The solution $s_2(\beta')$ is decreasing in $\beta'$, and vanishes at $\beta'=\frac{q\expec[W]}{\expec[W^2]}$.
If we continue to increase $\beta'$ beyond this point, we  enter a region with two solutions.
\begin{proposition}[Two solutions for straight lines below or equal $G_1$]\label{prop:soln_above_h}
    For any $\beta'\geq \frac{q\expec[W]}{\expec[W^2]}$,
    the equation  \eqref{eq:stat-cond-t}  has two solutions: $0=t_1(\beta')<t_3(\beta') < \infty$.
\end{proposition}
\begin{proof}
    As $t_1(\beta')=0$ is always a solution, we need to prove the existence and uniqueness of the solution $t_3(\beta')$. If $\beta'\geq \frac{q\expec[W]}{\expec[W^2]}$ then the line $t/\beta'$ does not intersect ${\mathscr{F}}$ for $t\in(0,t_*)$, since its slope is less than the slope of $\mathscr{F}$ at $t=0$, and $\mathscr{F}$ is strictly convex in $(0,t_*)$.  For $t > t_*$ the  line $t/\beta'$ is bound to intersect ${\mathscr{F}}$ exactly one time since
    ${\mathscr{F}}(t_*) > \frac{t_*}{\beta'}$, $\mathscr{F}$ is strictly concave and $\mathscr{F}(t)<1$.
\end{proof}

\subsection{Optimal solution}
\label{sec:optimal-solution}

After identifying stationary points, the next step is to determine which is the optimal solution. Calling $s_i(\beta') = \frac{1}{\beta'}  t_i(\beta')$ for $i\in\{1,2,3\}$, where the $t_i(\beta')$ have been described in the previous section, the following lemma establishes the local nature of the stationary points of $p_{\beta',0}$. 

\begin{lemma}[Nature of stationary points] 
\label{lemma:nature-stat-points} 
As we vary $\beta'$ in the interval $[0,\infty)$ we have the following:
\begin{enumerate}
\item[\text{(i)}] 
for any $\beta' < \frac{1}{\frac{d\mathscr{F}}{dt}(t_0)}$ the point  $s_1(\beta')=0$ is a maximum of $p_{\beta',0}$;
\item[\text{(ii)}] 
for $\beta'\in (\frac{1}{\frac{d\mathscr{F}}{dt}(t_0)}, \frac{q\expec[W]}{\expec[W^2]})$,
among the points  $0=s_1(\beta')<s_2(\beta') <s_3(\beta')$, the first and the last are local maxima of $p_{\beta',0}$, while $s_2(\beta')$ is a minimum;
\item[\text{(iii)}] 
for $\beta' >\frac{q\expec[W]}{\expec[W^2]}$
the point $s_1(\beta')=0$ is a minimum of $p_{\beta',0}$, while $s_3(\beta')>0$ is a maximum.
\end{enumerate}
\end{lemma}

\begin{proof}
Starting from \eqref{eq:pressure-p-beta-zero} we compute
the first derivative 
\begin{eqnarray}
\frac{dp_{\beta,0}}{ds}(s) & = &
\beta' \frac{(q-1)}{q} \mathbb{E}[W]
 \left(\expec \left [\frac{W}{\expec[W]}\frac{\e^{\beta'Ws'}-1}{\e^{\beta'Ws'}+q-1}
    \right ]-s'\right) \nonumber\\
    & = &
\beta'\frac{(q-1)}{q} \mathbb{E}[W]\left({\mathscr{F}}(\beta' s) - s\right), 
\end{eqnarray}
and the second derivative
\begin{eqnarray}
\frac{d^2p_{\beta,0}}{ds^2}(s) & = &
\beta'\frac{(q-1)}{q} \mathbb{E}[W]\left(\beta'\frac{d{\mathscr{F}}}{ds}(\beta's) - 1\right). 
\end{eqnarray}
In particular, using \eqref{eq:derivative-F}, we get
\begin{eqnarray}\label{eq:2ndder}
\frac{d^2p_{\beta,0}}{ds^2}(0) & = &
\beta'\frac{(q-1)}{q} \mathbb{E}[W]\left(\beta' \frac{\mathbb{E}[W^2]}{q\mathbb{E}[W]} - 1\right), 
\end{eqnarray}
which implies that $s_1(\beta')=0$ is a (local) maximum for $\beta' <\frac{q\expec[W]}{\expec[W^2]}$, and a minimum otherwise. This proves the nature of the point $s_1(\beta')=0$ in all three items. 

As a consequence, being $s_1(\beta') < s_2(\beta') < s_3(\beta')$ in item (ii), it follows that
$s_2(\beta')$ is a minimum and $s_3(\beta')$ is a maximum. Similarly, being $s_1(\beta') < s_3(\beta')$ in item (iii), it follows that
$s_3(\beta')$ is a maximum.
\end{proof}
Lemma \ref{lemma:nature-stat-points} tells us that
$s=0$ is the best solution for any $\beta' < \frac{1}{\frac{d\mathscr{F}}{dt}(t_0)}$, whereas
$s_3(\beta)$ is best for $\beta'> \frac{q\expec[W]}{\expec[W^2]}$. The optimal value
in the region with three solutions
remains to be established.

\begin{proposition}[Optimal value and ``equal area'' construction]
\label{prop:optimal-value}
There exists $\beta'_c\in (\frac{1}{\frac{d\mathscr{F}}{dt}(t_0)}, \frac{q\expec[W]}{\expec[W^2]})$ such that when $0\le\beta'<\beta'_c$  the optimal value in \eqref{eq:pressure_evaluated_at_soln-a-B=0}
is $s_1(\beta')=0$ and when $\beta'>\beta'_c$ the  optimal value is
$s_3(\beta')$>0.
\end{proposition}
\begin{proof}
Since we already have the results of Lemma \ref{lemma:nature-stat-points}, we only need to 
compare $p_{\beta',0}(0)$ to $p_{\beta',0}(s_3(\beta'))$ when 
$\beta'\in (\frac{1}{\frac{d\mathscr{F}}{dt}(t_0)}, \frac{q\expec[W]}{\expec[W^2]})$.
To this aim, we define 
\begin{align}
R(\beta')=\int_0^{t_3(\beta')} (\mathscr{F}(t)-\frac{t}{\beta'})dt \label{eq:def_R}   
\end{align}
and observe that
\begin{eqnarray}
p_{\beta',0}(s_3(\beta')) - p_{\beta',0}(0) & = &
\beta' \frac{(q-1)}{q} \mathbb{E}[W]
\int_{0}^{s_3(\beta')}  \left(\expec[\frac{W}{\expec[W]}\frac{\e^{\beta'Ws'}-1}{\e^{\beta'Ws'}+q-1}
    ]-s'\right) ds'\nonumber\\
    & = &
\frac{(q-1)}{q} \mathbb{E}[W] 
R(\beta')
\end{eqnarray}
We now prove that the mapping
$\beta'\mapsto R(\beta')$ is monotonically increasing on the compact set 
$\beta'\in [\frac{1}{\frac{d\mathscr{F}}{dt}(t_0)}, \frac{q\expec[W]}{\expec[W^2]}]$.
Furthermore we prove that 
$R({\frac{1}{\frac{d\mathscr{F}}{dt}(t_0)}})<0$ 
and $R(\frac{q\expec[W]}{\expec[W^2]}) >0$.
Therefore, by the intermediate value theorem, there exists a unique $\beta_c' \in (\frac{1}{\frac{d\mathscr{F}}{dt}(t_0)}, \frac{q\expec[W]}{\expec[W^2]})$ such that $R(\beta_c')=0$
and the claim of the proposition follows.

To prove monotonicity we compute
\begin{eqnarray}
\label{eq:derivative-R}
\frac{d R}{d \beta'}(\beta')
&=&
\left (\mathscr{F}(t_3(\beta'))-\frac{t_3(\beta')}{\beta'}\right ) \frac{d t_3}{d\beta'}(\beta')
+
\int_0^{t_3(\beta')} \frac{t}{{\beta'}^2}\, dt  
\nonumber\\
&=&
\frac{t_3(\beta')^2}{2 {\beta'}^2}>0
\end{eqnarray}
where we observe that the term inside the round bracket is zero because $t_3(\beta')$ is a stationary point
and the differentiability of $t_3(\beta')$ is guaranteed 
by the implicit function theorem.
Thus \eqref{eq:derivative-R} shows that $R(\beta')$ is strictly increasing for $\beta'\in [\frac{1}{\frac{d\mathscr{F}}{dt}(t_0)}, \frac{q\expec[W]}{\expec[W^2]}]$.

To see that $R({\frac{1}{\frac{d\mathscr{F}}{dt}(t_0)}})<0$, use that $t_3({\frac{1}{\frac{d\mathscr{F}}{dt}(t_0)}})=t_0$ and recall that the tangent line $G_2(t)=t\frac{d\mathscr{F}}{dt}(t_0)$ completely lies above the graph of $\mathscr{F}$, with the only points of intersections being $(0,0)$ and $(t_0,\mathscr{F}(t_0))$. Therefore $\mathscr{F}(t)-t\frac{d\mathscr{F}}{dt}(t_0) <0$ for all $t\in(0,t_0)$.

To see that $R(\frac{q\expec[W]}{\expec[W^2]}) >0$
the argument is similar. First we note that the tangent line $G_1(t)=\frac{\expec[W^2]}{q\expec[W]} t$,
is below the graph of $\mathscr{F}$ in the interval $[0,t_*]$. This is because in this interval $\mathscr{F}$ is convex, and $G_1(t)$ is a tangent line to the graph of $\mathscr{F}$ at a point inside this interval. 
Next, consider the line $L(t)$ joining the points $(t_*,\mathscr{F}(t_*))$ and $(t_3(\beta_1'),\mathscr{F}(t_3(\beta_1')))$. Due to concavity of $\mathscr{F}$ in the interval $[t_*,t_3(\beta_1')]$, the line $L(t)$ is strictly below the graph of $\mathscr{F}$. Since $G_1(t_*)<L(t_*)=\mathscr{F}(t_*)$, this means the line $G_1(t)$ is below the graph $\mathscr{F}$ in the interval $[t_*,t_3(\frac{q\expec[W]}{\expec[W^2]})]$. Combining everything, we get that $G_1(t)$ is less or equal $\mathscr{F}(t)$ for all $t\in[0,t_3(\frac{q\expec[W]}{\expec[W^2]})]$. Consequently, $R(t_3(\frac{q\expec[W]}{\expec[W^2]})>0$.
\end{proof}

\invisible{
\Gib{The following is the former proof of Lemma 5.10 (restricted to continuity).}
\begin{proof}
Note that for fixed $t$, $\mathscr{F}(t)-t/\beta'$ is a monotone increasing continuous function of $\beta'$ on $[\beta_0',\beta_1']$. Hence if we can show that $t_3(\beta')$ is also continuous and monotone increasing as a function of $\beta'$, with the definition $t_3(\beta_0')=t_0$, we can conclude $R(\beta')$ is also continuous in the same domain.  

This is because, for $\beta_1 < \beta_2$, where $\beta_1,\beta_2 \in [\beta_0',\beta_1']$, one can write 
\begin{align*}
    |R(\beta_1)-R(\beta_2)|=&\left|\int_{0}^{t_3(\beta_1)}(\mathscr{F}(t)-t/\beta_1)dt-\int_{0}^{t_3(\beta_2)}(\mathscr{F}(t)-t/\beta_2)dt \right|\\&\leq \int_{0}^{t_3(\beta'_1)}\mathbf{1}_{[t_3(\beta_1),t_3(\beta_2)]}\left(|\mathscr{F}(t)-t\mathscr{F}'(t_0)|\vee |\mathscr{F}(t)-\frac{t\expec[W^2]}{q\expec[W]}| \right)dt.
\end{align*}
In particular, bounding from above the continuous function $\left(|\mathscr{F}(t)-t\mathscr{F}'(t_0)|\vee |\mathscr{F}(t)-\frac{t\expec[W^2]}{q\expec[W]}| \right)$ of $t$ in the compact interval $[0,t_3(\frac{q\expec[W]}{\expec[W^2]})]$ by its maximum value, we note that when $\beta_2-\beta_1 \to 0$, the r.h.s.\ above, and hence $|R(\beta_1)-R(\beta_2)|$ also, tend to zero.
{The proof that $t_3(\beta')$ is continuous and increasing is a tedious analysis argument, and can be found in Section \ref{sec-t3-mon}.}


Note that from the definition of $R$, it is clear that it is increasing in $\beta'$ if $t_3$ is so. Further, we claim that $R(\frac{1}{\mathscr{F}'(t_0)})<0$. To see this, recall from the proof of Proposition \ref{prop:sol_lim_vec_below_CR} that the tangent line $G_2(t)=t\mathscr{F}'(t_0)$ completely lies above the graph of $\mathscr{F}(t)$ in the interval $[0,t_0]$, with the only points of intersections being $(0,0)$ and $(t_0,\mathscr{F}(t_0))$. Next, we claim that $R(\frac{q\expec[W]}{\expec[W^2]})>0$. The argument is similar. First we note that the tangent line $G_1(t)=\frac{t\expec[W^2]}{q\expec[W]}$ is below the graph of $\mathscr{F}$ in the interval $[0,t_*]$. This is because in this interval $\mathscr{F}$ is convex, and $G_1(t)$ is a tangent line to the graph of $\mathscr{F}$ at a point inside this interval. 

Consider the line $L_2(t)$ joining $(t_*,\mathscr{F}(t_*))$ and $(t_3(\frac{q\expec[W]}{\expec[W^2]}),\mathscr{F}(t_3(\frac{q\expec[W]}{\expec[W^2]})))$. Due to concavity of $\mathscr{F}$ in the interval $[t_*,t_3(\frac{q\expec[W]}{\expec[W^2]})]$, the line $L_2(t)$ is strictly below the graph of $\mathscr{F}$. Since $G_1(t_*)<L_2(t_*)=\mathscr{F}(t_*)$, this means the line $G_1(t)$, which is also the line joining $(t_*,G_1(t_*))$ and $(t_3(\frac{q\expec[W]}{\expec[W^2]}),\mathscr{F}(t_3(\frac{q\expec[W]}{\expec[W^2]})))$ is below the graph $\mathscr{F}$ in the interval $[t_*,t_3(\frac{q\expec[W]}{\expec[W^2]}))]$. Combining everything, we get that the line $G_2(t)$ is always below the graph of $\mathscr{F}$ in the interval $[0,t_3(\frac{q\expec[W]}{\expec[W^2]}))]$. Consequently, $R(t_3(\frac{q\expec[W]}{\expec[W^2]})))>0$. Since $R$ is also monotone increasing, as proved earlier, this means by the intermediate value theorem, $R$ has a unique root $\beta_c'\in (\beta_0',\beta_1')$. 
In particular, writing out $R(\beta_c')=0$, we note that $\beta'_c$ solves
\begin{align}
    \beta'_c=\frac{(t_3(\beta_c))^2}{2\int_0^{t_3(\beta'_c)}\mathscr{F}(t)dt}.
\end{align}

\end{proof}

}

\subsection{Conclusion of the proof of Theorems \ref{thm-nature-pt-a} and \ref{thm-nature-pt-b}}
\label{sec:conclusion-proof}

First we prove that $\beta_c$ is well defined, i.e.
if $s(\beta_1)>0$ and $\beta_2>\beta_1$, then also $s(\beta_2)>0$. From Proposition \ref{prop:optimal-value}
it is enough to prove that $s_3(\beta')$ is increasing in $\beta'$, and then (see \eqref{betaprime}) in $\beta$. Since $$s_3(\beta')=\frac{t_3(\beta')}{\beta'}={\mathscr{F}(t_3(\beta'))}$$ we have
\be
\frac{ds_3}{d \beta'}(\beta')= 
\frac{d{\mathscr{F}}}{d t}(t_3(\beta'))
\frac{d{t_3}}{d \beta'}(\beta')
\ee
The derivative of $t_3(\beta')$ can be computed by applying the implicit function theorem:   
\be
\frac{d t_3}{d \beta'}(\beta')= - \frac{t_3(\beta')}{\beta'^2(\frac{d\calF}{dt}(t_3(\beta'))-\frac{1}{\beta'})}.
\ee
Then, since  $\frac{d{\mathscr{F}}}{d t}(t)>0$
and $\frac{d\calF}{dt}(t_3(\beta'))-\frac{1}{\beta'}<0$, we have \be
\frac{ds_3}{d \beta'}(\beta')>0.
\ee
This shows that $s(\beta')$ is strictly increasing.
The existence of a phase transition has been established
in Proposition \ref{prop:optimal-value}, together with the critical inverse temperature $\beta_c'$.
The equal area construction implies that 
$(s(\beta_c),\beta_c)$ are the solutions to the {\em stationarity condition} \eqref{FOPT-1-B1-B>0} and the {\em criticality condition} \eqref{crit-cond-B>0}.
Finally $s(\beta_c)=\lim_{\beta\searrow \beta_c} s(\beta)$ is clearly positive, implying the phase transition is first order.

\begin{figure}
\centering
\includegraphics[width=10cm]{fig1v2.jpg}
\caption{The function $t\mapsto \mathscr{F}(t)$ for
  the Erd\H{o}s–R\`enyi model with {\color{red}$w=5$ and $q=21$}.
  The two tangent lines $G_1(t)$ and $G_2(t)$ are shown.
  The slope of the line associated to the ``equal-area construction'' identifies the critical  temperature $\frac{1}{\beta_c}$.}
  \label{fig:second-deriv}
\end{figure}

\begin{figure}
\centering
\includegraphics[width=10cm]{fig1v3.jpg}
\caption{The function $t\mapsto \mathscr{F}(t)$ for
  the Erd\H{o}s–R\`enyi model with {\color{red}$w=5$ and $q=21$}.
  The two tangent lines $G_1(t)$ and $G_2(t)$ are shown.
  The slope of the line associated to the ``equal-area construction'' identifies the critical  temperature $\frac{1}{\beta_c}$.}
  \label{fig:second-deriv}
\end{figure}
\vspace{1.cm}


}

\invisible{

OLD TEXT WRITTEN BY CLAUDIO

\subsection{Analysis of the fixed-point equation for $s(\beta)$}

Recall \eqref{eq:for_optimal_s-B1-rep}. We change variables $t=s\beta'$, and look for solutions in $[0,\beta']$ of
\be \label{eq:for_optimal_s_2}
\expec \left[\frac{W}{\expec[W]}\frac{\e^{tW}-1}{\e^{tW}+q-1} \right]=\frac{t}{\beta'}.
\ee
Note that there are no solutions in $t$ of \eqref{eq:for_optimal_s_2} in $(\beta',\infty)$ since the l.h.s.\ of \eqref{eq:for_optimal_s_2} is bounded by $1$.
Recall $\mathscr{F}(t)$ in \eqref{def-Fcal}.
Note that $\mathscr{F}(0)=0$ and $\lim_{t \to \infty}\mathscr{F}(t)=1$. Further,
\be
\frac{d}{dt}\mathscr{F}(t)=\expec\left[\frac{W^2}{\expec[W]}\frac{q\e^{tW}}{(\e^{tW}+q-1)^2}\right]>0
\ee
which means $\mathscr{F}$ is increasing.
Further,
\be\label{eq:second_derivative_gen-rep}
\frac{d^2}{dt^2}\mathscr{F}(t)=q\expec\left[\frac{W^3}{\expec[W]}\frac{\e^{tW}[(q-1)^2-\e^{2tW}]}{(\e^{tW}+q-1)^4}\right]=\frac{q}{\expec[W]}\expec\left[W^3\frac{(q-1)\e^{tW}-\e^{2tW}}{(\e^{tW}+q-1)^3}\right].
\ee
Observe that for $q\leq 2$, $\mathscr{F}$ is always concave since $\frac{d^2}{dt^2}\mathscr{F}(t)<0$, while for $q>2$, $\mathscr{F}$ is convex near $0$, while concave near $\infty$. This behavior change explains the change from a second-order phase transition in Ising ($q=2$) and site percolation ($q=1$) to a first-order transition in Potts with $q\ge 3$.

Further observe that because of this form of the function $\mathscr{F}$, when we look for roots of $\mathscr{F}(t)-\frac{t}{\beta'}$, for small values of $\beta'$, there is only one root, namely $t=0$, since the line $t/\beta'$ only intersects the curve $\mathscr{F}(t)$ once at $t=0$. While for large $\beta'$, there are potentially more points of intersection of $\mathscr{F}(t)$ with $t/\beta'$, and the belief is that the largest root is relevant - this root will give us the optimal $t$, and hence the optimal $s$. 

\invisible{\RvdH{TO DO 8: We should work with $B>0$, and show that $s(\beta,B)$ is increasing, so that $s(\beta)=\lim_{B\searrow 0} s(\beta, B)$ is well defined. Claudio will do this.}}

In order to study the solution $s(\beta,B)$ to equation \eqref{FOPT-1-B1-B>0} as a function of $B$, let us denote
\be\label{eq:fun_optimal_B}
{\tt F}(s,\bpr,B)= \mathbb{E}\left[ \frac{W}{\expec[W]} \frac{\e^{\beta' W s+B}-1}{\e^{\beta' Ws+B}+q-1}\right] - s,
\ee
and write \eqref{FOPT-1-B1-B>0}  as
\be\label{FOPT-1-B1-bis}
{\tt F}(s(\beta,B),\bpr,B)=0.
\ee
\begin{proposition}
For all
 $\bpr < \frac{\expec[W]}{\expec[W^2]}$ the solution $s(\beta,B)$ to equation \eqref{FOPT-1-B1-B>0} is increasing in $B$, so that $s(\beta')=\lim_{B\searrow 0} s(\beta,B)$ is well defined. 
\end{proposition}
\begin{proof}
  (I assume here that the existence of the solution to \eqref{FOPT-1-B1-bis} has been already proven elsewhere.) Let us compute the derivative w.r.t. $B$ of \eqref{FOPT-1-B1-bis}, we get
  \be\label{eq:der-eq-B}
  \frac{\partial {\tt F} }{\partial s}(s(\beta,B),\bpr,B)\cdot \frac{\partial s}{\partial B}(\beta,B) + \frac{\partial {\tt F} }{\partial B}(s(\beta,B),\bpr,B)=0,
  \ee
  where the partial derivatives of $\tt F$, that clearly exist, are given by:
  \be\label{eq:der-f-s}
   \frac{\partial {\tt F} }{\partial s}(s,\bpr,B)= 
    \mathbb{E}\left[ \frac{W^2}{\expec[W]}  \frac{\bpr q\ \e^{\beta' W s+B}}{(\e^{\beta' Ws+B}+q-1)^2}\right] - 1,
  \ee
  and 
    \be
   \frac{\partial {\tt F} }{\partial B}(s,\bpr,B)= 
    \mathbb{E}\left[ \frac{W}{\expec[W]}  \frac{ q\ \e^{\beta' W s+B}}{(\e^{\beta' Ws+B}+q-1)^2}\right].
  \ee
  Since $\frac{\partial {\tt F} }{\partial B}(s,\bpr,B)$ is clearly positive, we now want to prove that $ \frac{\partial {\tt F} }{\partial s}(s,\bpr,B)$ is negative provided that $\bpr$ is sufficienlty small. These facts along with \eqref{eq:der-eq-B}, show that $s(\bpr, B)$ is increasing.
  Since $q\,x/(x+q-1)<1$ for all $q>2$ and $x\ge 0$,
  we bound the average in \eqref{eq:der-f-s} as follows:
 \eqan{ 
   \mathbb{E}\left[ \frac{W^2}{\expec[W]}  \frac{\bpr q\ \e^{\beta' W s+B}}{(\e^{\beta' Ws+B}+q-1)^2}\right] <  
   \bpr \frac{\expec[W^2]}{\expec[W]}.
  }
  Therefore, for $\bpr < \frac{\expec[W]}{\expec[W^2]}$ 
  we have
  $\frac{\partial {\tt F} }{\partial s}(s,\bpr,B)<0$,
  as claimed. This proves the lemma.
\invisible{\eqan{
  &\frac{W^2}{\expec[W]}  \frac{\bpr q\ \e^{\beta' W s+B}}{(\e^{\beta' Ws+B}+q-1)^2}< \frac{W^2}{\expec[W]}  \frac{\bpr q\ \e^{\beta' W s+B}}{\e^{\beta' Ws+B} (\e^{\beta' Ws+B}+q-1)}<
   \frac{ W^2 \bpr q}{\expec[W] (\beta' Ws+q)}=\nonumber \\
   & \frac{W^2 \bpr q}{\expec[W]  W \bpr q (\frac{s}{q}+1)}
   }
   }
\end{proof}
From now on, we work with $t(\beta)=s(\beta)/\beta'$, which satisfies
\be\label{eq:for_optimal_t-0}
    \mathbb{E}\left[ \frac{W}{\expec[W]} \frac{\e^{tW}-1}{\e^{t W }+q-1}\right] = \frac{t}{\beta'}.
    \ee
\subsection{Identifying the critical value}
Recall the function $\mathscr{F}(t)$ from \eqref{def-Fcal}. Throughout this section, we assume that $t\mapsto \frac{d^2}{dt^2}\mathscr{F}(t)$ is first positive and then negative. This has the following consequence:
\begin{lemma}[Inflection point]
\label{lem-inflection}
Assume that $t\mapsto \frac{d^2}{dt^2}\mathscr{F}(t)$ is first positive and then negative. Then, there is a point $t_* \in (0,\infty)$, which is the inflection point of $\mathscr{F}$: $\frac{d^2}{dt^2}\mathscr{F}(t_*)=0$, and $\frac{d^2}{dt^2}\mathscr{F}(t)>0$ in $(0,t_*)$, while $\frac{d^2}{dt^2}\mathscr{F}(t)<0$ in $(t_*,\infty)$. 
\end{lemma}

Throughout this section, we further assume that $\expec[W^2]<\infty$. 
\invisible{Since
    \eqn{
    \frac{d}{dt}\mathscr{F}(0)=\frac{\expec[W^2]}{q\expec[W]}\equiv 1/\beta_0',
    }
there are exactly {\em two} solutions to $\mathscr{F}(t)=t/\beta'$ whenever $\beta'<\beta_0'$.}

We next identify the value $\beta_1'$ such that for $\beta'>\beta'_1$, the unique solution to \eqref{eq:for_optimal_t-0} is $t=0:$

\begin{proposition}[First non-zero intersection point]\label{prop:tangents}
Assume that $\expec[W^2]<\infty$, and that $t\mapsto \frac{d^2}{dt^2}\mathscr{F}(t)$ is first positive and then negative. Then, the equation $\mathscr{F}(t)=t(\frac{d}{dt}\mathscr{F}(t))$ has a unique positive root $t_0>0$. Further, there are only two tangent lines to the curve $\mathscr{F}(t)$ which also pass through $(0,0)$. One of them is $G_1(t)=\left (\frac{q\expec[W]}{\expec[W^2]}\right )^{-1}t$, which is the tangent to $\mathscr{F}(t)$ at $t=0$. The other one is $G_2(t)=(\frac{d}{dt}\mathscr{F}(t_0))t$, which is the tangent to $\mathscr{F}(t)$ 
at a point $(t_0,\mathscr{F}(t_0))$, with $t_0>t_*$. Further, $\beta_c'>1/t_0,$ so that, in particular, $\beta_c>0$.
\end{proposition}

\begin{proof} Let us consider the intercept $q(t)=\calF(t)-t\calF^\prime(t)$ of the tangent to $\calF(t)$ for $t \in [0, t_*]$ and $t \in (t_*, \infty)$. We want to check that 
besides $t=0$ there exists a unique $t_0\ne 0$ such that $q(t_0)=0$. Since $q^\prime(t)=-t \calF^{\prime \prime}(t)$, and thanks to Lemma \ref{lem-inflection}, we have that $q(t)$ is decreasing in $(0,t_*)$ and increasing in $(t_*, \infty)$.
Thus, from the fact that $q(0)=0$ (this corresponds to the tangent $G_1(t)$) we have that  $q(t_*)\le q(t)<0$ for $0<t \le t_*$, which means that in this interval there are no tangent lines passing through $(0,0)$. On $t>t_*$ the intercept $q(t)$ increases monotonically from $q(t_*)<0$ to $1$ (since $q(t)\to 1$ as $t\to +\infty$). Then by continuity there is a unique $t_0 > t_*$ such that $q(t_0)=0$, which corresponds to a tangent line passing through the origin.
\end{proof}

We define the region $\frac{1}{\mathscr{F}'(t_0)}<\beta'<\frac{q\expec[W]}{\expec[W^2]}$ as the \textit{triple-intersection region}.

Recall the pressure $p(s,\beta')$ from \eqref{eq:pressure_evaluated_at_soln} which we now consider  as a function of $t$, i.e. we introduce $r(t):=p(\frac{t}{\beta'}, \beta')$ and solve the optimization problem w.r.t. the variable $t$. The first and second derivatives of $r(t)$, can be expressed in terms of  $\calF(t)$:
\begin{align}
& r'(t)=\frac{q-1}{q}\, \expec[W] \left ( \calF(t)- \frac{t}{\beta'} \right)\label{eq:derfr},\\
& r''(t)=\frac{q-1}{q} \expec[W]\, 
\left ( \calF'(t)-\frac{1}{\beta'} \right).\label{eq:dersr}
\end{align}
Therefore, to optimize $t\mapsto p(\frac{t}{\beta'}, \beta'),$ we need to look for solutions of $\calF(t)- \frac{t}{\beta'}=0$. In the following three propositions, we identify the number of such solutions in the regions $\beta'<\frac{1}{\mathscr{F}'(t_0)},
\beta'\in (\frac{1}{\mathscr{F}'(t_0)}, q\expec[W]/\expec[W^2])$ and $\beta'>q\expec[W]/\expec[W^2]$, respectively, which we call the {\em $\beta'$ below}, {\em in} and {\em above the triple-intersection regions}, respectively. 

\begin{proposition}\label{prop:zeromax}
    For any $\bpr$ below and in the triple-intersection region (i.e. $0<\bpr<q\expec[W]/\expec[W^2])$ the point $t=0$ is a (possibly local) maximum of the pressure. For $\bpr$ above the triple-intersection regions (i.e. $\beta'>q\expec[W]/\expec[W^2]$) $t=0$ is a miniumum.
\end{proposition}
\begin{proof}
The fact that $t=0$ satisfies the stationarity condition \eqref{eq:for_optimal_s_2} (i.e. $r'(t)=0$) is obvious. The thesis is easily obtained by observing, see \eqref{eq:dersr}, that
    $$
    r''(0)=\frac{q-1}{q} \expec[W]
\left ( \frac{\expec[W^2]}{q\expec[W]}-\frac{1}{\beta'} \right).
    $$   
\end{proof}
In the following we will see that, depending on $\bpr$, $t=0$ (and hence $s=0$) may be the global or a local maximizer of the pressure.

We start showing that there is a unique solution for $\beta'$ above the triple-intersection region. Recall that $t_0$ is the unique positive solution to $\mathscr{F}(t_0)=\mathscr{F}'(t_0)t_0$.
\begin{proposition}[Unique solution for $\beta'$ below triple-intersection region]\label{prop:sol_lim_vec_below_CR}
    For any $\beta'\leq \frac{1}{\mathscr{F}'(t_0)}$, let $t_1(\beta')$ be the largest solution to $\mathscr{F}(t)=\frac{t}{\beta'}$. Then $t_1(\beta')=0$ when $\beta'<\frac{1}{\mathscr{F}'(t_0)}$  and $t_1(\frac{1}{\mathscr{F}'(t_0)})=t_0>0$. Further, for any $\beta'\leq \frac{1}{\mathscr{F}'(t_0)}$ the optimizer of \eqref{eq:laplace} has the form \eqref{eq:soln_form}, with $s(\beta')=0$.
\end{proposition}
\begin{proof} Recall $t_*$ is the inflection point of $\calF(t)$. Because of its strict concavity on $(t_*, +\infty)\ni t_0$, the graph of $\calF(t)$ lies below the tangent line $G_2$ at the point $(t_0, \calF(t_0))$. In particular, we have that $G_2(t_*)>\calF(t_*)$. The line $L$ passing through the origin and the point $(t_*, \calF(t_*))$ is above the graph of $\calF(t)$ on the interval $(0,t_*)$ because of convexity. Since $G_2$ is above $L$ (due to $G_2(t_*)>\calF(t_*)$) we conclude that $G_2$ is above the graph of $\calF(t)$ for all $t>0$. Hence, any line passing through the origin with slope $1/\beta'> \calF'(t_0)$, being above $G_2$, intersects the graph of $\calF(t)$ only at $t_1(\beta')=0$. Then, when $1/\beta'> \calF'(t_0)$, 
by Proposition \eqref{prop:zeromax} the point $t_1(\beta')=0$ is the (unique) maximum point. 
When the slope is $1/\beta'=\calF'(t_0)$ there are two intersection points: at $t=0$ and at $t_1(\bpr)={\mathscr{F}'(t_0)}=t_0>0$. Since also in this case $t=0$ is a maximum point and no other stationarity points besides $t_1(\bpr)$ exist, we conclude that  $t_1(\bpr)$ cannot be a maximum; hence $t=0$ is the global optimizer also for $1/\beta'=\calF'(t_0)$.
By Theorem \ref{thm-soln_form},  the optimizer
of \eqref{eq:laplace} is of the form \eqref{eq:soln_form} parametrized by $s(\beta')$ that for all $\beta'\leq \frac{1}{\mathscr{F}'(t_0)}$ is equal to $0$ (recall that $t=s \bpr$). 
\end{proof}



\begin{proposition}[Three solutions for $\beta'$ inside the triple-intersection region]\label{prop:soln_inside_h}
For any $\beta'$ in the triple-intersection region, the equation $\mathscr{F}(t)=\frac{t}{\beta'}$ has three solutions: $0=t_1(\beta')<t_2(\beta')<t_3(\beta')$. Further, both $t_1(\beta')$ and $t_3(\beta')$ are local maximum points of $p(\frac{t}{\beta'},\beta')$, while $t_2(\beta')$ is a local minimum point of $p(\frac{t}{\beta'},\beta')$. 
\end{proposition}
We emphasize that we do not yet identify what is the {\em optimal} solution inside the triple-intersection region. This will be done in the next section.

\begin{proof}
The solution $t_1(\beta')=0$ exists for any $\beta'$ and is a local maximum point of $r(t)$. Indeed $r'(0)$ is trivially $0$ and $r''(0)=\frac{q-1}{q} \expec[W] (\calF'(0)-\frac{1}{\beta'})=\frac{q-1}{q} \expec[W] 
(\frac{\expec[W^2]}{q\expec[W]} -\frac{1}{\beta'})<0$ since in the triple-intersection region $\beta'<\frac{q\expec[W]}{\expec[W^2]}$.\\ 
Now, recalling that $t_0$ is the value that defines the "upper" tangent line  $G_2(t)$  (see Proposition \ref{prop:tangents}), we introduce the  line $L(t)=t/\beta'$ with $\frac{\expec[W^2]}{q\expec[W]} <\frac{1}{\beta'}<\calF'(t_0)$. In order to prove the existence of the solutions $t_2(\beta')$ and $t_3(\beta')$ to the equation $\calF(t)=L(t)$, we consider separately the three cases: $\frac{1}{\beta'} < \frac{1}{t_*} \calF(t_*)$, $\frac{1}{\beta'} = \frac{1}{t_*} \calF(t_*)$ and $\frac{1}{\beta'} > \frac{1}{t_*} \calF(t_*)$, where $t_*$ is the unique inflection point whose existence is proven in Lemma \ref{lem-inflection}.
In the following we denote $\Delta (t)=L(t)-\calF(t)$. 

\smallskip
$\bullet$ {\bf Case ${\frac{1}{\beta'} < \frac{1}{t_*} \calF(t_*)}$}. In this case we have that $\Delta(t_*)<0$, while $\Delta(t)>0$ for $t>0$ sufficiently small, since $\Delta(0)=0$ and $\Delta'(0)>0$. Thus there exists a point $\hat{t}>0$ smaller than $t_*$ such that 
$\calF(\hat{t})=L(\hat{t})$. We want to prove that this point is unique in $(0,t_*)$.
Arguing by contradiction, we assume that there are several solutions  $0<\hat{t}_1<\hat{t}_2<\ldots$ to the equation $\Delta(t)=0$. Then by Lagrange's mean value theorem there should exist at least two  points $\xi_1 \in (0, \hat{t}_1)$ and $\xi_2 \in (\hat{t}_1, \hat{t}_2)$ such that $\calF'(\xi_i)=L'(\xi_i)\equiv\frac{1}{\beta'}$, but this is not possible because $\calF'(t)$ is strictly increasing on $(0,t_*)$. This proves that $\hat{t}$ is the unique solution  in the open interval $(0,t_*)$; we denote it $t_2(\beta')$.\\
For further reference we observe also that $\calF'(t_*)>\frac{1}{\beta'}$. Indeed, again by the mean value theorem there exists a $\xi \in (0, t_2(\beta'))$ such that $\calF'(\xi)=\bpinv$ and since 
$\calF'(t)$ is increasing on $(0,t_*)$ it follows that $\bpinv=\calF'(\xi)<\calF'(t_2(\beta'))<\calF'(t_*)$. The first inequality shows that $r''(t_2(\beta'))=\frac{q-1}{q} \expec[W] (\calF'(t_2(\beta'))-\frac{1}{\beta'})>0$, proving that $t_2(\beta')$ is a local minimum of the pressure $p(\frac{t}{\beta'})$.\\
Now we consider the solution on $(t_*,+\infty)$. Since $\Delta(t_*)<0$ and $\Delta(t)\to +\infty$ as $t\to +\infty$ (because $\calF(t)$ is bounded) we have that there exist solutions to $\Delta(t)=0$ in $(t_*, +\infty)$; let $\hat{t}$ the smallest one. In order to prove that this solution is unique in $(t_*, +\infty)$ we consider the function $\Delta_1(t)=L_1(t)-\calF(t)$, where $L_1(t)= \frac{1}{\beta'}(t-t_*)+\calF(t_*)$. In a right neighborhood of $t_*$ we have that $\Delta_1(t)<0$, since  $\Delta_1(t_*)=0$ and $\Delta'_1(t_*)<0$ (recall that $\calF'(t_*)>\frac{1}{\beta'}$). From the fact that $\Delta_1(t)\to +\infty$ as $t\to +\infty$ we conclude that there exist solutions to $\Delta_1(t)=0$ on $(t_*,+\infty)$.  Let us denote by $\tilde{t} > t_*$ the smallest of such solutions. We observe that $\tilde{t}$ is smaller than $\hat{t}$, because $\Delta(t_*)<\Delta_1(t_*)=0$ and $\Delta'_1(t)=\Delta'(t)$. 
On the other hand, by mean value theorem there exists a $\xi \in (t_*,\tilde t)$ such that $\calF'(\xi)=L_1'(t)\equiv \frac{1}{\beta'}$. Then for $t>\xi$, we have $\calF'(t)<\calF'(\xi)=\frac{1}{\beta'}$ because $\calF'(t)$ is strictly decreasing on $(t_*, +\infty)$. Then, since $\hat{t}>\tilde t> \xi$ we have also $\calF'(\hat{t})<\calF'(\xi)=\frac{1}{\beta'}$. \\
Now suppose that $\hat{t}_1>\hat{t}$ is a further solution to $\Delta(t)=0$. Then by the standard argument that makes use of the mean value theorem, we have a point $\xi_1\in (\hat{t}, \hat{t}_1)$ such that $\calF'(\xi_1)=\frac{1}{\beta'}$. But this
is not possible because  $\calF'(t)<\calF'(\xi)=\frac{1}{\beta'}$ for $t>\xi$, as previously observed. This contradiction proves that $\hat{t}$ is the unique solution in $(t_*, +\infty)$. We call it $t_3(\beta')$. This is a local maximum point of the pressure, in fact $r''(t_3(\beta'))=\frac{q-1}{q} \expec[W] (\calF'(\hat{t})-\frac{1}{\beta'})<0$, thanks to the inequality proven above. \\
$\bullet$ {\bf Case ${\frac{1}{\beta'} = \frac{1}{t_*} \calF(t_*)}$}. Obviously  the point $t=t_*$ is a solution and no other solutions exist in the open interval $(0,t_*)$ since, being $\calF(t)$ strictly convex on the interval, its graph lies below $L(t)$. By repeating the argument used in the previous case, we can show that $\calF'(t_*)>\bpinv\equiv \frac{1}{t_*} \calF(t_*)$. Then, being  $\Delta'(t_*)<0$ and $\Delta(t_*)=0$, with 
$\Delta(t)\to +\infty$ as $t\to +\infty$
we conclude that there exist solutions to $\Delta(t)=0$ in the open interval $(t_*, +\infty)$. Let $\hat{t}$ be the smallest of such solutions. The mean value theorem implies that there exists $\xi \in (t_*,\hat{t})$ such that $\calF'(\xi)=L'(\xi)\equiv \frac{1}{t_*}\calF(t_*)$ and, since $\calF'(t)$ is decreasing in $(t_*, +\infty)$, we conclude that $\calF'(\hat{t}) < \frac{1}{\beta'} \equiv \frac{1}{t_*}\calF(t_*)$. A further solution, say $\hat{t}_1>\hat{t}$, would imply the existence of a $\xi_1 \in (\hat{t}, \hat{t}_1)$ such that $\calF'(\xi_1)=\frac{1}{t_*}\calF(t_*)$ then, being $\calF'(t)$ decreasing we would have $\calF'(\hat{t})>\frac{1}{t_*}\calF(t_*)$ which leads to a contradiction that proves the uniqueness of the solution in $(t_*,+\infty)$. In this case we set $t_2(\beta'):=t_*$ and $t_3(\beta'):=\hat{t}$. The inequalities proven above give $r''(t_2(\beta'))=\frac{q-1}{q} \expec[W] (\calF'(t_*)-\frac{1}{\beta'})>0$ and $r''(t_3(\beta'))=\frac{q-1}{q} \expec[W] (\calF'(\hat{t})-\frac{1}{\beta'})<0$, showing that the two points are respectively a local minimum and a local maximum of the pressure.\\
$\bullet$ {\bf Case ${\frac{1}{\beta'} > \frac{1}{t_*} \calF(t_*)}$} (with the further condition that $\bpinv < \calF'(t_0)$). In this case there are no solutions to $\Delta(t)=0$ in the open interval $(0,t_*)$. In fact, the line $L(t)$ is above the line $G_*(t)=\frac{1}{t_*} \calF(t_*)\,t$ which, by strict convexity is, in its turn, above $\calF(t)$ on $(0,t_*)$. Now we observe that in the present case $\Delta(t_*)>0$, while $\Delta(t_0)=L(t_0)-\calF(t_0)=L(t_0)-G_2(t_0)<0$, since $L(t)<G_2(t)$ for all $t>0$. Then there exists at least a solution $\hat{t}\in (t_*, t_0)$. Since $\calF'(t)$ is decreasing for $t>t_*$, we have on the same interval: $\Delta'(t)=\frac{1}{\beta'}-\calF'(t)<\frac{1}{\beta'}-\calF'(t_*)<\calF'(t_0) - \calF'(t_*)<0$ which implies that $\hat{t}$ is the unique in $(t_*,t_0)$. Since $\Delta(t_0)<0$ and $\Delta(t)\to +\infty$ for $t\to +\infty$, we conclude that there exists at least a solution in $(t_0,+\infty)$. With the same argument used above one can prove that this solution, that we call $t_3(\beta')$, is unique. The solution $\hat{t}$ will be denoted as $t_2(\beta')$. We observe also that there exists a $\xi \in (t_2(\beta'),t_3(\beta'))$ such that $\calF'(\xi)=\frac{1}{\beta'}$, then since $\calF'(t)$ is decreasing in $(t_*,+\infty)$, we have: 
$\calF'(t_2(\beta'))>\calF'(\xi)=\frac{1}{\beta'}>\calF'(t_3(\beta'))$. As a consequence we have 
$r''(t_2(\beta'))=\frac{q-1}{q} \expec[W] (\calF'(t_2(\beta'))-\frac{1}{\beta'})>0$, which shows that 
$t_2(\beta')$ is a local minimum point and 
$r''(t_3(\beta'))=\frac{q-1}{q} \expec[W] (\calF'(t_3(\beta'))-\frac{1}{\beta'})<0$ which shows that 
$t_3(\beta')$ is a local maximum point for the pressure.
\end{proof}
\begin{remark}\label{rem-sol-bound-triple}
    Let us observe that at the boundary of the triple-intersection region the three solutions $t_1(\beta'), t_2(\beta'), t_3(\beta')$ still exists with $t_2(\frac{1}{\mathscr{F}'(t_0)})= t_3(\frac{1}{\mathscr{F}'(t_0)})=t_0$ and $t_1(\frac{q\expec[W]}{\expec[W^2]})= t_2(\frac{q\expec[W]}{\expec[W^2]}))=0$.
    
\end{remark}
For later use we study some properties of $t_3(\beta')$.
\begin{lemma}\label{lem-t3_cont_monot}
The function $\beta' \to t_3(\beta')$ is increasing and continuously differentiable in 
 $\beta' \in (\frac{1}{\mathscr{F}'(t_0)},\frac{q\expec[W]}{\expec[W^2]})$. Moreover, also the function $\beta' \to s(\beta')=\frac{t_3(\beta')}{\beta'}$ is increasing in the same interval. 
\end{lemma}
\begin{proof}
    Recall that for any  $b \in (\frac{1}{\mathscr{F}'(t_0)},\frac{q\expec[W]}{\expec[W^2]})$, $t_3(\beta')$ satisfies $G(b,t_3(b))=0$, where $G(\beta',t):=\mathscr{F}(t)-\frac{t}{\beta'}$. Since $\frac{\partial G}{\partial t}(b,t_3(b))=\calF'(t_3(b))-\frac{1}{b}<0$, by the Implicit Function Theorem there exists a unique continuous function $\tau(\beta')$ defined in a neighbourhood $U$ of $b$ such that $\tau(b)=t_3(b)$ and that satisfies to $G(\beta',\tau(\beta'))=0$ for $\beta' \in U$. Moreover, since $G$ is continuously  differentiable w.r.t. $\beta'$, $\tau(\beta')$ is also continuously differentiable in the same neighbourhood and, being $\tau'(b)=- [\frac{\partial G}{\partial t}(b,t_3(b))]^{-1}\cdot \frac{t_3(b)}{b^2}>0$, $\tau(\beta')$ is increasing in $b$. By uniqueness of the function defined implicitly, we have that $\tau(\beta')\equiv t_3(\beta')$ for $\beta' \in U$; this proves the lemma. The fact that $\calF'(t_3(\beta'))-\frac{1}{\beta'}<0$ for $\beta'$ inside the triple-intersection region is shown in the proof of Proposition \ref{prop:soln_inside_h}.
    
    Next we prove that $s(\beta')=\frac{t_3(\beta')}{\beta'}$ is increasing in $\beta'$, and then (see \eqref{betaprime}) in $\beta$. Since
    $$\frac{ds(\beta')}{d \beta'}= \frac{1}{\beta'^2} \left ( \frac{d t_3(\beta')}{d \beta'}\beta'-t_3(\beta')\right),$$
    we want to show that the second factor in the r.h.s. is positive.
    The derivative of $t_3(\beta')$, that has been computed above, is:   
    $$\frac{d t_3(\beta')}{d \beta'}= - \frac{t_3(\beta')}{\beta'^2(\calF'(t_3(\beta'))-\frac{1}{\beta'})}. $$
    Then, recalling that $\calF'(t_3(\beta'))-\frac{1}{\beta'}<0$,  we have
    $$\beta'\frac{d t_3(\beta')}{d \beta'}= \frac{t_3(\beta')}{\beta'(\frac{1}{\beta'}-\calF'(t_3(\beta')))}> t_3(\beta') $$
    since, being $\calF'(t)>0$, the denominator in the middle term is positive and less than 1. This shows that $s(\beta')$ is increasing in the triple intersection region.
\end{proof}
\begin{proposition}[Two solutions for $\beta'$ above the triple-intersection region]\label{prop:soln_above_h}
    For any $\beta'\geq \frac{q\expec[W]}{\expec[W^2]}$, let $t_3(\beta')$ be the largest solution to $\mathscr{F}(t)=\frac{t}{\beta'}$. Then $t_3(\beta')>0$, and further, the optimizer of \eqref{eq:laplace} has the form \eqref{eq:soln_form} (up to permutations), with $s(\beta)=\frac{t_3(\beta')}{\beta'}$.
\end{proposition}
\begin{proof}
    Recall $t_*$ is the inflection point of $\mathscr{F}(t)$. Further, the quantity $\left(\frac{q\expec[W]}{\expec[W^2]}\right)^{-1}$ is the slope of $\mathscr{F}(t)$ at $t=0$. So, if the slope $\frac{1}{\beta'}$ of the line $\frac{t}{\beta'}$ is less than that of the slope of $\mathscr{F}$ at $t=0$, (which is precisely the condition $\beta'\geq \frac{q\expec[W]}{\expec[W^2]}$) since $\mathscr{F}$ is convex in $(0,t_*)$, this means that the line $\frac{t}{\beta'}$ does not intersect $\mathscr{F}(t)$ in $(0,t_*)$. On the other hand, we note that $\mathscr{F}(t)<1$, and $\frac{t}{\beta'} \to 1$, as $t \to \beta'$. This means, there is a $\delta>0$, such that $\frac{\beta'-\delta}{\beta'}=\mathscr{F}(\beta'-\delta)$, and also let $\beta'-\delta$ be the smallest such point of intersection between $\frac{t}{\beta'}$ and $\mathscr{F}(t)$. As argued already, it must be that $\beta'-\delta>t_*$. Since the function $\mathscr{F}$ is concave on $[\beta'-\delta,\beta')$, the line $\frac{t}{\beta'}$ stays above of the curve $\mathscr{F}$ in the interval $(\beta'-\delta,\beta')$. Hence it must be that $\beta'-\delta$ is the unique non-zero intersection point of $\frac{t}{\beta'}$ and $\mathscr{F}(t)$. In particular, $t_3(\beta')=\beta'-\delta>0$. 
    Next, we claim that $t_3(\beta')$ is the global maximizer of $p(\frac{t}{\beta'})$. To see this, note that $p'(\frac{t}{\beta'}, \beta')=\expec{W}\left(\mathscr{F}(t)-\frac{t}{\beta'} \right)$, so the stationary points of $p(\frac{t}{\beta'}, \beta')$ are $t_1(\beta')=0$ and $t_3(\beta')$. Further, note that since $\mathscr{F}(t)$ is concave at $t_3(\beta')$, (as $t_3(\beta')>t_*$) where the line $\frac{t}{\beta'}$ intersects it, it must be that the slope $\mathscr{F}'(t_3(\beta'))$ of $\mathscr{F}$ at $t_3(\beta')$ is smaller than the slope $\frac{1}{\beta'}$ of the line $\frac{t}{\beta'}$. In other words, $p''(\frac{t_3(\beta')}{\beta'},\beta')=\expec{W}\left(\mathscr{F}'(t_3(\beta'))-\frac{1}{\beta'} \right)<0$, i.e. $t_3(\beta')$ is indeed the global maximum of $p(\frac{t}{\beta'},\beta')$.

    In particular, the optimizer of \eqref{eq:laplace} has the form \eqref{eq:soln_form} (up to permutations), with $s(\beta)=t_3(\beta')/\beta'$.
\end{proof}

\subsection{The optimal $s$ inside the triple-intersection region: identification of $\beta_c$}
Recall from Proposition \ref{prop:soln_inside_h} that when $\beta' \in (\frac{1}{\mathscr{F}'(t_0)},\frac{q\expec[W]}{\expec[W^2]})$, there are three solutions to $\mathscr{F}(t)=\frac{t}{\beta'}$: $0=t_1(\beta')<t_2(\beta')<t_3(\beta')$. For any $\beta' \in (\frac{1}{\mathscr{F}'(t_0)},\frac{q\expec[W]}{\expec[W^2]})$, we define 
\begin{align}
R(\beta')=\int_0^{t_3(\beta')} (\mathscr{F}(t)-\frac{t}{\beta'})dt. \label{eq:def_R}   \end{align}
Observe that, thanks to Remark \ref{rem-sol-bound-triple}, $R(\beta')$ is well defined also at the boundary points of the triple-intersection region.
\begin{lemma}[Comparing the pressure in $0$ and in $t_3(\beta')$]
Let $\beta_0'=1/{\mathscr{F}'(t_0)}$ and $\beta_1'=q\expec[W]/\expec[W^2]$.
\invisible{Define \[R(\beta_0')=\int_0^{t_0} (\mathscr{F}(t)-t\mathscr{F}'(t_0))dt,\] and \[R(\frac{})=\int_0^{t_3(\frac{q\expec[W]}{\expec[W^2]})} (\mathscr{F}(t)-\frac{t\expec[W^2]}{q\expec[W]})dt,\] where recall from Proposition \ref{prop:soln_above_h} that $t_3(\frac{q\expec[W]}{\expec[W^2]})$ is the largest solution to $\mathscr{F}(t)=\frac{\expec[W^2]t}{q\expec[W]}$. With these definitions,}
$\beta'\mapsto R(\beta')$ is a differentiable and monotone increasing function in $\beta'$ on the compact set $[\beta_0',\beta_1']$, with $R(\beta_0')<0$ and $R(\beta_1')>0$. Therefore $R(\beta')$ has a unique root $\beta'_c \in (\beta_0',\beta_1')$.
\end{lemma}
\Gib{I (Claudio) changed slightly the statement of the lemma asserting that $R(\beta')$ is in fact differentiable. A proof (a bit shorter than the original one) follows. }
\begin{proof}
Observe that, being $t_3(\beta')$ differentiable by Lemma \ref{lem-t3_cont_monot}, $R(\beta')$ is differentiable because the integrand in \eqref{eq:def_R} is also differentiable  
w.r.t. $\beta'$. Then we compute the derivative
$$
\frac{d R(\beta')}{d \beta'}=
\int_0^{t_3(\beta')} \frac{t}{{\beta'}^2}\, dt + \left (\mathscr{F}(t_3(\beta'))-\frac{t_3(\beta')}{\beta'}\right )=\frac{t_3(\beta')^2}{2 {\beta'}^2}>0,$$
which shows that $R(\beta')$ is increasing on  $[\beta_0',\beta_1']$.\\
Further, we claim that $R(\beta_0')<0$. To see this, recall from the proof of Proposition \ref{prop:sol_lim_vec_below_CR} that the tangent line $G_2(t)=t\mathscr{F}'(t_0)$ completely lies above the graph of $\mathscr{F}(t)$ in the interval $[0,t_0]$, with the only points of intersections being $(0,0)$ and $(t_0,\mathscr{F}(t_0))$. Next, we claim that $R(\beta_1')>0$. The argument is similar. First we note that the tangent line $G_1(t)=\frac{t}{\beta_1'}$ is below the graph of $\mathscr{F}$ in the interval $[0,t_*]$. This is because in this interval $\mathscr{F}$ is convex, and $G_1(t)$ is a tangent line to the graph of $\mathscr{F}$ at a point inside this interval. 
Consider the line $L(t)$ joining the points $(t_*,\mathscr{F}(t_*))$ and $(t_3(\beta_1'),\mathscr{F}(t_3(\beta_1')))$. Due to concavity of $\mathscr{F}$ in the interval $[t_*,t_3(\beta_1')]$, the line $L(t)$ is strictly below the graph of $\mathscr{F}$. Since $G_1(t_*)<L(t_*)=\mathscr{F}(t_*)$, this means the line $G_1(t)$, which is also the line joining the points $(t_*,G_1(t_*))$ and $(t_3(\beta_1'),\mathscr{F}(t_3(\beta_1')))$ is below the graph $\mathscr{F}$ in the interval $[t_*,t_3(\beta_1'))]$. Combining everything, we get that the line $G_2(t)$ is always below the graph of $\mathscr{F}$ in the interval $[0,t_3(\beta_1'))]$. Consequently, $R(t_3(\beta_1')))>0$. 
By the intermediate value theorem and the monotonicity, there exists a unique $\beta_c' \in (\beta_0',\beta_1')$ such that $R(\beta_c')=0$.
\end{proof}

\invisible{
\Gib{The following is the former proof of Lemma 5.10 (restricted to continuity).}
\begin{proof}
Note that for fixed $t$, $\mathscr{F}(t)-t/\beta'$ is a monotone increasing continuous function of $\beta'$ on $[\beta_0',\beta_1']$. Hence if we can show that $t_3(\beta')$ is also continuous and monotone increasing as a function of $\beta'$, with the definition $t_3(\beta_0')=t_0$, we can conclude $R(\beta')$ is also continuous in the same domain.  

This is because, for $\beta_1 < \beta_2$, where $\beta_1,\beta_2 \in [\beta_0',\beta_1']$, one can write 
\begin{align*}
    |R(\beta_1)-R(\beta_2)|=&\left|\int_{0}^{t_3(\beta_1)}(\mathscr{F}(t)-t/\beta_1)dt-\int_{0}^{t_3(\beta_2)}(\mathscr{F}(t)-t/\beta_2)dt \right|\\&\leq \int_{0}^{t_3(\beta'_1)}\mathbf{1}_{[t_3(\beta_1),t_3(\beta_2)]}\left(|\mathscr{F}(t)-t\mathscr{F}'(t_0)|\vee |\mathscr{F}(t)-\frac{t\expec[W^2]}{q\expec[W]}| \right)dt.
\end{align*}
In particular, bounding from above the continuous function $\left(|\mathscr{F}(t)-t\mathscr{F}'(t_0)|\vee |\mathscr{F}(t)-\frac{t\expec[W^2]}{q\expec[W]}| \right)$ of $t$ in the compact interval $[0,t_3(\frac{q\expec[W]}{\expec[W^2]})]$ by its maximum value, we note that when $\beta_2-\beta_1 \to 0$, the r.h.s.\ above, and hence $|R(\beta_1)-R(\beta_2)|$ also, tend to zero.
{The proof that $t_3(\beta')$ is continuous and increasing is a tedious analysis argument, and can be found in Section \ref{sec-t3-mon}.}


Note that from the definition of $R$, it is clear that it is increasing in $\beta'$ if $t_3$ is so. Further, we claim that $R(\frac{1}{\mathscr{F}'(t_0)})<0$. To see this, recall from the proof of Proposition \ref{prop:sol_lim_vec_below_CR} that the tangent line $G_2(t)=t\mathscr{F}'(t_0)$ completely lies above the graph of $\mathscr{F}(t)$ in the interval $[0,t_0]$, with the only points of intersections being $(0,0)$ and $(t_0,\mathscr{F}(t_0))$. Next, we claim that $R(\frac{q\expec[W]}{\expec[W^2]})>0$. The argument is similar. First we note that the tangent line $G_1(t)=\frac{t\expec[W^2]}{q\expec[W]}$ is below the graph of $\mathscr{F}$ in the interval $[0,t_*]$. This is because in this interval $\mathscr{F}$ is convex, and $G_1(t)$ is a tangent line to the graph of $\mathscr{F}$ at a point inside this interval. 

Consider the line $L_2(t)$ joining $(t_*,\mathscr{F}(t_*))$ and $(t_3(\frac{q\expec[W]}{\expec[W^2]}),\mathscr{F}(t_3(\frac{q\expec[W]}{\expec[W^2]})))$. Due to concavity of $\mathscr{F}$ in the interval $[t_*,t_3(\frac{q\expec[W]}{\expec[W^2]})]$, the line $L_2(t)$ is strictly below the graph of $\mathscr{F}$. Since $G_1(t_*)<L_2(t_*)=\mathscr{F}(t_*)$, this means the line $G_1(t)$, which is also the line joining $(t_*,G_1(t_*))$ and $(t_3(\frac{q\expec[W]}{\expec[W^2]}),\mathscr{F}(t_3(\frac{q\expec[W]}{\expec[W^2]})))$ is below the graph $\mathscr{F}$ in the interval $[t_*,t_3(\frac{q\expec[W]}{\expec[W^2]}))]$. Combining everything, we get that the line $G_2(t)$ is always below the graph of $\mathscr{F}$ in the interval $[0,t_3(\frac{q\expec[W]}{\expec[W^2]}))]$. Consequently, $R(t_3(\frac{q\expec[W]}{\expec[W^2]})))>0$. Since $R$ is also monotone increasing, as proved earlier, this means by the intermediate value theorem, $R$ has a unique root $\beta_c'\in (\beta_0',\beta_1')$. 
In particular, writing out $R(\beta_c')=0$, we note that $\beta'_c$ solves
\begin{align}
    \beta'_c=\frac{(t_3(\beta_c))^2}{2\int_0^{t_3(\beta'_c)}\mathscr{F}(t)dt}.
\end{align}

\end{proof}

}

\begin{proposition}[$\beta_c'$ is the critical value]
    Let $\beta' \in (\beta_0',\beta_1')$. When $\beta'\in (\beta_0',\beta'_c)$, the optimizer of \eqref{eq:laplace} has the form \eqref{eq:soln_form}, with $s(\beta')=0$. When $\beta'\in [\beta'_c,\beta_1')$, the optimizer of \eqref{eq:laplace} has the form \eqref{eq:soln_form}, with $s(\beta)=\frac{t_3(\beta')}{\beta'}$. At $\beta'=\beta'_c$, both $t_1(\beta'_c)=0$ and $t_3(\beta'_c)>0$ are global maximizers of $p(\frac{t}{\beta'_c}, \beta_c)$. In particular, $p(0, \beta)=p(\frac{t_3(\beta'_c)}{\beta'_c},\beta)$.  
\end{proposition}
\begin{proof}
    Recall $p_\beta'(s)=\mathscr{F}(s\beta')-s$. Note that from the previous proof $R(\beta')<0$ for $\beta \in [\beta_0',\beta_c')$, and $R(\beta')>0$ for $\beta \in (\beta_c',\beta_1']$. Consequently, since $R(\beta')=\int_{0}^{t_3(\beta')}p'(\frac{t}{\beta'},\beta)dt$, we have $p(\frac{t_3(\beta')}{\beta'})<p(0)$ for $\beta \in [\beta_0',\beta_c')$, and $p(\frac{t_3(\beta')}{\beta'},\beta)>p(0,\beta)$ for $\beta \in (\beta_c',\beta_1']$. Since from Proposition \ref{prop:soln_inside_h} we know that $0=t_1(\beta')$ and $t_3(\beta')$ are the only local maximum points of $p(\frac{t}{\beta'},\beta)$, we note that in the interval $[\beta_0',\beta_c')$, the global maximizer of $p(\frac{t}{\beta'},\beta)$ is $0=t_1(\beta')$, while in the interval $(\beta_c',\beta_1']$, the global maximizer of $p(\frac{t}{\beta'},\beta)$ is $t_3(\beta')$. This proves the first part of the statement of the Proposition. Further, recalling that $R(\beta_c')=0$, we note that $p(0,\beta)=p(\frac{t_3(\beta_c')}{\beta'_c},\beta)$, which means that for $\beta=\beta'_c$, both $0=t_1(\beta'_c)$ and $t_3(\beta'_c)>0$ are global maximizers of $p(\frac{t}{\beta'_c},\beta_c)$.
\end{proof}
}

\section{The order parameter: Proof of Theorem \ref{thm-proportions}}
\label{sec-proportions}

\begin{proof}[Proof of Theorem \ref{thm-proportions}]
\invisible{Recall that $\by^\star(\beta)$ denotes the optimizer of \eqref{var-problem-psi} that is of the form in \eqref{eq:soln_form}. Further, we let $\mathscr{Y}_\beta$ be the set of all permutations of the vector $\by^\star(\beta)$ in \eqref{eq:soln_form}. By symmetry, each of the solutions in $\mathscr{Y}_\beta$ is an optimizer of \eqref{eq:laplace}. Let $\mathscr{Y}_\beta(\vep)$ be the $\vep$-neighborhood or $\mathscr{Y}_\beta$. As proved above, recall also \eqref{x-k-form}, for all $\vep>0$, 
    \eqn{
    \max_{\by \not\in \mathscr{Y}_\beta(\vep)} \left [  {F_\beta}(\boldsymbol{0}; \by)    -\frac{\beta'}{2 \mathbb{E}[W]}  \sum_{i=1}^q y_i^2 \right ]
    <\sup_{s\in[0,1]} p_\beta(s).
    }
Recall that $\mathscr{X}_\beta$ denotes the set of all permutations of the vector $\bx^\star(\beta)$ defined in and below \eqref{x1-def}, and that $\mathscr{X}_\beta(\vep)$ denotes the $\vep$-ball around $\mathscr{X}_\beta$. We investigate
    \eqn{
    \mu_{n,\beta,0}\Big((X_k(\beta))_{k\in[q]}\not \in \mathscr{X}_\beta(\vep)\Big),
    }
and show that it is exponentially small for every $\vep>0$.

\RvdH{TO DO 4: Now we need to go to proportions!}

\RvdH{Below is a new try.}}

We consider the weighted color proportions $Y_k=Y_k(\bsigma)$ as defined in \eqref{Y-k-number-def}.
\invisible{We prove that $(Y_k/n)_{k\in [q]}$ undergoes a phase transition, in that it behaves differently when $s(\beta,0)=0$ versus $s(\beta,0)>0.$}
We fix $k=1$, and consider the moment generating function of $nY_1$. Let
    \eqn{
    \expec[Z_n(\beta,B;t)]
    =\expec\Bigg[\sum_{\bs \in \omqn} \exp \Big[ \beta \sum_{i<j}  I_{i,j} \indic{\s_i=\s_j} + B\sum_{i\in[n]} \indic{\s_i=1}+tnY_1\Big ]\Bigg].
    }
Then, following the analysis in Section \ref{sec-Hub-Strat}, we obtain that
    \eqan{
    \varphi_{\beta,B}(t)\equiv\lim_{n\to \infty} \frac 1 n \log \expec[Z_n(\beta,B;t)]
& = \max_{\by\in \R^q} \left [  F_{\beta,B}(t; \by)    -\frac{\beta'}{2 \mathbb{E}[W]}  \sum_{k=1}^q y_k^2 \right ], 
    }
where now
    \eqn{
    F_{\beta,B}(t; \by)=
    \expec\left [ \log \left ( \e^{(\beta' y_1 +t)  \frac{W}{\expec[W]}+B}+\cdots + \e^{ \beta' y_q  \frac{W}{ \expec[W]}} \right )\right ].
    }
Replace $y_1+t/\beta'$ by $y_1 $ to obtain
    \eqan{
    \varphi_{\beta,B}(t)
& = \max_{\by\in\R^q} \left [F_{\beta,B}(0; \by)    -\frac{\beta'}{2 \mathbb{E}[W]}  \sum_{k=1}^q (y_k-t/\beta'\indic{k=1})^2 \right ]\nn\\
&=\max_{\by\in\R^q} \left [  {F_{\beta,B}}(0; \by)    -\frac{\beta'}{2 \mathbb{E}[W]}  \sum_{k=1}^q y_k^2
+\frac{1}{\expec[W]}(ty_1-(t^2/(2\beta'))\right ].
\label{VP-weighted-spins}
    }
Let $\by^\star(\beta,B; t)$ be an optimizer of the optimization problem. Then, 
    \eqan{
    \varphi_{\beta,B}(t)&= F_{\beta,B}(0; \by^\star(\beta,B;t))    -\frac{\beta'}{2 \mathbb{E}[W]}  \sum_{i=1}^q (y_i^\star(\beta,B;t))^2
+\frac{1}{\expec[W]}(ty_1^\star(\beta,B;t)-(t^2/(2\beta')).
\label{VP-weighted-spins-2}
    }
We note that, defining for all $t\geq 0$,
    \eqn{
    \label{varphi-beta-t}
    \lim_{n\rightarrow \infty}\frac{1}{n}\log\expec_{\mu_{n;\beta,B}}[\e^{tY_1}]
    =\psi_{\beta,B}(t),
    }
we have
    \eqn{
    \psi_{\beta,B}(t)=\varphi_{\beta,B}(t)-\varphi_{\beta,B}(0).
    \label{phi-t-0}
    }
Since the function to be maximized in \eqref{VP-weighted-spins} is continuous in $t$ and bounded in $\by$, as $t\searrow 0$, every subsequence of $\by^\star(\beta,B;t)$ contains a further converging subsequence. Thus, when the optimizer $\by^\star(\beta,B)$ for $t=0$ is unique, $\by^\star(\beta,B;t)$ converges as $t\searrow 0$ to it. Therefore, using that $\varphi_{\beta,B}(t)$ is a supremum and replacing  $\by^\star(\beta,B;t)$ in \eqref{VP-weighted-spins-2} with $\by^\star(\beta,B)$, we obtain from \eqref{phi-t-0} that
    \eqn{\label{lower-b-psi}
    \frac{1}{t}\psi_{\beta,B}(t)\geq \frac{y_1^\star(\beta,B)}{\expec[W]}-t/(2\beta'\expec[W]).
    } 

Further, also $\varphi_{\beta,B}(0)$ is a supremum, which is obtained from \eqref{VP-weighted-spins} by computing the r.h.s.\ for $\by^\star(\beta,B)$ and setting $t=0$. If, instead, we substitute $\by^\star(\beta,B;t)$, we obtain a lower bound for $\varphi_{\beta,B}(0)$ that, by \eqref{phi-t-0}, yields
    \eqn{\label{upper-b-psi}
    \frac{1}{t}\psi_{\beta,B}(t)\leq 
    \frac{y_1^\star(\beta,B;t)}{\expec[W]}-t/(2\beta'\expec[W]).
    }
Thus, from \eqref{lower-b-psi} and \eqref{upper-b-psi}, it follows that 
$$\frac{d\psi_{\beta,B}(0)}{dt}= \frac{y_1^\star(\beta,B)} {\expec[W]},$$ i.e., \eqref{psi-mGF} follows.

\invisible{When differentiating w.r.t.\ $t$, we note that the term multiplying the $\frac{d}{dt}y_i^\star(\beta,t)$ equals zero, since $\by^\star(\beta,t)$ is an optimizer. Therefore,
    \eqn{
    \lim_{n\to \infty} \frac{d}{dt}\frac 1 n \log \expec[Z_n(\beta,t)]
    =\frac{y_1^\star(\beta,t)}{\expec[W]}-\frac{t}{2\expec[W]\beta'}.
    }
\RvdH{We need a proof here!}
We apply this to $t=0$, and using that then $\by^\star(\beta,0)=\by^\star(\beta)$. Thus, when $\mathscr{Y}_\beta$ only contains one element, 
    \eqan{
    \lim_{n\to \infty}\mu_{\beta,B}\Big(Y_i/n\Big)=\lim_{n\to \infty} \frac{d}{dt}\frac 1 n \log \expec[Z_n(\beta,t)]\Big|_{t=0}
    =\frac{y_1^\star(\beta)}{\expec[W]},
    }
which equals $1/q$ since $y_1^\star(\beta)+\cdots+y_q^\star(\beta)=\expec[W]$, and they are all the same. This shows that the weighted spin averages are all the same when $\mathscr{Y}_\beta$ contains only one element.}
\smallskip

\invisible{\RvdH{Are we here not implicitly assuming that the optimizer over $s$ is unique when $s(\beta)>0$?}

For the other case, when $\mathscr{Y}_\beta$ contains $q$ optimizers, let $y_1^\star(\beta)$ be the largest element of one such optimizer. We have just proved that there exists a unique such largest value, and it is the same for all the optimizers. Let $\by^\star(\beta)$ be the corresponding (now unique) optimizer for which the first element is the largest, and all other elements are the same smaller value (see the proof of Theorem \ref{thm-soln_form-positive}). Then, as before, substituting $\by=\by^\star(\beta)$ in the optimization problem in \eqref{VP-weighted-spins} provides a lower bound, so that 
    \eqan{
    \frac{1}{t}\varphi_\beta(t)
    &\geq
\frac{y_1^\star(\beta)}{\expec[W]}-t/(2\beta'\expec[W]).
    }
Thus, 
    \eqn{
    \liminf_{t\searrow 0} \frac{1}{t}\varphi_\beta(t)
    \geq \frac{y_1^\star(\beta)}{\expec[W]}>\frac{y_i^\star(\beta)}{\expec[W]},
    }
for all $i\in \{2, \ldots, q\}$. Thus, since $y_1^\star(\beta)+\cdots+y_q^\star(\beta)=\expec[W],$
    \eqn{
    \liminf_{t\searrow 0}\frac{1}{t}\varphi_\beta(t)>1/q.
    }
This shows that there is a phase transition in the weighted sum of spins, depending on whether $s(\beta)>0$ or not. The argument suggests that the weighted sum of spins is close to an element in $\mathscr{Y}_\beta$. Further, we note that, again as above
    \eqn{
    \liminf_{t\searrow 0} \frac{1}{t}\varphi_\beta(t)
    \leq \liminf_{t\searrow 0} \frac{y_1^\star(\beta,t)}{\expec[W]}.
    }
Since this has to be {\em at least} $y_1^\star(\beta)/\expec[W]$, and the sequence $\by^\star(\beta,t)$ converges as $t\searrow 0$ to one of the maximizers $\by^\star(\beta,0)=\by^\star(\beta)$ for $t=0$, we have that $y_1^\star(\beta,t)$ has to converge to the largest coordinate of the each of the $q$ maximizers. Thus,
    \eqn{
    \lim_{t\searrow 0} \frac{1}{t}\varphi_\beta(t)
    =\frac{y_1^\star(\beta)}{\expec[W]},
    }
where $y_1^\star(\beta)$ is the largest coordinate of any maximizer. Thus, \eqref{s(beta)>0-varphi} follows.}
\smallskip

For \eqref{der-x1-equation}-\eqref{x1-equation}, we note that, by convexity of $B\mapsto \varphi_{n}(\beta,B)\equiv \frac1n \log\expec[Z_n(\beta,B)],$ for every $t>0$,
    \eqn{
    \frac{1}{n}\expec_{\mu_{n;\beta,B}}\left[\sum_{i=1}^n\indic{\s_i=1}\right] \leq  \frac{1}{t}[\varphi_{n}(\beta,B+t)-\varphi_{n}(\beta,B)],
    } 
    since the l.h.s. coincides with $\frac{\partial}{\partial B}\varphi_{n}(\beta,B)$.
Letting $n\rightarrow \infty$ gives that
    \eqn{
    \limsup_{n\rightarrow \infty} \frac{1}{n}\expec_{\mu_{n;\beta,B}}\left[\sum_{i=1}^n\indic{\s_i=1}\right] \leq  \frac{1}{t}[\varphi(\beta,B+t)-\varphi(\beta,B)],
    }
and letting $t\searrow 0$, we obtain
    \eqn{
    \limsup_{n\rightarrow \infty} \frac{1}{n}\expec_{\mu_{n;\beta,B}}\left[\sum_{i=1}^n\indic{\s_i=1}\right] \leq  \frac{\partial}{\partial B}\varphi(\beta,B).
    }
This proves the upper bound in \eqref{der-x1-equation} for $B>0$.

For the lower bound in \eqref{der-x1-equation} for $B>0$, we follow exactly the same proof, now starting, for every $t>0$, from 
    \eqn{
    \frac{1}{n}\expec_{\mu_{n;\beta,B}}\left[\sum_{i=1}^n\indic{\s_i=1}\right] \geq  \frac{1}{t}[\varphi_{n}(\beta,B)-\varphi(\beta,B-t)].
    }
For \eqref{x1-equation}, we take $\bB=(B,0, \ldots,0)$, and differentiate \eqref{eq:laplace} w.r.t.\ $B$. Using the fact that the derivatives of $\by^\star(\beta,B)$ w.r.t.\ $B$ do not contribute by the stationarity condition \eqref{eq:stat0}, we obtain that
    \be
    \label{x1-equation-rep}
    x^\star_1(\beta,B) = \frac{\partial}{\partial B}F_{\beta,\bB}( \by^\star(\beta,B)),
    \ee
which, by \eqref{F-pressure}, indeed gives \eqref{x1-equation}.
\end{proof}

\invisible{\RvdH{To complete this argument, we would have to show that 
    \eqn{
    \lim_{n\to \infty} \frac{d}{dt}\frac 1 n \log \expec[Z_n(\beta,t)]\Big|_{t=0}
    =\frac{y_1^\star(\beta)}{\expec[W]}.
    }}}


\invisible{\section{Critical temperature: Proof of Theorem \ref{thm-crit-value}}
In this section, we identify the critical value.}

\section{Zero-crossing condition: Proof of Theorems \ref{thm-critical-value}, \ref{thm-density-Guido}, \ref{thm-Pareto} and \ref{thm-two-atoms}}
\label{sec-unique-zero-crossing-or-not}
In this section, we investigate the zero-crossing condition in Condition \ref{cond-zer-cross}. In Section \ref{sec-critical-value}, we identify the critical value and prove  Theorem \ref{thm-critical-value}. In Section \ref{sec-unique-zero-crossing}, we 
investigate the special examples in \eqref{Guido(3)}, and prove Theorem \ref{thm-density-Guido}. We then zoom in into the Pareto case in Section \ref{sec-Pareto-(3,4)}, where we prove Theorem \ref{thm-Pareto}.

\subsection{Critical value for the zero field case: Proof of Theorem \ref{thm-critical-value}}
\label{sec-critical-value}
In this section, we identify the critical value when $B=0$ and  Condition \ref{cond-zer-cross} holds.

\invisible{\begin{proof}[Proof of Theorem \ref{thm-critical-value}]
We compute that 
    \eqn{
    \frac{p_{\beta,0}(t/\beta')-p_{\beta,0}(0)}{\expec[W]}=
    \frac{1}{\expec[W]}\expec\left[\log{(\e^{tW}+q-1)}\right]
    -\frac{q-1}{2q} \frac{t^2}{\beta'}-\frac{t}{q}.
    }
We substitute
    \eqn{
    \mathscr{F}_0(t_c)=\frac{t_c}{\beta_c'},
    }
which holds by the stationarity condition \eqref{FOPT-1-B1-B>0}.
Then, we arrive at the fact that $t_c$ solves
    \eqn{
    \mathscr{K}(t_c)=0,
    }
where
    \eqn{
    \mathscr{K}(t)=
    \frac{1}{\expec[W]}\expec\left[\log{(\e^{tW}+q-1)}\right]
    -\frac{q-1}{2q} t\mathscr{F}_0(t)-\frac{t}{q}.
    }
We note that 
    \eqn{
    \mathscr{K}(0)=\mathscr{K}'(0)=\mathscr{K}''(0)=0,
    }
while $\mathscr{K}'''(0)<0.$ Further,
    \eqn{
    \lim_{t\rightarrow \infty} \mathscr{K}(t)=\infty.
    }
We compute that
    \eqn{
    \label{kprime-funct}
    \mathscr{K}'(t)=\frac{q-1}{2q} [\mathscr{F}_0(t)-t\mathscr{F}'_0(t)],
    }
and 
    \eqn{
    \label{Kdoubleprime-funct}
    \mathscr{K}''(t)
    =-\frac{q-1}{2q} t\mathscr{F}''_0(t).
    }
Thus, if $\mathscr{F}_0$ satisfies the unique zero-crossing condition in Condition \ref{cond-zer-cross}, then so does $\mathscr{K}$. In particular, $t\mapsto \mathscr{K}(t)$ will first be negative, and then become positive.
We conclude that $\mathscr{K}(t_c)=0$ has a unique solution. 

Since $t=s\beta'$, and $\mathscr{F}_0(t_c)=t_c/\beta_c'$, the two equalities for $\beta_c$ and $s(\beta_c)$ follow. 
\end{proof}

\Gib{Claudio: I have slightly rewritten (see below) the proof in order to better highlight the fact that in this case one of the two maxima is always $s=0$. Moreover the factor $1/q$ was lacking in $\log()$ of equation (7.1). Also, the details (7.5) and (7.6) seem unnecessary in the proof.}}

\begin{proof}[Proof of Theorem \ref{thm-critical-value}]
Recall from Theorems \ref{thm-nature-pt-a} and \ref{thm-nature-pt-b} (see also Lemma \ref{lemma:nature-stat-points-bis}) that the phase transition occurs at a $\beta_c$ at which the {\em criticality condition} $p_{\beta_c,0}(s_1(\beta_c,0))=p_{\beta_c,0}(s_3(\beta_c,0))$ holds. Here, by Lemma \ref{lemma:nature-stat-points-bis}, $s_1(\beta_c,0)$ and $s_3(\beta_c,0)$ are the local maxima of $p_{\beta_c,0}$, with $s_1(\beta,0)\equiv 0$ when $B=0$.   
Thus, in order to impose the {\em criticality condition} we compute that
    \eqn{
    \frac{p_{\beta,0}(t/\beta')-p_{\beta,0}(0)}{\expec[W]}=
    \frac{1}{\expec[W]}
    \expec
    \left[
    \log \left(\frac{\e^{tW}+q-1}{q}\right)\right]
    -\frac{q-1}{2q} \frac{t^2}{\beta'}-\frac{t}{q}.
    }
Further, in order to incorporate the {\em stationarity condition} \eqref{FOPT-1-B1-B>0}, we substitute $\mathscr{F}_0(t)=t/\beta'$ in the previous equation, 
since this relation holds at the critical point $t_c$, i.e., $\mathscr{F}_0(t_c)={t_c}/{\beta_c'}$.
Then, we arrive at the fact that $t_c$ solves
    \eqn{
    \mathscr{K}(t_c)=0,
    }
where
    \eqn{\label{kappa}
    \mathscr{K}(t)=
    \frac{1}{\expec[W]}\expec\left[\log \left (\frac{\e^{tW}+q-1}{q}\right )\right]
    -\frac{q-1}{2q} t\mathscr{F}_0(t)-\frac{t}{q}.
    }
\invisible{We note that 
    \eqn{
    \mathscr{K}(0)=\mathscr{K}'(0)=\mathscr{K}''(0)=0,
    }
while $\mathscr{K}'''(0)<0.$ 
We observe that
    \eqn{
    \lim_{t\rightarrow \infty} \mathscr{K}(t)=\infty.
    }
    }
We also observe that
    \eqn{
    \label{kprime-funct}
    \mathscr{K}'(t)=\frac{q-1}{2q} [\mathscr{F}_0(t)-t\mathscr{F}'_0(t)],
    }
and 
    \eqn{
    \label{Kdoubleprime-funct}
    \mathscr{K}''(t)
    =-\frac{q-1}{2q} t\mathscr{F}''_0(t).
    }
\invisible{Thus, since $\mathscr{F}_0$ satisfies the unique zero-crossing condition in Condition \ref{cond-zer-cross}, then so does $\mathscr{K}(t)$. In particular, $t\mapsto \mathscr{K}(t)$ is first negative and then positive.
We conclude that $\mathscr{K}(t)=0$ has a unique solution $t_c$.}
Recall that $\mathscr{F}_0(0)=0$ and that $\mathscr{F}^{''}_0(t)$ is first positive and then negative. It thus follows from \eqref{kappa},\eqref{kprime-funct} and \eqref{Kdoubleprime-funct}
that $ \mathscr{K}(0)=\mathscr{K}^{'}(0)=\mathscr{K}^{''}(0)=0$, and that $\mathscr{K}^{''}(t)$, and hence $\mathscr{K}(t)$, is negative for small $t>0$. Furthermore, we have from \eqref{kappa} that $\mathscr{K}(t)$ is positive for large $t>0$. Since $\mathscr{F}_0(t)$, and hence $\mathscr{K}(t)$, has a unique inflection point $t_*$, we conclude that $\mathscr{K}(t)$ has a unique solution $t_c$.
Further, since $t=s\beta'$, and $\mathscr{F}_0(t_c)=t_c/\beta_c'$, the two equalities \eqref{beta-s-crit} for $\beta_c$ and $s(\beta_c,0)$ follow. 
\end{proof}


\subsection{Zero-crossing: Proof of Theorem \ref{thm-density-Guido}}
\label{sec-unique-zero-crossing}
For $\e^{-B}(q-1)>1$, recall from \eqref{eq:second_derivative_gen} that we need to consider
	\eqn{
	\label{Guido(1)}
	\frac{d^{2}}{d t^{2}} \mathscr{F}_B(t)=\frac{q}{\mathbb{E}[W]}
	\mathbb{E}\left[W^{3} \frac{\e^{-B}(q-1) \e^{t W}-\e^{2 t W}}{\left(\e^{t W}+\e^{-B}(q-1)\right)^{3}}\right],
	}
where $W$ is a non-negative random variable with $\mathbb{E}[W]>0$, and there are $\tau>2, c>0$ such that \eqref{main-assumption-Fn} holds, and whose density has the form in \eqref{Guido(3)}. We note that \eqref{Guido(1)} only depends on $B$ and $q$ through $\e^{-B}(q-1)$. We will, below, thus replace $\e^{-B}(q-1)$ by $q-1$, and simply assume that $q-1>1$.
We thus concentrate on $\mathscr{F}(t)\equiv \mathscr{F}_0(t)$ with $q>2$.


\invisible{
	\eqn{
	\label{main-assumption-Fn}
	1-\prob(W \leq w) \leq \frac{D}{w^{\tau-1}},\qquad w \geq 1.
	}
We assume that $W$ has a pdf $g(w)$  of the form
	\eqn{
	\label{Guido(3)}
	g(w)=C w^{p} \e^{-\varphi(w)} \chi_{[b, c)}(w), \qquad w \geq 0,
	}
with $p \in \mathbb{R}, \varphi(w) \geq 0,0 \leq b<c \leq \infty, C \geq 0$ such that $\mathbb{E}[W]>0$ and for which \eqref{main-assumption-Fn} holds.}
\smallskip

With the density $f_{\sss W}(w)$ of $W$ of the form \eqref{Guido(3)}, we can write
	\eqan{
	\label{Guido(4)}
	\frac{d^2}{d t^2} \mathscr{F}(t)&=\frac{q C}{\mathbb{E}[W]} \int_{0}^{\infty} w^{3} \frac{(q-1) \e^{t w}-\e^{2 t w}}{\left(\e^{t w}+q-1\right)^{3}} w^{p} \e^{-\varphi(w)} \indicwo{[b,c)}(w)dw\nn\\
	&=\frac{q C}{t^{p+4} \mathbb{E}[W]} \int_{b t}^{c t} x^{p+3} \e^{-\varphi(x/t)} a(x) d x,
	}
where
	\eqn{
	\label{Guido(5)}
	a(x)=\frac{(q-1) \e^{x}-\e^{2 x}}{\left(\e^{x}+q-1\right)^{3}}, \qquad x \geq 0.
	}

We are particularly interested in densities $f_{\sss W}(w)$  in \eqref{Guido(3)} such that $\frac{d^{2}}{d t^{2}} \mathscr{F}(t)$ has exactly one positive zero $t_{0}$. For this we must consider
\be
	\Phi(t)=\Phi(t ; b, c, r ; \varphi)=\int_{b t}^{c t} x^{r} \e^{-\varphi(x/t)} a(x) d x,
\ee
since the last line of \eqref{Guido(4)} implies that
\eqn{\label{eq:Der2FPhi}
\frac{d^{2}}{d t^{2}} \mathscr{F}(t)=\frac{q C}{t^{p+4} \mathbb{E}[W]}\, \Phi(t),}
with $r=p+3$. 
We split this analysis into several cases, depending on whether $b=0$ or $b>0$, and $c<\infty$ or $c=\infty$. The proofs will be similar, with specific adaptations for the various cases. 
\medskip


\paragraph{\bf Case $0<b<c<\infty$.} We will prove the following result:

\begin{proposition} 
\label{prop-Guido-1}
Let $r \in \mathbb{R}, 0<b<c<\infty$. Assume that $\varphi \in C^{2}([b, c])$ and that $w\varphi^{\prime}(w)$ is non-negative and non-decreasing in $w \in[b, c]$. Then there is exactly one $t_{0}\in (0, \infty)$ such that $\Phi\left(t_{0} ; b, c, r ; \varphi\right)=0$.
\end{proposition}
\medskip

\paragraph{\bf Case $0<b<c=\infty$.}  We will prove the following result:

\begin{proposition}
\label{prop-Guido-3}
Let $r \in \mathbb{R}, 0<b<c=\infty$. Assume that $\varphi \in C^{2}([b, \infty))$ and that w $\varphi^{\prime}(w)$ is non-negative and non-decreasing in $w\in [b, \infty)$, with $\lim _{w \rightarrow \infty} w\varphi^{\prime}(w)=\infty$. Then there is exactly one $t_{0}\in(0,\infty)$ such that $\Phi\left(t_{0} ; b, \infty, r ; \varphi\right)=0$.
\end{proposition}

The condition $w \varphi^{\prime}(w) \rightarrow \infty$ as $w \rightarrow \infty$, assures that $\exp (-\varphi(w)$ ) $=O\left(w^{-\tau}\right)$ as $w \rightarrow \infty$, for all $\tau>2$. Hence \eqref{main-assumption-Fn} holds for the density $f_{\sss W}(w)$ in \eqref{Guido(3)}.
\medskip

\paragraph{\bf Case $0=b<c<\infty$.}  We will prove the following result:

\begin{proposition} 
\label{prop-Guido-2}
Let $r=p+3>2,0=b<c<\infty$. Assume that $\left.\varphi \in C^{2}\left([0, c\right]\right)$ and that $w\varphi^{\prime}(w)$ is non-negative and non-decreasing in $w \in[0, c]$. Then there is exactly one $t_{0}\in(0, \infty)$ such that $\Phi(t ; 0, c, r ; \varphi)=0$.
\end{proposition}

Since $b=0$, we have to restrict $p$ in the form \eqref{Guido(3)} of the density $f_{\sss W}(w)$ to the range $p>-1$.
\medskip

\paragraph{\bf  Case $0=b<c=\infty$.}  We will prove the following result:

\begin{proposition}
\label{prop-Guido-4}
Let $r=p+3>2, b=0, c=\infty$. Assume that $\varphi \in C^{2}([0, \infty))$ and that $w\varphi^{\prime}(w)$ is non-negative and non-decreasing in $w \in[0, \infty)$ with $\lim _{w \rightarrow \infty} w \varphi^{\prime}(w)=\infty$. Then there is exactly one $t_{0}\in(0,\infty)$ such that $\Phi\left(t_{0} ; 0, \infty, r ; \varphi\right)=0$.
\end{proposition}

The main difficulty in the proofs of Propositions \ref{prop-Guido-3} and \ref{prop-Guido-4} is to show that $\Phi(t)=\Phi(t ; b, c, r ; \varphi)>0$ for small positive $t$. For this, the condition $w\varphi^{\prime}(w) \rightarrow \infty$ is sufficient.


We next prove the various propositions.

\begin{proof}[Proof of Proposition \ref{prop-Guido-1}] In this case $0<b<c<\infty$. The function $a(x)$ has a single zero in $x \geq 0$ at $x_{0}=\log (q-1)$, and, for $x \geq 0$, 
	\eqn{
	\label{Guido(2-2)}
	   a(x)>0 \text{ for }0 \leq x<x_{0}, \qquad \text{while} \qquad  a(x)<0\text{ for }x>x_{0}.
	}
As a consequence, for $t>0$,
	\eqn{
	\label{Guido(2-3)}
	   \Phi(t)>0\text{ for }t \leq x_{0} / c, \qquad \text{while}
    \qquad \Phi(t)<0\text{ for }t \geq x_{0} / b.
	}
Hence, $\Phi(t)$ has at least one zero in $(x_{0}/c, x_{0} / b)$. Let $t_{0}$ be the smallest such zero, so that $\Phi(t)>0$ for $t \in[x_{0} / c, t_{0})$. 
In order to prove that $t_0$ is the unique zero of $\Phi(t)$, we shall show that  $\Phi'(t) \leq 0$ if $t \in(x_{0}/c, x_{0} / b)$ is such that $\Phi(t) \leq 0$.
We compute
	\eqan{
	\label{Guido(2-4)}
	\Phi^{\prime}(t) & =\frac{d}{d t}\left[\int_{b t}^{c t} x^{r} \e^{-\varphi(x / t)} a(x) d x\right] \\
	& =-\left.b x^{r} \e^{-\varphi(x / t)} a(x)\right|_{x=b t}+\left.c x^{r} \e^{-\varphi(x / t)} a(x)\right|_{x=c t}\nn\\
	& \qquad+\int_{b t}^{c t} x^{r} \e^{-\varphi(x / t)} \cdot \frac{x}{t^{2}} \varphi^{\prime}\left(\frac{x}{t}\right) a(x) d x .\nn
	}
We have $b t<x_{0}< ct$, and so $a(b t)>0>a(c t)$. Therefore, the two expressions on the second line of \eqref{Guido(2-4)} are
both negative, and we see that
	\eqn{
	\label{Guido(2-5)}
	\Phi^{\prime}(t)<\frac{1}{t} \int_{b t}^{c t} x^{r} \e^{-\varphi(x / t)} \frac{x}{t} \varphi^{\prime}\left(\frac{x}{t}\right) a(x) d x.
	}
Since $w \varphi^{\prime}(w)$ is non-decreasing in $w \in[b, c]$, by \eqref{Guido(2-2)},
	\eqn{
	\label{Guido(2-6)}
	\left(\frac{x}{t} \varphi^{\prime}\left(\frac{x}{t}\right)-\frac{x_{0}}{t} \varphi^{\prime}\left(\frac{x_{0}}{t}\right)\right) a(x) \leq 0, 
	\qquad bt \leqslant x \leqslant c t.
	}
Therefore,
	\eqan{
	\label{Guido(2-7)}
	\Phi'(t)\,t
    &<\nn\int_{b t}^{c t} x^{r} \e^{-\varphi(x / t)} \frac{x}{t} \varphi^{\prime}\left(\frac{x}{t}\right) a(x) d x
	=\int_{b t}^{c t} x^{r} \e^{-\varphi(x / t)} \frac{x_{0}}{t} \varphi^{\prime}\left(\frac{x_{0}}{t}\right) a(x) d x\\
	&\qquad+ \int_{b t}^{c t} x^{r} \e^{-\varphi(x / t)}
	\left(\frac{x}{t} \varphi^{\prime}\left(\frac{x}{t}\right)-\frac{x_{0}}{t} \varphi^{\prime}\left(\frac{x_{0}}{t}\right)\right) a(x) d x\nn\\
	&\leq  \frac{x_{0}}{t} \varphi^{\prime}\left(\frac{x_{0}}{t}\right) \int_{b t}^{c t} x^{r} \e^{-\varphi(x/t)} a(x) d x=\frac{x_{0}}{t} \varphi^{\prime}\left(\frac{x_{0}}{t}\right) \Phi(t) \leq 0,
	}
since $\varphi^{\prime}\left(\frac{x_{0}}{t}\right) \geq 0$ and $\Phi(t) \leqslant 0$. Hence, $\Phi^{\prime}(t)<0$.

Now, arguing by contradiction, suppose that there are points $t \in\left(t_{0}, x_{0}/b\right)$ such that $\Phi(t)=0$, and let

	\[
	t_{1}=\inf \left\{t \in (t_{0}, x_{0}/b) \mid \Phi(t)=0\right\}.
	\]
Then $\Phi\left(t_{1}\right)=0$. Furthermore, from $\Phi\left(t_{0}\right)=0$, we have $\Phi^{\prime}\left(t_{0}\right)<0$, and so $\Phi(t)<0$ in a right neighbourhood of $t_{0}$. Therefore, $t_{1}>t_{0}$, and we have $\Phi(t)<0$ for $t_{0}<t<t_{1}$. It follows that $\Phi'(t)<0$ for $t \in\left[t_{0}, t_{1}\right]$. Hence $\Phi(t)$ is strictly decreasing in $t \in[t_{0}, t_{1}]$, contradicting $\Phi(t_{0})=0=\Phi(t_1)$.
\end{proof}

\begin{proof}[Proof of Proposition \ref{prop-Guido-3}]
In this case $b>0$ and $c=\infty$. We start by showing that $\e^{-\varphi(w)}=O\left(w^{-\tau}\right)$, $w \rightarrow \infty$, for all $\tau>2$. Let $\psi(w)=w \varphi^{\prime}(w)$ for $w \geq b$.  Then $\psi(w) \uparrow \infty$, $w\to\infty$ and, for $z \geq b^{2}$,
	\eqan{
	\label{Guido(3-2)}
	\varphi(z)&=\varphi(b)+\int_{b}^{z} \varphi^{\prime}(w) d w=\varphi(b)+\int_{b}^{z} \frac{\psi(w)}{w} d w\nn\\
	& \geqslant \varphi(b)+\int_{z^{1 / 2}}^{z} \frac{\psi(w)}{w} d w \ge\varphi(b)+\psi\left(z^{1 / 2}\right) \int_{z^{1 / 2}}^{z} 
	\frac{d w}{w} \\
	& \geq \varphi(b)+\frac{1}{2} \psi\left(z^{1 / 2}\right) \log z.\nn
	}
Therefore, for any $\tau>2$,
	\eqn{
	\label{Guido(3-3)}
	\exp (-\varphi(z)) \leqslant \exp (-\varphi(b)) z^{-\frac{1}{2} \psi\left(z^{1 / 2}\right)} \leqslant \exp (-\varphi(b)) z^{-\tau},
	}
when $\frac{1}{2} \psi(z^{1/2}) \geq \tau$, which holds when $z$ is large enough.

We next note that $\Phi(t)<0$ when $b t \geq x_{0}$, since $a(x)<0$ when $x>x_{0}$. To see that $\Phi(t)>0$ for small $t>0$, we return to the first expression for $\frac{d^{2}}{d t^{2}} \mathscr{F}(t)$ in \eqref{Guido(4)}, so that
	\eqn{
	\label{Guido(3-4)}
	\frac{d^{2}}{d t^{2}} \mathscr{F}(t)=\frac{q C}{\expec[W]} 
	\int_{b}^{\infty} \frac{(q-1) \e^{t w}-\e^{2 t w}}{\left(\e^{t w}+q-1\right)^{3}}\, w^{p+3} \e^{-\varphi(w)} d w, 
	\qquad t>0.
	}
We have $\exp (-\varphi(w))=O\left(w^{-\tau}\right), w \rightarrow \infty$, for any $\tau>2$. Therefore, 
	\eqn{
	\label{Guido(3-5)}
	\int_{b}^{\infty} w^{p+3} \e^{-\varphi(w)} d w<\infty,
	}
so that $\frac{d^{2}}{d t^{2}} \mathscr{F}(t)$ is continuous in $t \geq 0$, and, in particular,
	\eqn{
	\label{Guido(3-6)}
	\lim _{t \searrow 0} \frac{d^{2}}{d t^{2}} \mathscr{F}(t)=\frac{q C}{\expec[W]} 
	\int_{b}^{\infty} \frac{q-2}{q^{3}} w^{p+3} \e^{-\varphi(w)} d w
	}
is finite and positive. Thus, by continuity, $\frac{d^{2}}{d t^{2}} \mathscr{F}(t)$ is positive for small $t>0$, and so is $\Phi(t)$.

The proof can now be completed in the same way as the proof of Proposition \ref{prop-Guido-1}. Thus, from $\Phi(t)>0$ for small $t>0$ and $\Phi(t)<0$ for $t \geq x_{0} / b$, we conclude that $\Phi(t)$ has at least one positive zero in the range $(0, x_{0}/ b)$. Furthermore, when $t$, for $0<t<x_{0} / b$, is such that $\Phi(t) \leq 0$, we have that $\Phi^{\prime}(t)<0$. For, in the present case, the term on the second line of \eqref{Guido(2-4)} of the proof of Proposition \ref{prop-Guido-1} that involves $b$ is negative while the term involving $c=\infty$ is zero. We omit further details.
\end{proof}

\begin{proof}[Proof of Proposition \ref{prop-Guido-2}]
In this case $b=0$ and $c<\infty$. Since $b=0$, we have to restrict $p$ in the form of the density $f_{\sss W}(w)$ in \eqref{Guido(3)} to the range $p>-1$. Therefore, $r=p+3>2$. We have that $c t \leq x_{0}$ implies that $\Phi(t)>0$ as in the proof of Proposition \ref{prop-Guido-1}. We shall now show that
	\eqn{
	\label{Guido(4-2)}
	\lim _{t \rightarrow \infty} \Phi(t)
	=\lim _{t \rightarrow \infty} \int_{0}^{c t} x^{r} \e^{-\varphi(x/t)} a(x) d x
	=\e^{-\varphi(0)} \int_{0}^{\infty} x^{r} a(x) d x<0 .
	}
The function $u(x)=x^{r},$ for $x \geq 0$, is increasing in $x \geq 0$. Hence, 
	\eqan{
	\label{Guido(4-3)}
	\int_{0}^{\infty} u(x) a(x) d x&=u\left(x_{0}\right) \int_{0}^{\infty} a(x) d x
	+\int_{0}^{\infty}\left(u(x)-u\left(x_{0}\right)\right) a(x) d x\\
	&\leqslant u\left(x_{0}\right) \int_{0}^{\infty} a(x) d x,\nn
	}
where we have used that $\left(u(x)-u\left(x_{0}\right)\right) a(x) \leq 0$ for $0 \leq x<\infty$. Now, from
	\eqn{
	\label{Guido(4-4)}
	a(x)=\frac{d}{d x}\left[\frac{\e^{x}}{\left(\e^{x}+q-1\right)^{2}}\right],\qquad x \geq 0,
	}
we compute
	\eqn{
	\label{Guido(4-5)}
	\int_{0}^{\infty} a(x) d x=\left.\frac{\e^{x}}{\left(\e^{x}+q-1\right)^{2}}\right|_{0} ^{\infty}=-\frac{1}{q^{2}}<0.
	}
Hence, \eqref{Guido(4-2)} holds, and so $\Phi(t)<0$ for large $t>0$.

The proof can now be completed in the same way as that of Proposition \ref{prop-Guido-1}. 
Thus, from $\Phi(t)>0$ for $t \leq x_{0}/c$ and $\Phi(t)<0$ for large $t$, we conclude that $\Phi(t)$ has at least one positive zero
in the range $(x_{0}/c, \infty)$. Furthermore, when $t>x_{0}/c$ is such that $\Phi(t) \leqslant 0$, we have that $\Phi^{\prime}(t)<0$. For, in the present case, the term on the second line in \eqref{Guido(2-4)} 
 of the proof of Proposition \ref{prop-Guido-1} that involves $b=0$ is zero while the term involving $c<\infty$ is negative. We omit further details.
\end{proof}

\begin{proof}[Proof of Proposition \ref{prop-Guido-4}]
In this case, $b=0$ and $c=\infty$. As in the proof of Proposition \ref{prop-Guido-2}, we have to restrict to $p>-1$, i.e., to $r=p+3>2$. Also, as in the proof of Proposition \ref{prop-Guido-3}, we can show that $\Phi(t)>0$ for small $t>0$, and, as in the proof of Proposition \ref{prop-Guido-2}, we can show that $\Phi(t)<0$ for large $t>0$. Thus, $\Phi(t)$ has at least one positive zero in $(0, \infty)$. Furthermore, $\Phi^{\prime}(t)<0$ when $t>0$ is such that $\Phi(t) \leq 0$. Indeed, in the present case, both terms on the second line of \eqref{Guido(2-4)} of the proof of Proposition \ref{prop-Guido-1} vanish, while the term on the last line is negative. For the latter, instead of \eqref{Guido(2-6)} in the proof of Proposition \ref{prop-Guido-1}, we now have
	\eqn{
	\label{Guido(4-7)}
	0 \not\equiv\left(\frac{x}{t} \varphi^{\prime}\left(\frac{x}{t}\right)
	-\frac{x_{0}}{t} \varphi^{\prime}\left(\frac{x_{0}}{t}\right)\right) a(x) \leqslant 0,
	\qquad 0 \leqslant x<\infty,
	}
since $w \varphi^{\prime}(w) \rightarrow \infty$ as $w=x/t \rightarrow \infty$. We omit further details.
\end{proof}

\subsection{The Pareto case: Proof of Theorem \ref{thm-Pareto}} 
\label{sec-Pareto-(3,4)}
We next take $b=1, c=\infty, r \in \mathbb{R}, \varphi \equiv 0,$ so that we reduce to the Pareto case in \eqref{Pareto}. In this case, \eqn{
	\label{Guido(7)}
	\Phi(t)=\Phi(t ; b=1, c=\infty, r ; \varphi=0)=\int_{t}^{\infty} x^{r} a(x) d x, \qquad t>0.
	}
We require in view of \eqref{main-assumption-Fn} and \eqref{Guido(3)} that $p<-2$, and so $r=3-\tau=p+3<1$. The proof of Theorem \ref{thm-Pareto} follows from the following results:

\begin{proposition}
\label{prop-Guido-5} 
In order that $\Phi(t)$ has exactly one zero $t_{0}\in(0,\infty)$, it is necessary and sufficient that $\lim _{t \rightarrow 0} \Phi(t)>0$ (with $\infty$ allowed).
\end{proposition}

\begin{proposition}
\label{prop-Guido-6}
For $r \geq 0$,
	\eqn{
	\label{Guido(8)}
	\lim _{t \searrow 0} \Phi(t)=\int_{0}^{\infty} x^{r} a(x) d x<0.
	}
For r $\leqslant-1$,
	\eqn{
	\label{Guido(9)}
	\lim _{t \searrow 0} \Phi(t)=\int_{0}^{\infty} x^{r} a(x) d x = \infty.
	}
Further,
	\eqn{
	\label{Guido(10)}
	\lim _{t \searrow 0} \Phi(t)=\int_{0}^{\infty} x^{r} a(x) d x \rightarrow \infty, \qquad r \downarrow-1.
	}
\end{proposition}

We prove, furthermore, the following result:

\begin{proposition}
\label{prop-Guido-7} 
\invisible{Let $q=3,4, \ldots$.} There is a unique $r=r(q) \in(-1,0)$ such that
\eqn{
	\label{Guido(5-9)}
	r\in (-1, r(q)) \Rightarrow \lim_{t\searrow 0} \Phi(t)>0;\qquad  r\in (r(q),0) \Rightarrow  \lim_{t\searrow 0} \Phi(t)<0.
	}
Consequently, $r=r(q)$ satisfies
    \eqn{
    \label{r(q)-equality}
    \lim _{t \searrow 0} \Phi(t)=\int_{0}^{\infty} x^{r} a(x) dx=0.
    }
    
\end{proposition}
We note that $r=3-\tau,$ so $r(q)\in (-1,0)$ implies that $\tau(q)\in(3,4)$.

\begin{proposition}[Asymptotics of $q\rightarrow r(q)$]
\label{prop-r(q)-asymp}
As $q\rightarrow \infty$, 
    \eqn{
    r=r(q)\sim-\frac{\log{q}}{q}.
    }
More precisely, with $x_0=\log{(q-1)}$, 
    \eqn{
    \label{r(q)-asymptoics-sharp}
    r(q)=\frac{-x_{0}}{q-1-(q-1)^{-1}-2 x_{0}}\left(1+O\left(\frac{1}{x_{0}^{2}}\right)\right) \sim \frac{-\log q}{q}.
    }
\end{proposition}

We next prove these propositions:
\begin{proof}[Proof of Propositions \ref{prop-Guido-5}--\ref{prop-Guido-6}]
In the Pareto case that we consider $(b=1, c=\infty, r \in \mathbb{R}, \varphi \equiv 0)$,
	\eqn{
	\label{Guido(5-1)}
	\Phi(t)=\int_{t}^{\infty} x^{r} a(x) d x, \qquad t>0.
	}
We require $p<-3$, i.e., $r<0$.

Recall \eqref{Guido(5)}, and $a\left(x_{0}\right)=0$ with $x_{0}=\log (q-1)$, and
	\eqn{
	\label{Guido(5-2)}
	   a(x)>0\text{ for }0 \leq x<x_{0} \qquad
       \text{while} \qquad a(x)<0\text{ for }x>x_{0}.
	}
As a consequence, $\Phi(t)$ decreases strictly in $t \in(0, x_{0})$ and increases strictly in $t \in\left(x_{0}, \infty\right)$. Furthermore, $\Phi(t)<0$ for $t \geq x_{0}$. 
Thus, $\Phi(t)$ has at most one zero in $(0,\infty)$, and this zero lies in $(0,x_0)$ when it exists. Moreover,
\be
\Phi \text{ has a zero in } (0,x_0) \iff \lim_{t\searrow 0} \Phi(t) >0 \text{ (with $\infty$ allowed). }
\ee
This proves Proposition  \ref{prop-Guido-5}. Next, for $r\ge 0$,
\be
\lim_{t\searrow 0 } \Phi(t) = \int_0^\infty x^r a(x)dx <0,
\ee
where the inequality was shown in the course of the proof of  Proposition  \ref{prop-Guido-2}.
Furthermore, from $a(0) = \frac{q-2}{q^3} >0$ and smoothness of $a(x)$ at $x=0$ 
\be
r\le -1 \Longrightarrow \lim_{t\searrow 0 } \Phi(t) = \infty.
\ee
Moreover, for $r>-1$,
	\eqn{
	\label{Guido(5-7)}
	L_{q}(r):=\lim _{t \searrow 0} \Phi(t)=\int_{0}^{\infty} x^{r} a(x) d x, 
	}
and
	\eqn{
	\label{Guido(5-8)}
	L_{q}(r)=\frac{a(0)}{r+1}+O(1) \rightarrow \infty, \qquad r \downarrow-1.
	}
This proves Proposition  \ref{prop-Guido-6}.	
\end{proof}

\begin{proof}[Proof of Proposition \ref{prop-Guido-7}]
The function $L_{q}(r)$, defined in  \eqref{Guido(5-7)}, is continuous in $r \in(-1,0]$, and we have $L_{q}(r) \rightarrow \infty$ when $r \downarrow-1$, and $L_{q}(r) \rightarrow L_{q}(0)<0$ when $r\uparrow 0$. Hence, there is at least one $r$ $\in(-1,0)$ such that $L_{q}(r)=0$. Let $r \in(-1,0)$ be such that $L_{q}(r)=0$, and let $\varepsilon \in \mathbb{R}, \varepsilon \neq 0$, be such that $r+\varepsilon \in(-1,0)$. Then
	\eqan{
	\label{Guido(5-10)}
	L_{q}(r+\varepsilon) & =\int_{0}^{\infty} x^{r+\varepsilon} a(x) d x=\int_{0}^{\infty} x^{r} x^{\varepsilon} a(x) d x \\
	& =\int_{0}^{\infty} x^{r}\left(x^{\varepsilon}-x_{0}^{\varepsilon}\right) a(x) d x+x_{0}^{\varepsilon} \int_{0}^{\infty} x^{r} a(x) d x \nn\\
	& =\int_{0}^{\infty} x^{r}\left(x^{\varepsilon}-x_{0}^{\varepsilon}\right) a(x) d x,\nn
	}
since $L_{q}(r)=0$ by assumption. The function $x^{\varepsilon}-x_{0}^{\varepsilon}$, for $x>0$, vanishes at $x=x_{0}$, and is strictly increasing in $x>0$ when $\varepsilon>0$ and strictly decreasing in $x>0$ when $\varepsilon<0$. Therefore, when $\varepsilon>0$,
	\eqn{
	\label{Guido(5-11)}
	\left(x^{\varepsilon}-x_{0}^{\varepsilon}\right) a(x) < 0, \quad x>0, x\neq x_0,
	}
and, when $\varepsilon<0$,
	\eqn{
	\label{Guido(5-12)}
	\left(x^{\varepsilon}-x_{0}^{\varepsilon}\right) a(x) > 0, \quad x>0, x\neq x_0.
	}
Therefore, $L_{q}(r+\varepsilon)<0$ when $\varepsilon>0$ and $L_{q}(r+\varepsilon)>0$ when $\varepsilon<0$. Finally,  when $r=r(q)$, we have $\lim_{t\to 0} \Phi(t)=0$, and so, by Proposition \ref{prop-Guido-5}, $\Phi(t)$ has no positive zero.
\end{proof}

\begin{proof}[Proof of Proposition \ref{prop-r(q)-asymp}]
This is given in Appendix \ref{appB}.
\end{proof}


\subsection{Unique zero-crossing condition holds for weights with small support: proof of Theorem \ref{thm-two-atoms}}
\label{sec-W-bounded-support}
The goal of this section is to show that the \emph{unique zero-crossing} condition
holds for weight variables $W$ with a distribution that has support contained in a compact set that is not too large.

In Theorem \ref{thm-two-atoms} we assume that the support of $W$ is contained in a compact set of the form $[w',w]$. By scaling, we may assume (without loss of generality) that $w'=1$, that is, $\prob(W\geq 1)=1$, so that $\prob(W\in[1,\olW])=1$. Then
it is enough to prove the following:

\invisible{We note that assumption \eqref{eq:bounded_supp_W_assump} is invariant under scaling: if $W$ satisfies it, then so does $CW$ for any constant $C >0$:}
\invisible{
\begin{proposition}[Weights on small compacts]
\label{prop:unique_inflection_for_bdd_supp}
    Fix $0\leq B<\log(q-1)$. Let $W$ be a random variable satisfying $\prob(W\in [1,w])=1.$ Assume \begin{equation}\label{eq:bounded_supp_W_assump}
    \olW\leq 1+\frac{\log{(2+\sqrt{3})}}{\log{(q-1)}-B}.
    \end{equation} 
    Then, for such $W$, $t \mapsto\mathscr{F}_B(t)$ satisfies Condition \eqref{cond-zer-cross}(a).
\end{proposition}

{Recall from \eqref{eq:second_derivative_gen} that
\begin{align}
\label{H(t)-def}
   H(t):= \frac{d^2}{dt^2}\mathscr{F}_B(t)
    &=q \expec\left[\frac{W^3}{\expec[W]}\frac{\e^{tW+B}(q-1-\e^{tW+B})}{(\e^{tW+B}+q-1)^3}\right].
\end{align}
Before proving Proposition \ref{prop:unique_inflection_for_bdd_supp}, we first do some analysis when $W$ is a constant, say $W=C$. In this case $H(t)$ takes the form
\begin{align*}
   H(t)= q\,C^2\,\frac{\e^{tC+B}(q-1-\e^{tC+B})}{(\e^{tC+B}+q-1)^3}.  
\end{align*}
Writing $p=q-1$ and, for any $a>0$,
\begin{align*}
    t_B(a):=\frac{\log p-B}{a},
\end{align*}
we note that $H(t)$ in this case has a unique root at $t_B(C)$. 
\invisible{Note that $t_B(C)>0$ for any $C>0$ since $B<\log(q-1)$.}
We also observe that $t_B(C)>0$ for any $C>0$ since $B<\log(q-1)$, and that 
$H(t)>0$ for $t\in (0,t_B(C))$ and 
$H(t)<0$ for $t\in (t_B(C),\infty)$.
Further, we can easily check that
\begin{align*}
    H'(t)&=\frac{q\,C^3\,\e^{tC+B}}{(\e^{tC+B}+p)^4}[(\e^{tC+B})^2+p^2-4pe^{tC+B}]\\&=\frac{q\,C^3\,\e^{tC+B}}{(\e^{tC+B}+p)^4}[(\e^{tC+B}-p(2+\sqrt{3}))(\e^{tC+B}-p(2-\sqrt{3}))],
\end{align*}
so that $H$ has two critical points $t^{-}_B(C):=\frac{\log(p(2-\sqrt{3}))-B}{C}$ and $t^{+}_B(C):=\frac{\log(p(2+\sqrt{3}))-B}{C}$, and observe that $t^{-}_B(C)<t_B(C)<t^{+}_B(C)$, i.e., $H$ is strictly decreasing at its unique root $t_B(C)$. With these observations in hand, we are ready to give the proof of Proposition \ref{prop:unique_inflection_for_bdd_supp}.

\begin{proof}[Proof of Proposition \ref{prop:unique_inflection_for_bdd_supp}]
   \invisible{As in the last discussion, we write
   \begin{align*}
   H(t)=\frac{d^2}{dt^2}\mathscr{F}_B(t)= \frac{q}{\expec[W]}\expec\left[W^3\frac{p\e^{tW+B}-\e^{2(tW+B)}}{(\e^{tW+B}+p)^3}\right].   
   \end{align*}
   Further, recall for $a>0$, we use the notation $$t_B(a):=\frac{1}{a}(\log p-B),$$ where $p=q-1$.}
   Recall $H(t)$ from \eqref{H(t)-def}. 
   We note that $H(t)>0$ for $t \in (0,t_B(\olW))$, while $H(t)<0$ for $t \in (t_B(\ulW),\infty)$. This is because, for any fixed value of $W \in [1,w]$, this holds for $\frac{p\e^{tW+B}-\e^{2(tW+B)}}{(\e^{tW+B}+p)^3}$, as can be verified from the analysis done for $W=C$ above. Hence, by the intermediate value theorem, $H(t)$ has at least one root, and the set of all its roots is contained in the interval $(t_B(\olW),t_B(\ulW))$. We want to prove that under the assumption \eqref{eq:bounded_supp_W_assump}, the function $H(t)$ is strictly decreasing on the interval $(t_B(\olW),t_B(\ulW))$, which will then imply that it has a unique root.

We have as above
 \invisible{\begin{align*}
       H'(t)=\frac{d^3}{dt^3}\mathscr{F}_B(t)&=\frac{q}{\expec[W]}\expec\left[W^4 \frac{\e^{tW+B}(p^2-4p\e^{tW+B}+\e^{2(tW+B)})}{(\e^{tW+B}+p)^4}\right]\\&=\frac{q}{\expec[W]}\expec\left[W^4 \frac{\e^{tW+B}(\e^{tW+B}-(2-\sqrt{3})p)(\e^{tW+B}-(2+\sqrt{3})p)}{(\e^{tW+B}+p)^4}\right]\\&=\frac{q}{\expec[W]}\expec\left[W^4 \e^{tW+3B}\frac{(\e^{tW}-(2-\sqrt{3})p')(\e^{tW}-(2+\sqrt{3})p')}{(\e^{tW+B}+p)^4}\right],
   \end{align*}}
 \begin{align*}
       H'(t)=\frac{d^3}{dt^3}\mathscr{F}_B(t)&=\frac{q}{\expec[W]}\expec\left[W^4 \frac{\e^{tW+B}(p^2-4p\e^{tW+B}+\e^{2(tW+B)})}{(\e^{tW+B}+p)^4}\right]\\&=\frac{q}{\expec[W]}\expec\left[W^4 \frac{\e^{tW+B}(\e^{tW+B}-(2-\sqrt{3})p)(\e^{tW+B}-(2+\sqrt{3})p)}{(\e^{tW+B}+p)^4}\right]\\&=\frac{q}{\expec[W]}\expec\left[W^4 \e^{tW+3B}\frac{h(t)}{(\e^{tW+B}+p)^4}\right],
   \end{align*}
where we use the notation 
$$
h(t):= (\e^{tW}-(2-\sqrt{3})p')(\e^{tW}-(2+\sqrt{3})p'),
$$
and $p'=\e^{-B}p=\e^{-B}(q-1)>1$. Recall that the values of $W$ that contribute to the above expectation all lie in the interval $(\ulW,\olW)$, and also that we are only interested in $t \in (t_B(\olW),t_B(\ulW))
\equiv(\frac{\log p-B}{\olW},\log p-B)=(\frac{\log p'}{\olW},\log p')$.
We shall show that $h(t)<0$
for $\ulW \le W \le \olW$ and $\frac{\log p'}{\olW} \le t \le \log p'$.
\medskip

We distinguish two cases according to the value of $p'$:\\
$\bullet$ $1<p'\leq 2+\sqrt{3}$. In this case, $\e^{tW}-(2-\sqrt{3})p'>0$ for any $t>0$ and $W>0$ since $(2-\sqrt{3})p'\leq 1$. 
On the other hand, $\e^{tW}< \e^{w \log p'}\le (2+\sqrt{3})p'$ since $t<t_B(1)=\log p'$ and $W\leq w$, where in the last step we have used \eqref{eq:bounded_supp_W_assump}. Hence, $h(t)<0$ in $t \in (t_B(\olW),t_B(\ulW))$.
\medskip

\noindent
$\bullet$ $p'> 2+\sqrt{3}$. Let us first show the inequality  $\e^{tW}-(2-\sqrt{3})p'>0$. Observe that this will be implied by the inequality $1/w \geq 1+\frac{\log(2-\sqrt{3})}{\log p'}$, since $t>\frac{\log p'}{w}$ and $W\geq 1$. Since $(2-\sqrt{3})(2+\sqrt{3})=1$, we have
$$1+\frac{\log(2-\sqrt{3})}{\log p'}=1-\frac{\log(2+\sqrt{3})}{\log p'}.$$ Furthermore, by \eqref{eq:bounded_supp_W_assump} and the inequality $1\geq (1+x)(1-x)=1-x^2$ with $x=\frac{\log (2+\sqrt{3})}{\log p'}$, 
$$1/w\geq \frac{1}{1+\frac{\log (2+\sqrt{3})}{\log p'}}\geq 1-\frac{\log(2+\sqrt{3})}{\log p'}=1+\frac{\log(2-\sqrt{3})}{\log p'}.$$ The inequality $\e^{tW}-(2+\sqrt{3})p'<0$ is obtained as before by applying again \eqref{eq:bounded_supp_W_assump}: 
$\e^{tW} \le \e^{w \log p'} \le (2+\sqrt{3})p$. We conclude that 
$h(t)<0$ for all $t \in (t_B(\olW),t_B(\ulW))$ and, thus, that $H(t)$, being decreasing on the same interval, has a unique root. This proves that Condition \eqref{cond-zer-cross}(a) is satisfied. 
\end{proof}

}
}

\begin{proposition}[Weights on small compacts]
\label{prop:unique_inflection_for_bdd_supp}
    Fix $0\leq B<\log(q-1)$, and let $W$ be a random variable satisfying $\prob(W\in [1,w])=1.$ Assume \begin{equation}\label{eq:bounded_supp_W_assump}
    \olW\leq 1+\frac{\log{(2+\sqrt{3})}}{\log{(q-1)}-B}.
    \end{equation} 
    Then, for such $W$, $t \mapsto\mathscr{F}_B(t)$ satisfies Condition \ref{cond-zer-cross}(a).
\end{proposition}
\begin{proof}[Proof of Proposition \ref{prop:unique_inflection_for_bdd_supp}]
Recalling the convention in the beginning of Section \ref{sec-unique-zero-crossing} of replacing $\e^{-B}(q-1)$ by $q-1>1$, we can consider $\mathscr{F}_0$ instead of $\mathscr{F}_B$, and the condition in \eqref{eq:bounded_supp_W_assump} becomes
\begin{equation}\label{eq:bounded_supp_W_assumpp}
    \olW\leq 1+\frac{\log{(2+\sqrt{3})}}{\log p}\, ,\quad p:=q-1 .
    \end{equation} 
    Furthermore, 
    \begin{align}
\label{H(t)-def}
   H(t):= \frac{d^2}{dt^2}\mathscr{F}_0(t)
    &=\frac{q}{\expec[W]} \expec\left[{W^3}\,\frac{p\,\e^{tW}-\e^{2tW}}{(\e^{tW}+p)^3}\right].
\end{align}
In the case that $\prob(W=C)=1$ for some $C>0$, we have 

\begin{align*}
   H(t)= q\,C^2\,\frac{\e^{tC}(p-\e^{tC})}{(\e^{tC}+p)^3}.  
\end{align*}
Then $H(t)>0$ for $t<\frac{1}{C} \log p$ and $H(t)<0$ for $t>\frac{1}{C} \log p$, and so the result holds in this case.

For the general case, we assume, in addition to $\prob(W \in [1,w])=1$, that $\prob{(W=1)} \ne 1 \ne \prob(W=w)$. With these assumptions, we have for $t>0$
\begin{equation}
t \le \frac{1}{w} \log p\, \Rightarrow H(t)>0,\quad t\ge \log p\, \Rightarrow  H(t) <0.
\end{equation}
In particular, $H(\frac{1}{w} \log p) >0 > H(\log p)$. Hence, by the intermediate value theorem, $H(t)$ has at least one zero in $(\frac 1 w  \log p,\, \log p )$. We shall show that $H(t)$ is strictly decreasing in $t \in (\frac 1 w  \log p,\, \log p )$. We have from \eqref{H(t)-def}
 \begin{align*}
       H'(t)&=\frac{q}{\expec[W]}\expec\left[W^4\, \e^{tW}\,\frac{p^2-4\,p\,\e^{tW}+\e^{2tW}}{(\e^{tW+B}+p)^4}\right]\\&=\frac{q}{\expec[W]}\expec\left[W^4\, \e^{tW} \frac{(\e^{tW}-(2-\sqrt{3})p)(\e^{tW}-(2+\sqrt{3})p)}{(\e^{tW}+p)^4}\right]\,.
   \end{align*}
   We shall show that
   \begin{equation}
       h(t,W):= (\e^{tW}-(2-\sqrt{3})p)(\e^{tW}-(2+\sqrt{3})p)<0
   \end{equation}
   for $t \in (\frac 1 w  \log p,\, \log p )$ and $W\in [1,w]$.
   We first show that 
   \begin{equation}
       \label{eq:lowbounb1w}
       \frac{1}{w} \ge 1+ \frac{\log (2-\sqrt{3})}{\log p}\, .
   \end{equation}
   Indeed, by \eqref{eq:bounded_supp_W_assumpp} and the inequality $(1+x)(1-x)=1-x^2\le 1$ with $x=\frac{\log(2+\sqrt{3})}{\log p}>0$, we have
   \begin{equation}
       \frac 1 w \ge \frac{1}{1+\log (2+\sqrt{3})/\log p} \ge 1- \frac{\log (2+\sqrt{3})}{\log p}= 1+ \frac{\log (2-\sqrt{3})}{\log p}.
   \end{equation}
   From \eqref{eq:lowbounb1w}, we get when $t> \frac 1 w \log p$ and $W\ge 1$
   \begin{equation}
       t\, W > \frac{\log p}{w} \ge \log p + \log (2-\sqrt{3}) = \log (2-\sqrt{3})p\,,
   \end{equation}
   i.e. $\e^{tW} > (2-\sqrt{3})\,p$. Furthermore, from $t< \log p\,, W\le w$ and \eqref{eq:bounded_supp_W_assumpp}, we get
   \begin{equation*}
       t\,W < w \log p \le \log p + \log (2+\sqrt{3})= \log({2+\sqrt{3})\,p},
   \end{equation*}
   i.e., $\e^{tW}<(2+\sqrt{3})\, p$. Therefore, we have $h(t,W)<0$ when $t \in (\frac 1 w  \log p,\, \log p )$ and $W \in [1,w]$.
\end{proof}

\invisible{
\begin{remark}[Smallest positive value]
    {\rm We note that $\ulW$ can very well be $0$, in which case the condition \eqref{eq:bounded_supp_W_assump} is never satisfied. In fact, we can replace $\ulW$ in the assumption \eqref{eq:bounded_supp_W_assump} by the smallest strictly positive value it can take with positive probability (let us call this value $W_{>0}$). This is because we note that if we expand out the expectation in the expression of $H(t)$, the term corresponding to $W=0$ does not contribute, and the values of $W$ that contribute are between $W_{>0}$ and $\olW$, so that if $\frac{\olW}{W_{>0}}<\frac{\log{(2+\sqrt{3})p}}{\log{p}}$, Proposition \ref{prop:unique_inflection_for_bdd_supp} stays true with the same proof.}
    \hfill \ensymboldefinition
\end{remark}

\begin{remark}[On shifts]
    {\rm We note that if a random variable $W$ satisfies the assumption \eqref{eq:bounded_supp_W_assump}, then so does any arbitrary positive shift of it: we write the condition we need to check for the random variable shifted by $\delta$ for any arbitrary $\delta>0$ as,
    \begin{align*}
        \frac{\olW+\delta}{\ulW+\delta}<p_*,
    \end{align*}
    where $p_*=\frac{\log{(2+\sqrt{3})p}}{\log{p}}$. The last inequality is equivalent to, 
    \begin{align*}
        \olW-p_*<\delta(p_*-1),
    \end{align*}
    which is always true because the r.h.s.\ is positive, while the l.h.s.\ is negative by virtue of the fact that $W$ satisfies \eqref{eq:bounded_supp_W_assump}. The above argument also shows that we can allow for negative shifts $\delta<0$ as long as \begin{align*}
        \frac{\olW-p_*}{p_*-1}<\delta.
    \end{align*}
    }\hfill \ensymboldefinition
\end{remark}}

\section{General first-order phase transition: Proof of Theorem \ref{thm-general-first-order}}
\label{sec-FO-general}
In this section, we investigate the order of the phase transition when we do not assume that $t\mapsto \frac{d^2}{dt^2}\mathscr{F}_B(t)$ is first positive and then negative. This section is organised as follows. In Section \ref{sec-touching-points}, we introduce the set of touching points, and state a key lemma on their properties. For $B=0$, in Section \ref{sec-FO-PT-general}, we give the proof of the first-order phase transition in Theorem \ref{thm-general-first-order} and in Section \ref{sec-prop-touching-points} we prove the key result on touching points in Lemma \ref{lem-proof-touching-point}. The extension to $B>0$ is outlined in Section \ref{sxec-ext-B>0-general-first-order-phase transition}.

\subsection{Touching points}
\label{sec-touching-points}
In this section, we fix $B=0$ and abbreviate $\mathscr{F}(t)=\mathscr{F}_0(t)$. We start by defining the set of {\em touching points}:
\begin{definition}[Touching points]
\label{def-touching-points}
{\rm We say that $\beta'$ is a {\em touching point} when there exists  $t$ such that
    \eqn{
    \label{def-touching-point-t}
    \mathscr{F}(t)=\frac{t}{\beta'},
    \qquad\qquad
    \frac{d}{d t}
    \mathscr{F}(t)=\frac{1}{\beta'},
    }
and, moreover, for at least one $t=\theta(\beta')$ satisfying \eqref{def-touching-point-t}, the sign of $\mathscr{F}(t)-t/\beta'$ does not change in $[\theta(\beta')-\delta, \theta(\beta')+\delta]$ for some $\delta>0$ sufficiently small. The existence of such $\theta(\beta')$ follows when $\min\{k\geq 2\colon \mathscr{F}^{(k)}(\theta(\beta'))\neq 0\}$ is even, where $\mathscr{F}^{(k)}$ is the $k$th derivative of $\mathscr{F}$. 
}\hfill \ensymboldefinition
\end{definition}

\invisible{\Gib{Claudio: the symbol $t(\beta')$ used in this section to denote the solutions of \eqref{def-touching-points} is a little confusing because it is similar to that used in Sec.6 to denote the solutions of the first equation in \eqref{def-touching-points}. See also the following Remark.}
\RvdH{I agree, that is not good. Do we have an alternative? We have used so many letters already...Claudio: I introduced the alternative notation $\theta(\beta')$. }
}
\begin{remark}[Touching points under the unique zero-crossing condition]
{\rm When $t\mapsto \frac{d^2}{dt^2}\mathscr{F}(t)$ is first positive and then negative, as in Condition \ref{cond-zer-cross}(a), there are exactly {\em two} touching points. For any $\beta'$ in between them, there are 3 solutions $t_1(\beta'), t_2(\beta'), t_3(\beta')$ to $\mathscr{F}(t)=t/\beta'$. }\hfill \ensymboldefinition
\end{remark}

We now explore what happens if the unique zero-crossing condition is no longer satisfied. We denote $\TP$ the set of touching points introduced in Definition \ref{def-touching-points}, and also define 
$$
{\TTP}=\{t\ge 0\, :\, t=\theta(\beta') \mbox{ solves \eqref{def-touching-point-t}} \mbox{ for } \beta'\in \TP\},
$$
which is the set of the points $t$ occurring in the definition of touching points. The next lemma investigates relevant properties of the set of touching points:

\begin{lemma}[Touching points]
\label{lem-proof-touching-point}
Assume that $\expec[W^3]<\infty$. There are finitely many touching points. Further, for any $\beta'$ in between two touching points, there are $2k+1$ solutions of $\mathscr{F}(t)=t/\beta'$ for some $k\geq 1$, for any $\beta'<q\expec[W]/\expec[W^2]$.
\end{lemma}

\subsection{Proof of the first-order phase transition in Theorem \ref{thm-general-first-order}}
\label{sec-FO-PT-general}

In this section, we prove Theorem \ref{thm-general-first-order} subject to Lemma \ref{lem-proof-touching-point}:
\medskip

\noindent
{\it Proof of Theorem \ref{thm-general-first-order}.}
Denote $\beta_0'=1/{\mathscr{F}'(t_0)}$ and
$\beta_1'=q\expec[W]/\expec[W^2]$ the smallest and the largest touching point (their existence is proved in Lemma \ref{stetouch}). 
By Lemma \ref{lem-proof-touching-point}, for any $\beta'$ in between two touching points, there are $2k+1$ solutions for some $k\geq 1$. Hereafter, we denote the solutions, arranged in increasing order according to $\ell=1,\ldots, 2k+1$, by $t_\ell(\beta')$. When $\beta'>\beta_1'$, there are exactly two solutions of $\mathscr{F}(t)=t/\beta'$, and $t_2(\beta')$ optimizes $t\mapsto p(t,\beta')\equiv p_{\beta,0}(t/\beta')$ (see \eqref{var-problem-psi-B>0}). When $\beta'<\beta_0'$, there is exactly 1 solution of $\mathscr{F}(t)=t/\beta'$, which is $t_1(\beta')=0$, and $t_1(\beta')$ optimizes $t\mapsto p(t,\beta')$. We are left with dealing with $\beta'\in (\beta_0',\beta_1')$ that are not touching points.
\smallskip

Recall that Theorem \ref{thm-general-first-order} assumes that $\expec[W^3]<\infty$. This implies that $t\mapsto \mathscr{F}(t)$ is initially convex, since $t\mapsto \mathscr{F}''(t)$ is continuous with $\mathscr{F}''(0)>0$ (recall \eqref{eq:second_derivative_gen}).
Since $t\mapsto \mathscr{F}(t)$ is initially convex, 
$\mathscr{F}(t)<t/\beta'$ for all $t\in (t_1(\beta'), t_2(\beta')),$ so that $t_2(\beta')$ cannot be an optimizer of $t\mapsto p(t,\beta').$ When $\beta'$ is not a touching point, we have $\mathscr{F}(t)<t/\beta'$ for all $t\in (t_{2\ell-1}(\beta'), t_{2\ell}(\beta'))$, while
$\mathscr{F}(t)>t/\beta'$ for all $t\in (t_{2\ell}(\beta'), t_{2\ell+1}(\beta'))$. 

Thus, the only possible stationary points that can maximize $t\mapsto p(t,\beta')$ are $t_{2\ell+1}(\beta'),$ for all $\ell\geq 0.$ Indeed, following the previous discussion, we have for all $\ell \ge 0$:
$$
p(t_{2 \ell+1}(\beta'),\beta')-p(t_{2 \ell}(\beta'),\beta') = 
\frac{q-1}{q} \expec[W] \int_{t_{2 \ell(\beta')}}^{t_{2 \ell+1}(\beta')} (\mathscr{F}(t)-\frac{t}{\beta'} ) dt >0. 
$$

Further, let
    \eqn{
    \label{t-conv-def}
    t_{\rm conv}=\sup\{s \colon \mathscr{F}(u) \text{ is convex on }(0,s)\}>0,
    }
since $t\mapsto \mathscr{F}(t)$ is initially convex.
Then, we claim that $t_3(\beta')\geq t_{\rm conv}$ for every $\beta'$. 
\invisible{Indeed, if $\beta'\in (\beta_0',\beta_1')$, then $\mathscr{F}'(0)>1/\beta'$. Since $\mathscr{F}(t_2(\beta'))=t_2(\beta'))/\beta',$ the function $t\mapsto \mathscr{F}(t)-t/\beta'$ changes sign around $t_2(\beta')$. If $t_3(\beta')<t_{\rm conv},$ then $t\mapsto \mathscr{F}(t)-t/\beta'$ would be convex on $[0,t_3(\beta')],$ but this is inconsistent with $t\mapsto \mathscr{F}(t)-t/\beta'$ changing sign at least 3 times on $(0,t_3(\beta')).$} 
Indeed, if $t_3(\beta')<t_{\rm conv}$, then, by convexity, the graph of $\mathscr{F}(t)$ on $(0,t_3(\beta'))$ would be strictly below that of $t/\beta'$, but this is inconsistent with the definition of $t_2(\beta')$.
\smallskip

Thus, $t_3(\beta')\geq t_{\rm conv}$ for every $\beta'$ that is not a touching point. As a result, also $t_{2k+1}(\beta')\geq t_{\rm conv}$ for every $\beta'$ that is not a touching point. 
We conclude that if $t_{2k+1}(\beta')$ optimizes $t\mapsto p(t,\beta')$ for some $k\geq 1$, and $\beta'$ is not a touching point, then $t_{2k+1}(\beta')>t_{\rm conv}$. 

Let $t_{\rm opt}(\beta')$ be the optimizer of $t\mapsto p(t,\beta')=p_{\beta,0}(t/\beta')$. We conclude that $\beta'\mapsto t_{\rm opt}(\beta')$ has a first-order phase transition, where it jumps from $t_{\rm opt}(\beta')=0$ to $t_{\rm opt}(\beta')>t_{\rm conv}>0$. Since the optimal value $s^\star(\beta,0)=t_{\rm opt}(\beta')/\beta'$, this also proves that $\beta\mapsto s^\star(\beta,0)$ has a first-order phase transition.
\qed

\begin{remark}[Uniqueness of the jump of $\beta\mapsto t_{\rm opt}(\beta')$]
{\rm The above argument shows that the phase transition in our model is in general of first order, as expected. 
What we cannot show, however, is that the jump is {\em unique}. For that to occur, we further need that the jump is to the {\em largest} solution $t_{2k+1}(\beta')$ after the jump. This largest solution $t_{2k+1}(\beta')$ can be expected to depend continuously on $\beta'$, which would show that there is a {\em unique} jump.
}\hfill \ensymboldefinition
\end{remark}
\invisible{
\Gib{Claudio: the Implicit Function Theorem should imply that, except that at the touching points, the solutions are continuous and even more regular.  }
}

\invisible{Note that $t=0$ is {\em not} a solution to this equation. However, when $B\searrow 0,$ this smallest solution {\em does} converge to 0. 
\RvdH{Is this obvious? Somewhere, we should still use that $\beta'<\beta_c$, and this should be used to show that, for such $\beta$, $\lim_{B\searrow 0}t_1(\beta,B)=0.$}}

\invisible{\Gib{Claudio: The symbol $\beta_c$ is used to define $\beta_c$. Remco: Fixed, thanks.}}

\subsection{Properties of touching points: Proof of Lemma \ref{lem-proof-touching-point}}
\label{sec-prop-touching-points}
\invisible{
We list here the properties of the function defined in \eqref{def-Fcal}: 
\RvdH{TO DO 5: Turn this into a proof of Lemma \ref{lem-proof-touching-point}}

\medskip
\begin{itemize}
    \item[$(a)$] $0\le \calF(t)<1$ for all $t\ge 0$;
    \item[$(b)$] $\calF(0)=0$ and $\calF(t)\to 1$ as $t\to +\infty$;
    \item[$(c)$] $\calF(t)$ is strictly increasing;
    \item[$(d)$] $\calF(t)$ is real analytic, with $\calF'(0)= \expec[W^2]/q\expec[W]>0$;
    \item[$(e)$] $t\calF'(t)\to 0$ as $t\to +\infty$;
    \item[$(f)$] $\calF(t)$ is convex close to $t=0$ and concave for $t$ large.
\end{itemize}
\invisible{\Gib{Analyticity is required in order to avoid an infinite number of solutions. Maybe is not true in general.}
\RvdH{It is not true in general, since higher-order moments will appear.}}
\medskip
We denote $\TP$ the set of touching points introduced in Def.\ref{def-touching-points} and also define $\TTP=\{t(\beta')\, |\, \beta'\in \TP\}$, the set of the points $t$ occurring in the definition of touching point. In the following we use the notation 
$L_{\beta'}(t)=t/\beta'$ and $D_{\beta'}(t)=L_{\beta'}(t)-\calF(t)$. Then \eqref{def-touching-point-t} is written as 
\begin{equation}\label{eq:touching-bis}
D_{\beta'}(t)=0,\, D'_{\beta'}(t)=0.
\end{equation}

\begin{lemma}\label{lem-touchingp}
    There exists a $\tilde{\beta'}>0$ such that for any $\beta'\in (0,\tilde{\beta'})$ there are no solutions to $D_{\beta'}(t)=0$ aside from $t=0$; then   $(0,\tilde{\beta'})\cap TP=\emptyset$. Moreover, there are solutions to $D_{\beta'}(t)=0$ different from $t=0$ for any $\beta'$ sufficiently large.
\end{lemma}
\begin{proof}
    From properties $(c)$ and $(e)$ we have that $\calF'(t)>0$ and $\calF'(t)\to 0$ as $t\to +\infty$. Then, since there exists $M=\max \{ \calF'(t) | t \ge 0 \} < +\infty$, we can take a $\tilde{\beta'}<M^{-1}$. Thus, $D'_{{\beta'}}(t)=1/{\beta'}-\calF'(t) \ge 1/\tilde{{\beta'}}-M>0$, for all $0<\beta'<\tilde{\beta'}$. This fact implies that $D_{{{\beta'}}}(t)>D(0)=0$ for all $t>0$, that is: $\beta'\notin \TP$.\\
    The second statement can be proven by taking a $\hat{t}>0$ and $\hat{\beta'}=\hat{t}/\calF(\hat{t})$. Then, for any $\beta'>\hat{\beta'}$ we have that $L_{\beta'}(\hat{t})< \calF(\hat{t})$, i.e. $D_{\beta'}(\hat{t})<0$. Since $D_{\beta'}(t) \to +\infty$ as $t\to +\infty$, we conclude that for all $\beta'>\hat{\beta'}$ there exists a solution to $D_{\beta'}(t)=0$ in the set $(\hat{t}, +\infty)$.
\end{proof}
}
In what follows we will make use of the fact that the function $t\mapsto \calF_0(t)$ is real analytic in any interval that does not contain the point $t=0$. We start by proving the first property stated in Lemma \ref{lem-proof-touching-point}:
\begin{lemma}\label{stetouch}The set of touching points $\TP$ contains at least two points and is finite.
\end{lemma}
\begin{proof}
 The set $\TP$ is not empty since $\beta_1'=1/\calF'(0)=q\expec[W]/\expec[W^2]$ is obviously a touching point. We now show that there exists at least a second touching point. \\
    It is not difficult to prove that the set of $\beta'>0$ such that $\calF (t)-t/{\beta'}\ne 0$ for all $t> 0$ is bounded from above, so that

\begin{equation}\label{eq-beta0}
        \beta'_0=\sup \left\{ \beta'>0 \colon \calF (t)-t/{\beta'}
        \ne 0\,~ \forall   t> 0 \right \}
    \end{equation}
    is finite. Since the set in the previous display is open (because of the continuity of $\calF(t)-t/\beta'$  with respect to $(\beta',t)$), it does not contain $\beta'_0$. Thus, there exists a $t_0>0$ such that $\calF(t_0)-t_0/{\beta'_0}=0$, i.e.,
    the first condition in the definition of touching point \eqref{def-touching-point-t} is satisfied. This means that the line through the origin with slope $1/\beta'_0$ intersects the graph $\calF(t)$ at the point $(t_0, \calF(t_0))$.
\smallskip
    
    The intersection discussed above is, in fact, tangential, i.e., $\calF'(t_0) = 1/\beta'$. Otherwise, if the intersection were transversal (i.e., $\calF'(t_0) \ne 1/\beta'$), then the solution would be persistent for $\beta'$ in a small neighborhood of $\beta'_0$, which is not possible by definition of $\beta'_0$ (see \eqref{eq-beta0}).
    Since the non-transversality 
    is the further requirement in the definition touching point, see the second equation in \eqref{def-touching-point-t}, we conclude that $\beta'_0$ is a touching point. Finally observe that $\beta'_0 \ne \beta'_1$. Indeed, since $\calF(t)$ is convex close to $t=0$, for $\beta'$ in some small left neighborhood of $\beta'_1$ there exist solutions to  $\calF(t)-t/\beta'=0$ other than $t=0$. This shows that $\beta'_0 \ne \beta'_1$ because, by definition of $\beta'_0$ (see \eqref{eq-beta0}), there are no solutions $t>0$ to $\calF(t)-t/\beta'_0=0$ for $\beta'<\beta'_0$. This proves that $\TP$ contains at least two points.
    
    We now prove that $\TP$ is finite. 
    Let us observe that if  $\beta_1,\, \beta_2 \in \TP$ and $\beta_1 \ne \beta_2$ the corresponding points $\theta(\beta'_1) \in \TTP$ and $\theta(\beta'_2)\in \TTP$ are different.
    Indeed, supposing that 
$\theta(\beta'_1)=\theta(\beta'_2)$, by the first equation in \eqref{def-touching-point-t}, we would have
    $$
\frac{\theta(\beta'_1)}{\beta'_1}=\calF(\theta(\beta'_1))=\calF(\theta(\beta'_2))=\frac{\theta(\beta'_2)}{\beta'_2},
    $$
    which implies that $\beta_1'= \beta_2'$.
 From this, it follows that if  $\TP$ were infinite, then $\TTP$ would also be infinite. In order to prove that the latter is finite, we start by observing that $\TTP$ is bounded. Indeed, if  $t \in \TTP$, by \eqref{def-touching-point-t}, it satisfies the  equation
    \begin{equation}\label{eq:ftfprime}
    \calF(t)-t \calF'(t)=0,
       \end{equation}
whose possible solutions are contained in an open interval $(\underline{t},\overline{t})$. 
    Indeed, since $ \calF(t)-t \calF'(t) \to 1$ as $t\to +\infty$, there exists a $0<\overline{t}<+\infty$ such that  $\calF(t)-t \calF'(t) > 0$ for $t>\overline{t}$, while, by the convexity of $\calF(t)$ close to zero (this is a consequence of the assumption $\expec[W^3]< \infty$), the solutions are bounded from below by some $\underline{t}>0$. Since the set of solutions to \eqref{eq:ftfprime} is bounded, it follows that it is also finite because  $\calF(t)-t \calF'(t)$ is analytic on $(\underline{t},\overline{t})$.
    We conclude that, as claimed, since $\TTP$ is finite, also $\TP$ is finite.
\end{proof}
\invisible{\begin{remark}
    From the previous discussion, it follows that 
    $\beta'_0 = \min \TP$. Since $\TP$ is finite there exists also $\max \TP$: is  the point $\beta'_1=1/\calF'(0)$ this maximum?
\end{remark}}
Let us denote the touching points  by
$$
\beta^{(t)}_1<\beta^{(t)}_2<\cdots < \beta^{(t)}_m,
$$
with $\beta^{(t)}_1 \equiv \beta'_0$ given in \eqref{eq-beta0}. According to Definition \ref{def-touching-points}, for each $\ell$ there is a point $t=\theta(\beta^{(t)}_\ell)>0$ satisfying \eqref{def-touching-point-t}, at which the sign of $\calF(t)-t/\beta^{(t)}_\ell$ does not change in $[\theta(\beta^{(t)}_\ell)-\delta, \theta(\beta^{(t)}_\ell)+\delta]$ for some $\delta>0$. As observed in Definition \ref{def-touching-points}, a sufficient condition for this is that $\min\{k\geq 2\colon  \mathscr{F}^{(k)}(\theta(\beta^{(t)}_\ell))\neq 0\}$ is even. 
On the other hand, we observe that  the existence of a finite $k\ge 2$ such that  $\mathscr{F}^{(k)}(\theta(\beta^{(t)}_\ell))\neq 0$ is guaranteed by the fact that the analytic function $\calF(t)$ is not identically zero in a neighborhood of $\theta(\beta^{(t)}_\ell)$.
\invisible{We remark here that the analyticity of $\calF(t)$ implies that the condition is also necessary.} The minimum of such $k$ is necessarily even, otherwise the sign would change in $[\theta(\beta^{(t)}_\ell)-\delta, \theta(\beta^{(t)}_\ell)+\delta]$ for any small $\delta>0$. Then, we have proven the following lemma:
\begin{lemma}\label{lem-even-deriv}
  For $\beta^{(t)}_\ell\in \TP$, there is at least one solution  $t=\theta(\beta^{(t)}_\ell)>0$ to \eqref{def-touching-point-t} such that $\inf\{k\geq 2\colon \calF^{(k)}(\theta(\beta^{(t)}_\ell))\neq 0\}$ is even.
\end{lemma}
Let us introduce the set
$$
S_{\beta'}=\{t>0\, \colon\,   t/\beta' - \calF(t)=0\}\quad \mbox{for}\quad  \beta^{(t)}_1\le \beta' \le \beta^{(t)}_m,
$$
and denote its cardinality by $|S_{\beta'}|$.
Since $t=0$ is a solution to $t/\beta' - \calF(t)=0$ for any $\beta'$, in order to complete the proof of Lemma \ref{lem-proof-touching-point}, we need to show that the solutions can appear or disappear only at the touching points via bifurcations (saddle-nodes), where they are created or annihilated in pairs:
\invisible{Next, we want to show that, if $\beta^{(t)}_1< \beta'< \beta^{(t)}_m$ but $\beta' \notin \TP$,  there is an odd number of solutions
\begin{equation}\label{eq:intesect_t}
 t_\ell(\beta'),\quad \ell=1,\ldots, 2\,k+1   
\end{equation}
(for some $k$) to the equation 
\begin{equation}\label{eq:equaz}
 \mathscr{F}(t)=\frac{t}{\beta'},
\end{equation}
that we write also as
\begin{equation}\label{eq:equaz-bis}
    D_{\beta'}(t)=0.
\end{equation}
}
\begin{lemma}\label{lem:numsol}
 $|S_{\beta'}|=2k(\beta')$, for some integer $k(\beta')\ge 1$, for all $\beta'< \beta^{(t)}_m$ and  $\beta' \notin \TP$. More precisely, $k(\beta')$ is constant in each open interval $(\beta^{(t)}_i,\beta^{(t)}_{i+1} )$, $i=1,\ldots, m-1$.
\end{lemma}
\begin{proof}
    In the following, for the sake of notation, we introduce the function $D(t,\beta')=t/\beta'-\calF(t)$ and denote by $D^{(k)}_{\beta'}(t,\beta')$ its $k$-th derivative with respect to $t$.
    
        Let us start by observing that, for $\beta'=\beta^{(t)}_1\equiv \beta'_0$,  the function $D^{(1)}(t,{\beta'_0})$ is zero at each point of the set $S_{\beta'_0}$.
    If not, there would be, by the Implicit Function Theorem, positive solutions to $D(t,\beta')=0$  for $\beta'$ slightly smaller than $\beta'_0$, which is not possible by definition of $\beta'_0$. This means that there are no solutions to $D(t,\beta'_0)=0$ other than those corresponding to the touching point, see Definition \ref{def-touching-points}. 
    Thus, by Lemma \ref{lem-even-deriv}, for each $\hat{t} \in S_{\beta'}$ there is an {\em even} $\hat{k}=k(\hat{t}) \ge 2$ such that 
    \be\label{eq-lim-point-sing}
    D(\hat{t},\beta'_0)=D^{(1)}(\hat{t},{\beta'_0}) =\cdots  =D^{(\hat{k}-1)}(\hat{t},{\beta'_0}) = 0, D^{(\hat{k})}(\hat{t},\beta'_0) \ne 0, \frac{\partial D(\hat{t},\beta'_0)}{\partial \beta'} \ne 0.
    \ee
    These conditions imply that $D(t,\beta')$ has a {\em limit point singularity} at $(\hat{t},\beta'_0)$ \cite{Golubitsky1985SingularitiesAG} and, as a consequence, the equation $D(t,\beta')=0$ undergoes a saddle-node bifurcation at $\beta'=\beta'_0$. This means that there exist two branches of solutions $(t^+(\beta'),\beta')$ and $(t^-(\beta'),\beta')$ that originate from $(\hat{t},\beta'_0)$ (i.e., $t^\pm(\beta'_0)=\hat{t}$) and lie locally on one side of $\beta'=\beta'_0$. 
    Since there are no solutions for $\beta'<\beta'_0$, see \eqref{eq-beta0}, the two branches lie necessarily on a right neighborhood  $[\beta'_0, \beta'_0+\delta)$ of the touching point; in this  neighborhood $|S_{\beta'}|$ is even. As $\beta'$ is increased, $|S_{\beta'}|$ may change since some of the solutions originating from $\beta'_0$ may disappear, while  new solutions may appear. We elaborate this issue further below.
    
 Let us consider a solution emerging from $(\hat{t},\beta'_0)$, e.g., $t^+(\beta')$. By applying the Implicit Function Theorem, $t^+(\beta')$ can be continued up to some  $\beta'^+>\beta'_0$ at which $D^{(1)}(t^+(\beta'^+),\beta'^+)$ $=0$. This means that $t^+(\beta')$ satisfies \eqref{def-touching-point-t} for $\beta'=\beta'^+_\ell$. 
 Moreover, we observe that, since $k(\beta'^+)=\min \{k\ge 0 : D^{(k)}(t^+(\beta'^+),\beta'^+)\ne 0\}$ is even, $\beta'^+$ is a touching point. 
  The fact that $k(\beta'^+)$ is even can be inferred from the fact that, otherwise, by Proposition 9.1 in \cite{Golubitsky1985SingularitiesAG}, the equation $D(t,\beta')=0$ would be locally equivalent to a normal form for which
 the solution $t^+(\beta')$ could be continued for $\beta'$ larger than $\beta'^+$, which is not possible by definition of $\beta'^+$.
 This shows that any solution  $t^+(\beta')$ may disappear by collapsing on some other solution, via a  {\em limit point singularity}, only at the touching point $\beta^{(t)}_2$. \invisible{$t^{+}_r(\beta')$ or $t^{-}_r(\beta')$.} In any case, the number of non-zero solutions  obtained by continuation of those existing  at $\beta'=\beta'_0$ remains even. 
 We finally observe that no new solution can appear for $\beta' \in (\beta'_0 \equiv \beta^{(t)}_1 , \beta^{(t)}_2)$ because this would imply that $\beta'$ is a touching point. Therefore $|S_{\beta'}|$ is even and constant on $(\beta^{(t)}_1 , \beta^{(t)}_2)$.
 
The previous argument can be repeated for the other intervals $(\beta^{(t)}_{\ell-1}, \beta^{(t)}_{\ell}),\,  \ell=2,\ldots, m$, showing that there is an even number $2k(\beta')$ of non-zero solutions, which is constant  on each of them.
\end{proof}
\begin{remark}
[Largest touching point]
{\rm
Special care is required for the largest touching point $\beta'_1=1/\calF'(0)=q\expec[W]/\expec[W^2]$. Indeed, this point corresponds to a singularity at which a {\em transcritical bifurcation} takes place
since $D^{(2)}(0,\beta'_1)<0$ and $\partial^2 D(0,\beta'_1)/\partial t \partial \beta'<0$, see \cite{Golubitsky1985SingularitiesAG}. At $\beta'=\beta'_1$ a solution $t_2(\beta')$, which is positive for $\beta'<\beta'_1$, merges with $t_1(\beta')\equiv 0$ and continues to exist also for $\beta'$ slightly larger than $\beta'_1$, but taking negative values. Thus, for $\beta'$ larger than the touching point $\beta'_1$, the number of solutions of $\mathscr{F}(t)=t/\beta'$ is even.}\hfill\ensymboldefinition
\end{remark}
\begin{remark}[Multiple solutions at a touching point]{\rm
    In the proof of Lemma \ref{lem:numsol} we have implicitly assumed that for each
    touching point  $\beta^{(t)}_i$
    there is  only one solution $\theta(\beta^{(t)}_i)$ of \eqref{def-touching-points}. In the case of more solutions, the proof of the lemma does not change. Indeed, it is sufficient to note that at all $\theta(\beta^{(t)}_i)$'s corresponding to the same $\beta^{(t)}_i$, only a limit point singularity can occur; then a {\em couple} of solutions may appear or disappear. In any case, the number of solutions of $\mathscr{F}(t)=t/\beta'$ (excluding the trivial one $t=0$)  remains even.
    This remark also applies to the proof of Theorem \ref{thm-general-first-order}.}
    \hfill\ensymboldefinition
\end{remark}

\invisible{\footnote{For a smooth function $g(x,\lambda)$ to undergo a {\em trancritical bifurcation} at $(0,0)$ it is sufficient that at this point
\begin{equation}\label{eq:cns-transcritical}
 g=\frac{\partial }{\partial x}g=\frac{\partial }{\partial \lambda} g=0,\quad \mbox{} \quad \varepsilon:=\frac{\partial^{2} }{\partial x^{2}}g \ne 0,\; \delta:=\frac{\partial^2 }{\partial x \partial \lambda }g \ne 0.
 \end{equation}
}}

\invisible{
We start by explicitly computing the average in the stationarity condition \eqref{FOPT-1-B1-B>0} as
\eqan{
\label{Pareto-1}
\mathbb{E}\left[ \frac{W}{\expec[W]} \frac{\e^{\beta' W s}-1}{\e^{\beta' W s}+q-1}\right] &=(\tau-2)\int_1^\infty w^{-(\tau-1)}
\frac{\e^{\beta' w s}-1}{\e^{\beta' w s}+q-1}{\rm d}w.
}
We use that 
    \eqn{
    \beta' s\frac{\e^{\beta' w s}-1}{\e^{\beta' w s}+q-1}   =\frac{{\rm d}}{{\rm d} w} \Big[\frac{q}{q-1}\log(\e^{\beta' w s}+q-1)-\frac{s\beta'w}{q-1}\Big].
    }
Therefore, by partial integration,
    \eqan{
    \label{Pareto-2}
    \mathbb{E}\left[ \frac{W}{\expec[W]} \frac{\e^{\beta' W s}-1}{\e^{\beta' W s}+q-1}\right]
    &=\frac{\tau-2}{\beta's}\Big[w^{-(\tau-1)}\big[\frac{q}{q-1}\log(\e^{\beta' w s}+q-1)-\frac{s\beta'w}{q-1}\big]\Big]_{w=1}^\infty\nn\\
    &\qquad +\frac{(\tau-1)(\tau-2)}{\beta' s}\int_1^\infty w^{-\tau}
    \Big[\frac{q}{q-1}\log(\e^{\beta' w s}+q-1)-\frac{s\beta'w}{q-1}\Big]{\rm d}w\nn\\
    &=-\frac{\tau-2}{\beta's}\big[\frac{q}{q-1}\log(\e^{\beta' s}+q-1)-\frac{s\beta'}{q-1}\big]\nn\\
    &\quad +\frac{q(\tau-2)}{(q-1)\beta' s}\mathbb{E}\left[\log\left(\e^{\beta' W s}+q-1\right) \right]
    -\frac{\tau-1}{q-1}\nn\\
    &=-\frac{(\tau-2)q}{(q-1)\beta's}\log(\e^{\beta' s}+q-1)+\frac{\tau-2}{q-1}-\frac{\tau-1}{q-1}\nn\\
    &\quad+\frac{q(\tau-2)}{(q-1)\beta' s}\mathbb{E}\left[\log\left(\e^{\beta' W s}+q-1\right) \right]\nn\\
    &=-\frac{(\tau-2)q}{(q-1)\beta's}\log(\frac{\e^{\beta' s}+q-1}{q})-\frac{1}{q-1}\nn\\
    &\quad+\frac{q(\tau-2)}{(q-1)\beta' s}\mathbb{E}\left[\log\left(\frac{\e^{\beta' W s}+q-1}{q}\right) \right]
    }
By the stationarity condition \eqref{FOPT-1-B1-B>0}, this equals $s$. By the criticality condition \eqref{crit-cond-B>0}
\be
\label{FOPT-2-rep}
\mathbb{E}\left[\frac{1}{\mathbb{E}[W]}\log\left(\frac{\e^{\beta' W s}+q-1}{q}\right) \right] = \frac{\beta'}{q} \left(s + \frac12 (q-1) s^2\right).
\ee 
Combining this gives
\be
s= -\frac{(\tau-2)q}{(q-1)\beta's}\log(\frac{\e^{\beta' s}+q-1}{q})-\frac{1}{q-1}+\frac{q(\tau-1)}{(q-1)\beta' s}\frac{\beta'}{q} \left(s + \frac12 (q-1) s^2\right)
\ee
or
\be
\frac{(q-1)\beta's}{(\tau-2)q}s= -\log(\frac{\e^{\beta' s}+q-1}{q})-\frac{\beta's}{(\tau-2)q}+\frac{(\tau-1)}{(\tau-2)}\frac{\beta'}{q} \left(s + \frac12 (q-1) s^2\right)
\ee
We rewrite this equation as
    \eqn{
    \label{FOPT-Pareto}\log\left(\frac{\e^{\beta' s}+q-1}{q}\right)
    =\frac{\beta' s}{q}+\frac{(q-1)\beta' s^2}{2q} \frac{\tau-3}{\tau-2}.
    }
or, equivalently,
\be
\log\left(\frac{\e^{\beta' s}+q-1}{q}\right)
    =\frac{\beta' s}{q}+\frac{(q-1)\beta' s^2}{2q}\frac{\mathbb{E}[W]}{\mathbb{E}[W^2]}.
\ee
Apart from multiplicative constants, this is similar to the second equation \eqref{crit-cond-B>0} in the  homogeneous case. However, it needs to be complemented by a second equation that generalizes \eqref{FOPT-1-B1-B>0}.
\RvdH{If we use the rewrite in \eqref{Pareto-2} in the formula for $t_c$ using $\mathscr{K}(t)$, then \eqref{Pareto-2} obtains a clearer purpose, and we would obtain that $t_c$ (and thus $\beta_c$ and $s(\beta_c)$ satisfy ONE equation that looks a little like the homogeneous case. Would it even be the same?}
}

\invisible{\medskip

\paragraph{\bf On the nature of the fixed-point equation \eqref{FOPT-Pareto}.} Let
    \eqan{
    f(s)&=\log\left(\frac{\e^{\beta' s}+q-1}{q}\right),\\
    g(s)&=\frac{\beta' s}{q}\frac{\tau-3}{\tau-2}+\frac{(q-1)\beta' s^2}{2q} \frac{\tau-4}{\tau-2}.
    }
We note that
    \eqn{
    f'(s)=\frac{\beta'\e^{\beta' s}}{\e^{\beta' s}+q-1},\qquad
    f''(s)=\frac{(\beta')^2(q-1) \e^{-\beta' s}}{(1+(q-1)\e^{-\beta' s})^2}.
    }
As a result,
    \eqan{
    (f(s)-g(s))''&=\beta'(q-1)\Bigg[\frac{\beta'\e^{-\beta' s}}{(1+(q-1)\e^{-\beta' s})^2}-\frac{\tau-4}{q(\tau-2)}\Bigg]\nn\\
    &=\beta'(q-1)\Bigg[\frac{\beta'}{(\e^{\beta' s}+q-1)(1+(q-1)\e^{-\beta' s})}-\frac{\tau-4}{q(\tau-2)}\Bigg].
    }
Note that
    \eqn{
    (\e^{\beta' s}+q-1)(1+(q-1)\e^{-\beta' s})=\e^{\beta' s}+2(q-1)+(q-1)^2\e^{-\beta' s},
    }
which is first decreasing and then increasing. As a result, $s\mapsto (f(s)-g(s))''$ is first increasing and then decreasing. Further,
    \eqn{
    f''(0)-g''(0)=\beta'(q-1)\Big[\frac{\beta'}{q^2}-\frac{\tau-4}{q(\tau-2)}\Big]<0,
    }
since $\beta_c'<q(\tau-4)/(\tau-2).$ Thus, either $s\mapsto (f(s)-g(s))''$ is negative for all $s$, or starts negative, then becomes positive, and finally becomes negative again. Only in the latter case, there can be a positive solution to $f(s)-g(s)=0$. Call this solution $s(\beta).$ Then, $\beta_c$ is characterized by either $p(s(\beta),\beta)=p(0,\beta),$ or
$s(\beta_c)$ should solve
    \eqn{
    s(\beta_c)=\expec\Bigg[\frac{W}{\expec[W]}\frac{\e^{\beta' s(\beta_c)}-1}{\e^{\beta' s(\beta_c)}+q-1}\Bigg].
    }
It would be helpful if we could prove that $\beta\mapsto s(\beta)$ is decreasing, while $\beta\mapsto \beta s(\beta)$ is increasing.
\RvdH{TO DO 11: Can we move closer to identifying the critical value for Pareto distributions, using \eqref{FOPT-Pareto}?}}

\invisible{
For this, it turns out to be useful to investigate the homogeneous equation, arising when $\prob(W=w)=1$. Before we were taking $w=1$, which we can assume w.l.o.g.\, but the extra $w$ variable gives an extra degree of freedom. Then, 
\eqref{FOPT-1-B1-B>0}-\eqref{crit-cond-B>0} become
\be
\label{FOPT-1-w}
\frac{\e^{\beta' w s}-1}{\e^{\beta' w s}+q-1}= s,
\ee
and
\be
\label{FOPT-2-w}
\log\left(\frac{\e^{\beta'w s}+q-1}{q}\right) = \frac{\beta'}{q} \left(s + \frac12 (q-1) s^2\right).
\ee 
Redefining $\tilde{s}=s w$, this becomes
\be
\label{FOPT-1-w-2}
\frac{\e^{\beta' \tilde{s}}-1}{\e^{\beta' s}+q-1}= \frac{\tilde{s}}{w},
\ee
and
\be
\label{FOPT-2-w-2}
\log\left(\frac{\e^{\beta' \tilde{s}}+q-1}{q}\right) = \frac{\beta'}{q}\frac{\tilde{s}}{w} \left(1 + \frac12 (q-1) \frac{\tilde{s}}{w}\right).
\ee}

\subsection{General first-order phase transition: Extension to $B>0$}
\label{sxec-ext-B>0-general-first-order-phase transition}
In this section, we extend the analysis in the previous section to $B>0$. 
\invisible{For this, we define
    \eqn{
    \label{def-Fcal-B>0}
    \mathscr{F}_B(t)=\expec \left[\frac{W}{\expec[W]}\frac{\e^{tW+B}-1}{\e^{tW+B}+q-1} \right].}}
A significant difference with the case where $B=0$ is that now
    \eqn{
    \mathscr{F}_B(0)=\frac{\e^B-1}{\e^B+q-1}>0,
    }
rather than $\mathscr{F}(0)=\mathscr{F}_0(0)=0$. However, it is true that $\mathscr{F}_B(0)\searrow 0$ as $B\searrow 0$.
\invisible{\smallskip
We can compute that
    \be
    \label{eq:first_derivative_gen-B>0}
    \frac{d}{dt}\mathscr{F}_B(t)=\frac{q}{\expec[W]}\expec\left[W^2\frac{\e^{tW}}{(\e^{tW}+\e^{-B}(q-1))^2}\right],
    \ee
and
    \be
    \label{eq:second_derivative_gen-B>0}
    \frac{d^2}{dt^2}\mathscr{F}_B(t)=\frac{q}{\expec[W]}\expec\left[W^3\frac{\e^{-B}(q-1)\e^{tW}-\e^{2tW}}{(\e^{tW}+\e^{-B}(q-1))^3}\right].
    \ee}
From \eqref{eq:second_derivative_gen}, we conclude that $\frac{d^2}{dt^2}\mathscr{F}_B(0)>0$ when $\expec[W^3]<\infty$, $q>2$ and $B<\log(q-1)$, and thus
    \eqn{
    \label{t-conv-def-B>0}
    t_{\rm conv}(B)=\sup\{s \colon \mathscr{F}_B(u) \text{ is convex on }(0,s)\}>0.
    }
We note that $\lim_{B\searrow 0}t_{\rm conv}(B)=t_{\rm conv}>0.$
\smallskip

The optimal $s(\beta,B)$ satisfies the stationarity condition in \eqref{FOPT-1-B1-B>0}, which is equivalent to $s(\beta,B)=t(\beta,B)/\beta'$, where $t(\beta,B)$ solves
    \eqn{
    \label{F-B-fixed-point}
    \mathscr{F}_B(t)=\frac{t}{\beta'}.
    }
As before, there may be several such solutions, and we have to pick the one that optimizes the variational problem in \eqref{var-problem-psi-B>0}.
\smallskip

Since $\mathscr{F}_B(0)>0$, it is not immediately obvious what the first solution to \eqref{F-B-fixed-point} is. However, for $1/\beta'>\frac{d}{dt}\mathscr{F}_0(0)=\expec[W^2]/(q\expec[W])$ and $B>0$ sufficiently small, there will be a solution $t_1(\beta,B)>0$ that converges to zero when $B\searrow 0.$ This condition on $\beta'$ is equivalent to
    \eqn{
    \beta'< \frac{q\expec[W]}{\expec[W^2]}.
    }
\smallskip

We can extend the analysis in Section \ref{sec-FO-PT-general} to show that $\beta\mapsto t(\beta,B)$ makes a jump at some $\beta_c(B)$, leading to a first-order phase transition at this value of $\beta$. Further, the jump is beyond $t_{\rm conv}(B)>0$. However, for $B=0$, the jump is from $t_1(\beta)=0,$ while for $B>0$, the jump comes from  $t_1(\beta,B)$, which is the smallest positive solution to $\mathscr{F}_B(t)=t/\beta'$, whose existence we have proved above.

The above analysis prompts us to define the {\em critical inverse temperature} as
\eqn{
\label{criticalwithB}
\beta_c(B)=\inf\{\beta\colon  \lim_{\tilde\beta\nearrow\beta} s(\tilde\beta,B) \neq \lim_{\tilde\beta\searrow\beta} s(\tilde\beta,B)\}.
}
However, since we cannot show that the jump is {\em unique}, we cannot rule out that there exist other critical inverse temperatures larger than $\beta_c(B)$.

\bigskip

\noindent
{\bf Acknowledgments.}
C. Giardin\`a and C. Giberti are members of Gruppo Nazionale di Fisica Matematica (GNFM) of Istituto Nazionale di Alta Matematica (INdAM). C. Giardin\`a acknowledge 
support by the INdAM Project 2024, Cup E53C23001740001.
The work of RvdH was supported in part by the Netherlands Organisation for Scientific Research (NWO) through the Gravitation {\small{\sf NETWORKS}} grant no.\ 024.002.003. Most of the work of NM was done while he was a PhD student at the Department of Mathematics and Computer Science at Eindhoven University of Technology.


\newpage

\appendix

\section{Counterexample of weights with two atoms}
\label{appA}



{\color{black} Proposition \ref{prop:unique_inflection_for_bdd_supp} shows that for $W$ with support inside a {\em small} compact set, the unique zero-crossing Condition \ref{cond-zer-cross}(a), holds. Here we show that  the boundedness of the support is not sufficient for this condition  to hold. More precisely, for the case $B=0$, we give a counterexample to the statement that there is a unique $t_*>0$ such that $\frac{d^2}{dt^2}\mathscr{F}_0(t)>0$ for $0<t<t_*$ and $\frac{d^2}{dt^2}\mathscr{F}_0(t)<0$ for $t>t_*$.}
\smallskip


For the sake of notation, we write $H(t)=\frac{d^2}{dt^2}\mathscr{F}_0(t)$. By\eqref{eq:second_derivative_gen},
	\eqn{
	H(t)= q\expec\left[\frac{W^3}{\expec[W]}a(W t)\right], \qquad t \geq 0,
	}
where
	\eqn{
	a(s)=\frac{ p \e^{s}-\e^{2s}}{\left(\e^{s}+p\right)^{3}}, \qquad s \geq 0,
	}
and $p=q-1=2,3, \ldots$ 
\smallskip

 
We have that $a(s)=0$ precisely when $s=\log p$, and $a(s)>0$ for $s \in(0, \log p)$ and $a(s)<0$ for $s \in(\log p, \infty)$. We consider the case $p=6$, or $q=7$. We have, for $t>0$,
	\eqn{
	\label{Guido(6-2)}
	H(t)=
    q\expec\Big[\frac{W^3}{\expec[W]}a(Wt)\Big]=t^{-3}\expec\big[ F(W t)\big],\qquad
    \text{where}\qquad
	F(s)=s^{3} a(s).
	}
Thus,
 	\eqn{
	\label{Guido(6-4)}
	t^{3} H(t)=\expec\left[ F(W t)\right], \qquad t>0.
	}
We let the distribution of $W$ be $\prob(W=x_1)=c_{1}$ and $\prob(W=x_2)=c_{2},$ 
with $c_{1}, c_{2}>0, c_{1}+c_{2}=1$, and $0<x_{1}<x_{2}$. 
\invisible{While this $g$ is not a true pdf, we can approximate the two delta functions by a smooth, non-negative function $H$ supported by an interval around $x_{1}, x_{2}$ of length $10^{-10}$, say, and of unit integral. This should
approximate the integral of the right-hand side of \eqref{Guido(6-4)} with an absolute accuracy of order $10^{-10}$ (and better) for the relevant range $0 \leq t \leq 10$. 
\RvdH{Remco: It need not be a density, as long as it is a distribution.}}
With this choice of $W$,
	\eqn{
	\label{Guido(6-6)}
	t^{3} H(t)=c_{1} F\left(x_{1} t\right)+c_{2} F\left(x_{2} t\right).
	}
We assume $x_{1}=1, x_{2}=5$ and, {\color{black} observing that $F(1)>0, F(5)<0$, } we can choose $c_{1}, c_{2}>0$ with $c_{1}+c_{2}=1$ such that
	\eqn{
	\label{Guido(6-7)}
	H(1)=c_{1} F(1)+c_{2} F(5)=0.
	}
    More precisely, with
	\eqn{
	\label{Guido(6-8)}
	F(1)=0.013461797..., \quad F(5)=-0.717595354...,
	}
we find
	\eqn{
	\label{Guido(6-9)}
	c_{1}=0.98158584..., \quad c_{2}=1-c_{1}=0.018414151... .
	}
Since $F$ increases at $s=1$ and $s=5$, we  have that $t^3H(t)=c_{1} F(t)+c_{2} F(5 t)$ is positive for $t$ a little larger than 1, and negative for $t$ a little smaller than 1. Thus $H(t)$ has an upcrossing at $t=1$.

\invisible{\RvdH{Remove this table? The table is sort of replaced by Figure 2, right?}
A short table of values of $c_{1} F(t)+c_{2} F(5 t):$
$$
\begin{array}{|l|l|}
\hline
t & c_{1} F(t)+c_{2} F(5 t) \\
\hline
0.0 & 0.000000000 \\
0.3 & +7.862222746 \times 10^{-4} \\
0.5 & -1.637906998 \times 10^{-3} \\
0.8 & -6.466733494 \times 10^{-3} \\
0.9 & -4.043721218 \times 10^{-3} \\
1.0 & 0.000000000 \\
1.1 & +4.749780200 \times 10^{-3} \\
1.5 & +0.015335392 \\
1.8 & -2.63450994 \times 10^{-3} \\
2.0 & -0.034414816\\
\hline
\end{array}
$$
shows this.}
\smallskip


On the other hand, we note that $t^3H(t)=c_{1} F(t)+c_{2} F(5 t)>0$ for $0<t<\frac{1}{5} \log 6$ and $t^3H(t)= c_{1} F(t)+c_{2} F(5 t)<0$ for $t>\log 6$. This means that, since $H(1)=0$, $t^3H(t)$ crosses zero once between $(\frac{1}{5}\log 6, 1)$. 
{\color{black} We finally observe that this example does not violate Theorem \ref{thm-two-atoms}, since the small support condition \eqref{eq:small_support_condition} is not satisfied in this example, since   $5>1+\frac{\log (2+\sqrt{3})}{\log(q-1)}$ for $q=7$. }


\invisible{
$$
\begin{aligned}
& \operatorname{sgn}\left(t^{3} \mathscr{F}(t)\right) \\
&+~+~+~+~+~+~-~-~-~\cdots-~+~+~+~+~+~+~+~-~-~-~-\\
& x_{t} \in(0, \tfrac15 \log 6), \quad x_{t} \in(\tfrac15 \log 6, 1), \qquad x_{t} \in(1, \log 6), \qquad\quad  x_{t} \in(\log 6, \infty)
\end{aligned}
$$}


\section{Smoothing critical exponent: Proof of Proposition \ref{prop-r(q)-asymp}}
\label{appB}

\invisible{Here we restrict to $B=0$. Recall from \eqref{Guido(5-2)} that}
For $q>2$, and writing $x_0=\log{(q-1)}$, we have from \eqref{Guido(5)}
    \eqan{
    a(x)&=\frac{(q-1) \e^{x}-\e^{2 x}}{\left(\e^{x}+q-1\right)^{3}}=\frac{\e^{x_{0}+x}-\e^{2 x}}{\left(\e^{x}+\e^{x_{0}}\right)^{3}}\\
    &=\frac{\left(\e^{\frac{1}{2}\left(x_{0}-x\right)}-\e^{-\frac{1}{2}\left(x_{0}-x\right)}\right) \e^{\frac{1}{2} x_{0}+\frac{3}{2} x}}{\left(\e^{\frac{1}{2}\left(x-x_{0}\right)}+\e^{-\frac{1}{2}\left(x-x_{0}\right)}\right)^{3} \e^{\frac{3}{2} x+\frac{3}{2} x_{0}}}
    \nn\\
    &=\frac{1}{4(q-1)} \frac{\sinh \frac{1}{2}\left(x_{0}-x\right)}{\cosh ^{3} \frac{1}{2}\left(x_{0}-x\right)}.\nn
    }
By \eqref{r(q)-equality}, $r(q)$ satisfies
    \eqn{
    \int_0^\infty x^ra(x)dx=0.
    }
This can be rewritten as
    \eqn{
    \label{split-integral-a}
    0=\int_{0}^{2 x_{0}} x^{r} a(x) d x+\int_{2 x_{0}}^{\infty} x^{r} a(x) d x,
    }
so that the point-symmetry of $a(x)$ about the point $x=x_{0}$ can be taken advantage of.

It appears in the course of the proof of \eqref{r(q)-asymptoics-sharp} that we need to work with $r \in(-1,0)$ away from -1. Thus, we first present a lower bound $\ell(q)$ for $r(q)$ of the form
    \eqn{
    \label{equation-Guido-11}
    r(q)>\ell(q)=\left(\left(\frac{3}{2}+\frac{q-2}{q} x_{0}\right)^{2}-2\right)^{1 / 2}-\left(\frac{3}{2}+\frac{q-2}{q} x_{0}\right).
    }
The function $\ell(q)$ increases in $q$ from -1 at $q=2$ to 0 at
$q=\infty$. Hence, since for $q \geq 3$,
    \eqn{
    \label{equation-Guido-12}
    r(q)>\ell(q) \geq \ell(3)=\left(\left(\frac{3}{2}+\frac{1}{3} \log 2\right)^{2}-2\right)^{1 / 2}-\left(\frac{3}{2}+\frac{1}{3} \log 2\right)=-0.7327 \ldots ,
    }
we restrict attention to $r \geq-3 / 4$.

To show \eqref{equation-Guido-11}, for the first term of the right-hand side of \eqref{split-integral-a}, we have
    \eqn{
    \label{equation-Guido-13}
    \int_{0}^{2 x_{0}} x^{r} a(x) d x=x_{0}^{r} \int_{0}^{x_{0}} b(w)\left(\left(1-\frac{w}{x_{0}}\right)^{r}-\left(1+\frac{w}{x_{0}}\right)^{r}\right) d w,
    }
where, for $w\geq 0$,
    \eqn{
    \label{equation-Guido-14}
    b(w)=\frac{1}{4(q-1)} \frac{\sinh(w/2)}{\cosh(w/2)^3}.
    }
It follows from some elementary considerations that, for any $q=3,4, \ldots$ and all $0 \leq w \leq x_0$,
    \eqn{
    \label{equation-Guido-15}
    b(w) \geq \frac{1}{4(q-1)} \frac{w}{x_0} \frac{\sinh(x_0/2)}{\cosh(x_0/2)^3}=
    \frac{q-2}{q^3} \frac{w}{x_0}.
    }
Hence,
    \eqan{
    \label{equation-Guido-16}
    \int_{0}^{2 x_{0}} x^{r} a(x) dx 
    &\geq  \frac{q-2}{q^{3}} x_{0}^{r} \int_{0}^{x_{0}} \frac{w}{x_{0}}\left(\left(1-\frac{w}{x_{0}}\right)^{r}-\left(1+\frac{w}{x_{0}}\right)^{r}\right) d w\\
    &=\frac{q-2}{q^{3}} x_{0}^{r+1} \int_{0}^{1} t\left((1-t)^{r}-(1+t)^{r}\right) d t.\nn
    }
Then from
    \eqn{
    \label{equation-Guido-17}
    \int_{0}^{1} t\left((1-t)^{r}-(1+t)^{r}\right) d t=\frac{-r 2^{r+1}}{(r+2)(r+1)},
    }
we obtain
    \eqn{
    \label{equation-Guido-18}
    \int_{0}^{2 x_{0}} x^{r} a(x) d x \geqslant-r \frac{q-2}{q^{3}} \frac{\left(2 x_{0}\right)^{r+1}}{(r+2)(r+1)}.
    }

Next, for the second term at the right-hand
side of \eqref{split-integral-a}, we have 
    \eqan{
    \label{equation-Guido-19}
    \int_{2 x_{0}}^{\infty} x^{r} a(x) dx
    &=\int_{2 x_{0}}^{\infty} x^{r} \frac{d}{d x}\left(\frac{\e^{x}}{\left(\e^{x}+q-1\right)^{2}}\right) d x\\
    &=
    \frac{x^{r} \e^{x}}{\left(\e^{x}+q-1\right)^{2}}\Big|_{x=2x_0}^{\infty}-r \int_{2 x_{0}}^{\infty} \frac{x^{r-1} \e^{x}}{\left(\e^{x}+q-1\right)^{2}} d x>-\frac{\left(2 x_{0}\right)^{r}}{q^{2}},\nn
    }
where we have used that $\e^{x_{0}}=q-1$ and $r<0$.

Combining \eqref{split-integral-a}, \eqref{equation-Guido-18} and \eqref{equation-Guido-19}, we thus find
    \eqn{
    \label{equation-Guido-20}
    0>-r \frac{q-2}{q^{3}} \frac{\left(2 x_{0}\right)^{r+1}}{(r+2)(r+1)}-\frac{\left(2 x_{0}\right)^{r}}{q^{2}},
    }
and so
    \eqn{
    \label{equation-Guido-21}
    -r<\frac{\displaystyle\frac{\left(2 x_{0}\right)^r}{q^{2}}}{\displaystyle\frac{q-2}{q^{3}} \frac{\left(2 x_{0}\right)^{r+1}}{(r+2)(r+1)}}=\frac{q}{2(q-2) x_{0}}(r+2)(r+1),
    }
i.e.,
    \eqn{
    \label{equation-Guido-22}
    r^{2}+3 r+2+2\frac{q-2}{q} x_{0} r>0.
    }

The left-hand side of \eqref{equation-Guido-22} vanishes if and only if
    \eqn{
    \label{equation-Guido-23}
    r=r_{ \pm}= \pm\left(\left(\frac{3}{2}+\frac{q-2}{q} x_{0}\right)^{2}-2\right)^{1 / 2}-\left(\frac{3}{2}+\frac{q-2}{q} x_{0}\right).
    }
We have, for $q>2$,
    \eqn{
    \label{equation-Guido-24}
    r_{-}<-2, 
    \qquad r_{+}=\frac{-2}{\left(\frac{3}{2}+\frac{q-2}{q} x_{0}\right)+\left(\left(\frac{3}{2}+\frac{q-2}{q} x_{0}\right)^{2}-2\right)^{1 / 2}} \in(-1,0).
    }
Then from \eqref{equation-Guido-22} and $r=r(q) \in(-1,0)$, we obtain $r(q) \in$ $\left(r_{+}, 0\right)$, and this gives \eqref{equation-Guido-11}. Observe also from \eqref{equation-Guido-24} and $x_{0}=\log (q-1)$ that $r_{+}$ decreases in $q \geq 2$.
\smallskip

To proceed showing \eqref{r(q)-asymptoics-sharp}, we bound the two terms at the right-hand side of \eqref{split-integral-a} in yet another way. As to the first term at the right-hand side of \eqref{split-integral-a}, we observe, with an eye on \eqref{equation-Guido-13}, that, for $0 \leq w<x_{0}$,
    \eqan{
    \label{equation-Guido-25}
    \left(1-\frac{w}{x_{0}}\right)^{r}-\left(1+\frac{w}{x_{0}}\right)^{r}&=-2 r \frac{w}{x_{0}}-2 \sum_{k=1}^{\infty}\binom{r}{2 k+1}\left(\frac{w}{x_{0}}\right)^{2 k+1}\\
    &=-2 r \frac{w}{x_{0}}-2 r \sum_{k=1}^{\infty} d_{k}\left(\frac{w}{x_{0}}\right)^{2 k+1},\nn
    }
where
    \eqn{
    \label{equation-Guido-26}
    d_{k}=\frac{1}{r}\binom{r}{2 k+1}=\frac{(1-r)(2-r) \cdots (2 k-r)}{2 \cdot 3 \cdots (2 k+1)} \in(0,1), 
    \qquad k=1,2, \ldots.
    }
Therefore, from \eqref{equation-Guido-13},
    \eqn{
    \label{equation-Guido-27}
    \int_{0}^{2 x_{0}} x^{r} a(x) d x
    =-2 r x_{0}^{r-1} \int_{0}^{x_{0}} w b(w) d w-2 r x_{0}^{r} \int_{0}^{x_{0}} b(w) \sum_{k=1}^{\infty} d_{k}\left(\frac{w}{x_{0}}\right)^{2 k+1} d w.
    }
We evaluate
    \eqn{
    \label{equation-Guido-28}
    -2 r x_{0}^{r-1} \int_{0}^{x_{0}} w b(w) d w=-2 r x_{0}^{r}\left(\frac{q-2}{2 q(q-1) x_{0}}-\frac{1}{q^{2}}\right).
    }
Furthermore, 
we have, see \eqref{equation-Guido-14},
    \eqn{
    \label{equation-Guido-30}
    0<b(w) \leqslant \frac{1}{q-1} \e^{-w}, \qquad \frac{1}{2} x_{0} \leq w \leq x_{0}.
    }
Therefore, see \eqref{equation-Guido-25},
    \eqan{
    \label{equation-Guido-31}
    0&\leqslant-2 r \int_{x_{0} / 2}^{x_{0}} b(w) \sum_{k=1}^{\infty} d_{k}\left(\frac{w}{x_{0}}\right)^{2 k+1} dw\\
    &\quad=\int_{x_{0} / 2}^{x_{0}} b(w)\left(\left(1-\frac{w}{x_{0}}\right)^{r}-\left(1+\frac{w}{x_{0}}\right)^{r}+2 r \frac{w}{x_{0}}\right) d w\nn\\
    &\quad\leq
    \frac{1}{q-1} \e^{-x_{0} / 2} \int_{x_{0} / 2}^{x_{0}}\left(\left(1-\frac{w}{x_{0}}\right)^{r}-\left(1+\frac{w}{x_{0}}\right)^{r}+2 r \frac{w}{x_{0}}\right) d w\nn\\
    &\quad =\frac{x_{0}}{(q-1)^{3/2}} \int_{1 / 2}^{1}\left((1-t)^{r}-(1+t)^{r}+2 r t\right) d t.\nn
    }
We estimate
    \eqan{
    \label{equation-Guido-32}
    0 &\leq \int_{1 / 2}^{1}\left((1-t)^{r}-(1+t)^r+2 r t\right) d t \leq \int_{0}^{1}\left((1-t)^{r}-(1+t)^{r}+2 r t\right) dt\\
    &=r+2 \frac{1-2^{r}}{r+1}<r-\frac{2 r \log 2}{r+1}=-r\left(\frac{2 \log 2}{r+1}-1\right)<-r(8 \log 2-1),\nn
    }
where we have used that $1-2^{r}<-r \log 2$ and that we have restricted our analysis to the range $r\in \left(-\frac{3}{4}, 0\right)$. Hence,
    \eqn{
    \label{equation-Guido-33}
    0 \leq -2 r \int_{x_{0} / 2}^{x_{0}} b(w) \sum_{k=1}^{\infty} d_{k}\left(\frac{w}{x_{0}}\right)^{2 k+1} d w \leqslant-r \frac{x_{0}}{(q-1)^{3 / 2}}(8 \log 2-1).
    }
At the same time, for $0 \leq w \leq x_{0}/2,$
    \eqan{
    \label{equation-Guido-34}
    -2 r \sum_{k=1}^{\infty} d_{k}\left(\frac{w}{x_{0}}\right)^{2 k+1}&=-2 r\left(\frac{w}{x_{0}}\right)^{3} \sum_{k=1}^{\infty} d_{k}\left(\frac{w}{x_{0}}\right)^{2 k-2}\\
    &\leq -2 r\left(\frac{w}{x_{0}}\right)^{3} \sum_{k=1}^{\infty}\left(\frac{1}{2}\right)^{2 k-2}=-\frac{8}{3} r\left(\frac{w}{x_{0}}\right)^{3},\nn
    }
where we have used that $d_{k} \in(0,1)$, see \eqref{equation-Guido-26}. Therefore,
    \eqan{
    \label{equation-Guido-35}
    -2 r \int_{0}^{x_{0} / 2} b(w) \sum_{k=1}^{\infty} d_{k}\left(\frac{w}{x_{0}}\right)^{2 k+1} 
    &\leq -\frac{8}{3} r \int_{0}^{x_{0} / 2} \frac{\e^{-w}}{q-1}\left(\frac{w}{x_{0}}\right)^{3} d w
    \leq 
    \frac{-8 r}{3(q-1) x_{0}^{3}} \int_{0}^{\infty} \e^{-w} w^{3} dw\nn\\
    &=\frac{-16 r}{(q-1) x_{0}^{3}}.
    }

As to the second term of the right-hand side of \eqref{split-integral-a}, we have now, see \eqref{equation-Guido-19},
    \eqn{
    \label{equation-Guido-36}
    \int_{2 x_{0}}^{\infty} x^{r} a(x) d x=-\frac{\left(2 x_{0}\right)^{r}}{q^{2}}-r J,
    }
with
    \eqan{
    \label{equation-Guido-37}
    J&=\int_{2 x_{0}}^{\infty} x^{r-1} \frac{\e^{x}}{\left(\e^{x}+q-1\right)^{2}} d x \leq\left(2 x_{0}\right)^{r-1} \int_{2 x_{0}}^{\infty} \frac{\e^{x}}{\left(\e^{x}+q-1\right)^{2}} d x\\
    &=\left(2 x_{0}\right)^{r-1} \frac{-1}{\e^{x}+q-1}\Big|_{2 x_{0}} ^{\infty}=\frac{\left(2 x_{0}\right)^{r-1}}{(q-1) q},
    \nn
    }
where we have again used that $\e^{x_{0}}=q-1$.

Summarizing, we have now from \eqref{equation-Guido-27}, \eqref{equation-Guido-28}, \eqref{equation-Guido-33} and \eqref{equation-Guido-35} that
    \eqn{
    \label{equation-Guido-38}
    \int_{0}^{2 x_{0}} x^{r} a(x) d x=-2 r x_{0}^{r}\left(\frac{q-2}{2 q(q-1) x_{0}}-\frac{1}{q^{2}}\right)+E,
    }
with
    \eqn{
    \label{equation-Guido-39}
    0 \leq E \leq-r x_{0}^r\left(\frac{x_{0}}{(q-1)^{3 / 2}}(8 \log 2-1)+\frac{16}{(q-1) x_{0}^{3}}\right),
    }
and, by \eqref{equation-Guido-36} and \eqref{equation-Guido-37}, that
    \eqn{
    \label{equation-Guido-40}
    \int_{2 x_{0}}^{\infty} x^{r} a(x) d x=-\frac{\left(2 x_{0}\right)^{r}}{q^{2}}+F,
    }
with
    \eqn{
    \label{equation-Guido-41}
    0 \leq F \leq-r \frac{\left(2 x_{0}\right)^{r-1}}{(q-1) q}=-r x_{0}^{r} \frac{1}{2(q-1) q x_{0}}.
    }
Thus, from \eqref{split-integral-a}, we get
    \eqn{
    \label{equation-Guido-42}
    0=-2 r x_{0}^{r}\left(\frac{q-2}{2 q(q-1) x_{0}}-\frac{1}{q^{2}}\right)-\frac{\left(2 x_{0}\right)^r}{q^{2}}+E+F,
    }
with $E$ and $F$ bounded by \eqref{equation-Guido-39} and \eqref{equation-Guido-41}. Now note that, as $x_{0} \rightarrow \infty,$ $q=\e^{x_{0}}+1 \rightarrow \infty$, and
    \eqn{
    \label{equation-Guido-43}
    \frac{q-2}{2 q(q-1) x_{0}}-\frac{1}{q^{2}} \sim \frac{1}{2 q x_{0}},
    }
    \eqn{
    \label{equation-Guido-44}
    \frac{x_{0}}{(q-1)^{3 / 2}}(8 \log 2-1) \sim \frac{x_{0}}{q^{3 / 2}}(8 \log 2-1)=\frac{1}{2 q x_{0}} O\left(\frac{x_{0}^{2}}{q^{1 / 2}}\right),
    }
    \eqn{
    \label{equation-Guido-45}
    \frac{16}{(q-1) x_{0}^{3}} \sim \frac{16}{q x_{0}^{3}}=\frac{1}{2 q x_{0}} O\left(\frac{1}{x_{0}^{2}}\right),
    }
    \eqn{
    \label{equation-Guido-46}
    \frac{1}{2(q-1) q x_{0}} \sim \frac{1}{2 q^{2} x_{0}}=\frac{1}{2 q x_{0}} O\left(\frac{1}{q}\right).
    }
Therefore, $E$ and $F$ can be absorbed in the first term of the right-hand side of \eqref{equation-Guido-42}, at the expense of a relative error of order $1/ x_{0}^{2}$, and we get
    \eqn{
    \label{equation-Guido-47}
    0=-2 r x_{0}^r\left(\frac{q-2}{2 q(q-1) x_{0}}-\frac{1}{q^{2}}\right)\left(1+O\left(\frac{1}{x_{0}^{2}}\right)\right)-\frac{(2 x_0)^r}{q^{2}}.
    }
We find from \eqref{equation-Guido-47} that
    \eqan{
    \label{equation-Guido-48}
    -r 
    & =\frac{\left(2 x_{0}\right)^{r}}{q^{2}} /\left(2 x_{0}^{r}\left(\frac{q-2}{2 q(q-1) x_{0}}-\frac{1}{q^{2}}\right)\right) \left(1+O\left(\frac{1}{x_{0}^{2}}\right)\right)\\
    & =\frac{2^r}{\frac{(q-2) q}{q(q-1) x_{0}}-2}\left(1+O\left(\frac{1}{x_{0}^{2}}\right)\right)\nn\\
    & =\frac{2^r x_{0}}{q-1-(q-1)^{-1}-2 x_{0}}\left(1+O\left(\frac{1}{x_{0}^{2}}\right)\right)\nn.
    }
Hence, since $r<0$, we have that $0<2^r<1$, and so $-r=O\left(x_{0} / q\right)$. Consequently, $2^{r}=1+O\left(1/x_{0}^{2}\right)$, and we get the result in \eqref{r(q)-asymptoics-sharp} from \eqref{equation-Guido-46}.

\invisible{I have estimated $r(q) \approx-0.6$ for $q=3, r(q) \approx-0.3$, $q=7$, and I have shown rigorously that $L_{q}(r)>0$ when $r \in(-1,0)$ is fixed and q is large enough. The proof of the latter result is rather lengthy and technical.

Desiderata could include

\begin{itemize}
  \item[$\rhd$] computation of $L_q(r)$ and of $r(q)$,
  \item[$\rhd$] showing monotonicity of $r(q), q=3,4, \cdots$.
\end{itemize}}


\section{Computing the critical value via Theorem \ref{thm-critical-value}}
\label{appC}

In this section we discuss the computation of the critical temperature via Theorem \ref{thm-critical-value}. We first present an iterative procedure that holds for any distribution of $W$ satisfying the unique zero-crossing condition in Condition \ref{cond-zer-cross}(a). We then illustrate the actual computation of the critical value for a few special weight distributions.
\subsection*{Newton iteration}
According to Theorem \ref{thm-critical-value}, if $\expec[W^2]<\infty$, $q\ge 3$ and $B=0$,
there is a unique positive root $t=t_c$
of the equation $\mathscr{K}(t)=0$, with 
$\mathscr{K}(t)$ given by \eqref{K-function},
and then $\beta'_c$ and $s(\beta'_c)$ 
are determined by \eqref{beta-s-crit}. 
\invisible{In order to bound $t_c$,} We recall the following facts:
\begin{itemize}
\item[(i)] there exists a unique $t_*>0$, which is the unique positive zero of $\calF''(t)$ and then, by \eqref{Kdoubleprime-funct}, $t_*$ is the unique positive zero of $\mathscr{K}''(t)$;  
\item[(ii)] $\mathscr{K}(t)$ is concave on $[0,t_*]$ and convex on $[t_*,\infty)$;
\item[(iii)] there exists a $t_b>t_*$ which is the unique positive zero of 
$\calF_0(t)-t\calF_0'(t)$ as in Proposition \ref{prop:tangents}, with $B=0$. By \eqref{kprime-funct} and \eqref{Kdoubleprime-funct} $t_b$ is also the unique minimizer of $\mathscr{K}(t)$, with $\mathscr{K}(t_b)<0=\mathscr{K}'(t_b)$.
\end{itemize}
Therefore, since $0<t_*<t_b<t_c$, the unique zero $t_c$ of $\mathscr{K}(t)$ can be safely computed by Newton's method when one initializes the Newton iteration with a point $t^{(0)}>t_c$.

\begin{proposition}
\label{upper-bound-t_c}
Let $\expec[W^2]<\infty$, $q\ge 3$ and $B=0$. Then the critical value $t_c$ satisfies
\eqn{
t_c\le \frac{2q\log q}{(q-1)\mathbb{E}[W]}.
}
\end{proposition}
\begin{proof}
\invisible{We have
   \eqn{
       \label{K-eqn}
    \mathscr{K}(t)=
    \frac{1}{\expec[W]}\expec\left[\log\left(\frac{\e^{tW}+q-1}{q}\right)\right]
    -\frac{q-1}{2q} t\mathscr{F}(t)-\frac{t}{q}
    }
with
   \eqn{
      \label{F-eqn}
    \mathscr{F}(t)=\expec \left[\frac{W}{\expec[W]}\frac{\e^{tW}-1}{\e^{tW}+q-1} \right].
    }}
By \eqref{K-function},
    \eqan{
    \mathscr{K}(t)
    &=
-\frac{\log q}{ \mathbb{E}[W]}
    +\frac{1}{\expec[W]}\expec\left[\log\left(\e^{tW}+q-1\right)\right]
    -\frac{q+1}{2q}t + \frac{1}{2}(q-1) t \frac{1}{\expec[W]}
    \mathbb{E}\left[\frac{W}{\e^{tW}+q-1}\right]  
    \nonumber\\
& >
    -\frac{\log q}{ \mathbb{E}[W]}
    +\frac{1}{\expec[W]} \expec[tW]  -\frac{q+1}{2q}t 
    = 
    -\frac{\log q}{ \mathbb{E}[W]}
    +\frac{q-1}{2q}t.
    }
Hence $\mathscr{K}(t)>0$ when $t>\frac{2q \log q}{(q-1)\mathbb{E}[W]}$ and since $\mathscr{K}(t_c)=0$,
this gives the result.
\end{proof}
By Proposition \ref{upper-bound-t_c} and convexity of $\mathscr{K}$ on $[t_c,\infty)$, the Newton iteration
\eqn{
\label{newton}
t^{(0)} = \frac{2q \log q}{(q-1)\mathbb{E}[W]}; \qquad t^{(j+1)} = t^{(j)} - \frac{\mathscr{K}(t^{(j)})}{\mathscr{K}'(t^{(j)})}, \quad j=0,1,\ldots
}
converges monotonically to $t_c$.

We illustrate below the possibility of computing $t_c$ for the Erd\H{o}s-R\'enyi random graph and the random graph with power-law degree distribution. 

\subsection*{Homogeneous case}
In the homogeneous setup $W=1$. By \eqref{K-function} and \eqref{def-Fcal}, 
\eqn{
\label{K-ER}
    \mathscr{K}(t)=
    \log\left(\frac{\e^{t }+q-1}{q}\right)
    -\frac{(q+1)}{2q} t + \frac{1}{2} \frac{(q-1)t}{\e^{t}+q-1}.    
    }
It is well-known (cf.\ Remark \ref{remark-ER}) that in this case $t_c=2\log(q-1)$, and it is indeed readily verified that $\mathscr{K}(2\log(q-1))=0$. The upper bound for $t_c$ given in Proposition 
\ref{upper-bound-t_c} equals $\frac{2q\log q}{q-1}$ and it is tight especially when $q$ is large. As an example, we consider the case that $q=7$ so that $t_c= 2\log 6 = 3.583518938\ldots\;$. We have, for $t>0$,
\eqn{
\mathscr{K}(t) = \log\left(\frac{\e^t+6}{7}\right)-\frac{4}{7}t + \frac{3t}{\e^t+6}, 
\qquad
\mathscr{K}'(t) = 3\left(\frac{1}{7}- \frac{1}{\e^t+6}-\frac{t\e^t}{(\e^t+6)^2}\right).
}
The Newton iteration in \eqref{newton} with $t^{(0)}=\frac{7}{3}\log 7$ yields a $t^{(j)}$ which agrees with $t_c$ to 9 decimal places from $j=5$ onward.
\invisible{
\subsection{First-order phase transition for monomial weights: Proof of Theorem \ref{thm-pt-power-law}}
\begin{proposition}\label{prop:unique_inflection}
    For $W$ having a monomial distribution, the second derivative $\frac{d^2}{dt^2}\mathscr{F}(t)$ of $\mathscr{F}$ is first positive and then negative.
\end{proposition}
\begin{proof}
\RvdH{TO DO 11: Add this proof! Remco will do this.}
\end{proof}
\subsection{First-order phase transition for power law - exponential weights: Proof of Theorem \ref{thm-pt-exponential-law}}
\begin{proposition}\label{prop:unique_inflection_exp}
    For $W$ having distribution $f_W(w)=c\, w^p\, \e^{-b\, w^q} \1_{w\ge 0}$, with $p>0$, $b>0$, $q>0$, the second derivative $\frac{d^2}{dt^2}\mathscr{F}(t)$ of $\mathscr{F}$ is first positive and then negative.
\end{proposition}
\begin{proof}
Recalling eq. \eqref{eq:second_derivative_gen}, we write the average with respect to the  distribution $f_W(w)=c\, w^p\, \e^{-b\, w^q} \1_{w\ge 0}$
    \be\label{eq:second_derivative_expw}
    \frac{d^2}{dt^2}\mathscr{F}(t)=\frac{q\, c}{\expec[W]}\int_0^\infty w^{3+p}\,  a(t\, w)\, \e^{-b\, w^q} dw,
    \ee
    where the function
    \be\label{eq:a}
    a(x)= \frac{(q-1)\e^{x}-\e^{2x}}{(\e^{x}+q-1)^3} 
    \ee
vanishes for $x_0=\log (q-1)$ and is $a(x)>0$ for $x<x_0$ and $a(x)<0$ for $x>x_0$. We observe that the integral in \eqref{eq:second_derivative_expw} is positive for $t=0$ and negative for large $t$. The last fact can be shown by writing the integral as 
$$
\int_0^\infty w^{3+p}\,  a(t\, w)\, \e^{-b\, w^q} dw = \int_0^{x_0/t} w^{3+p}\,  a(t\, w)\, \e^{-b\, w^q} dw + \int_{x_0/t}^\infty w^{3+p}\,  a(t\, w)\, \e^{-b\, w^q} dw,
$$
and noting that the first term in the r.h.s. is positive and vanishing for $t\to +\infty$ (being the integrand a bounded function) and the second one is negative and decreasing as $t\to +\infty$. Then, by continuity, there are points $t$ at which the integral in \eqref{eq:second_derivative_expw} vanishes, we denote $t_0>0$ the smallest one.  We want to prove that this point is in fact the unique zero of the integral.\\
By changing variable, we can rewrite \eqref{eq:second_derivative_expw}  as 
$$
  \frac{d^2}{dt^2}\mathscr{F}(t)=\frac{q\, c}{\expec[W]\, t^{p+4}} \Phi(t),
$$
where
$$
\Phi(t)= \int_0^\infty x^{3+p}\,  a(x)\, \exp(-b\, (x/t)^q) dw.
$$
Now we claim that $\Phi'(t)<0$ whenever $\Phi(t)\le 0$. Then, as a consequence, we will have that $\Phi(t)< \Phi(t_0)=0$ for all $t>t_0$. This shows that $t_0$ is the unique zero of \eqref{eq:second_derivative_expw}. In order to prove the claim, consider the derivative:
$$
\Phi'(t)= \frac{b\, q}{t^{q+1}}  \int_0^\infty x^{3+p+q}\,  a(x)\, \exp(-b\, (x/t)^q) dx.
$$
The function $a(x)$ (see the line following \eqref{eq:a}) satisfies the condition $(x^q-x_0^q)x^{3+p} a(x)<0$ for all $x\ne x_0$ which implies that:
\begin{align}
\Phi'(t) &= \frac{b\, q}{t^{q+1}}  \int_0^\infty x^{3+p+q}\,  a(x)\, \exp(-b\, (x/t)^q) dx < \frac{b\, q}{t^{q+1}}  x_0^q\int_0^\infty x^{3+p}\,  a(x)\, \exp(-b\, (x/t)^q) dx\\
&< \frac{b\, q}{t^{q+1}}  x_0^q\  \Phi(t).
\end{align}
Therefore $\Phi(t)$ is decreasing when $\Phi(t)\le 0$, as claimed.
\end{proof}}
\subsection*{Pareto case with $\tau>4$}
In this section, we specialize to the Pareto distribution, for which the probability density function is given by \eqref{Pareto} and $
    \expec[W]=(\tau-1)/(\tau-2)$.
By \eqref{K-function} and \eqref{def-Fcal},
\eqn{
\mathscr{K}^{\text{Par}}(t)
= 
-\frac{\log q}{ \mathbb{E}[W]}
    -\frac{q+1}{2q}t 
    +\frac{1}{\expec[W]}\left(\expec\left[\log\left(\e^{tW}+q-1\right)\right]
    + \frac{1}{2}(q-1) t 
    \mathbb{E}\left[\frac{W}{\e^{tW}+q-1}\right]\right).
}
\invisible{We work out
\eqn{
\mathbb{E}\left[\frac{W}{\e^{tW}+q-1}\right] =
\int_{1}^{\infty}\frac{(\tau-1)w^{-\tau+1}}{\e^{tw}+q-1} dw.
}
}
By partial integration and basic manipulation, we compute
\invisible{\eqan{
\mathbb{E}[\log(\e^{tW}+q-1)]
&=& \int_{1}^{\infty}(\tau-1)w^{-\tau}\log(\e^{tw}+q-1)dw
\nonumber\\
&=& \int_{1}^{\infty}\log(\e^{tw}+q-1)d(-w^{-\tau+1})
\nonumber\\
&=& \log(\e^t+q-1) + t \int_{1}^{\infty} \frac{w^{-\tau+1}\e^{tw}}{\e^{tw}+q-1}dw 
\nonumber\\
&=& \log(\e^{t}+q-1)+\frac{t}{\tau-2} -(q-1)t\int_{1}^{\infty} \frac{w^{-\tau+1}}{\e^{tw}+q-1}dw,
}
}
\eqan{
\mathbb{E}[\log(\e^{tW}+q-1)]
&= \int_{1}^{\infty}(\tau-1)w^{-\tau}\log(\e^{tw}+q-1)dw
= \int_{1}^{\infty}\log(\e^{tw}+q-1)d(-w^{-\tau+1})\nonumber\\
&= \log(\e^t+q-1) + t \int_{1}^{\infty} \frac{w^{-\tau+1}\e^{tw}}{\e^{tw}+q-1}dw 
\nonumber\\
&= \log(\e^{t}+q-1)+\frac{t}{\tau-2} -(q-1)t\int_{1}^{\infty} \frac{w^{-\tau+1}}{\e^{tw}+q-1}dw,
}
and, since
\eqn{
\mathbb{E}\left[\frac{W}{\e^{tW}+q-1}\right] =
\int_{1}^{\infty}\frac{(\tau-1)w^{-\tau+1}}{\e^{tw}+q-1} dw,
}
we finally get
  \eqan{
  \label{K-Par}
    \mathscr{K}^{\text{Par}}(t)&=
    \frac{\tau-2}{\tau-1}\log\left(\frac{\e^{t }+q-1}{q}\right)
    -\left(\frac{q+1}{2q}-\frac{1}{\tau-1}\right) t\nonumber\\
    &\qquad + \frac{(\tau-3)(q-1)t}{2(\tau-1)} \int_{1}^{\infty}\frac{(\tau-2)w^{-\tau+1}}{\e^{t w}+q-1}dw.    
    }
Observe the close resemblance of $\mathscr{K}^{\text{Par}}$ in \eqref{K-Par} for the Pareto case with $\mathscr{K}= \mathscr{K}^{\text{ER}}$ in \eqref{K-ER} for the Erd\H{o}s-R\'enyi case.
Indeed, as $\tau\to\infty$, the three terms on the r.h.s.\ of \eqref{K-Par} converge to the corresponding terms on the r.h.s.\ of \eqref{K-ER}.
\smallskip

\begin{figure}
\centering
\includegraphics[width=10cm]{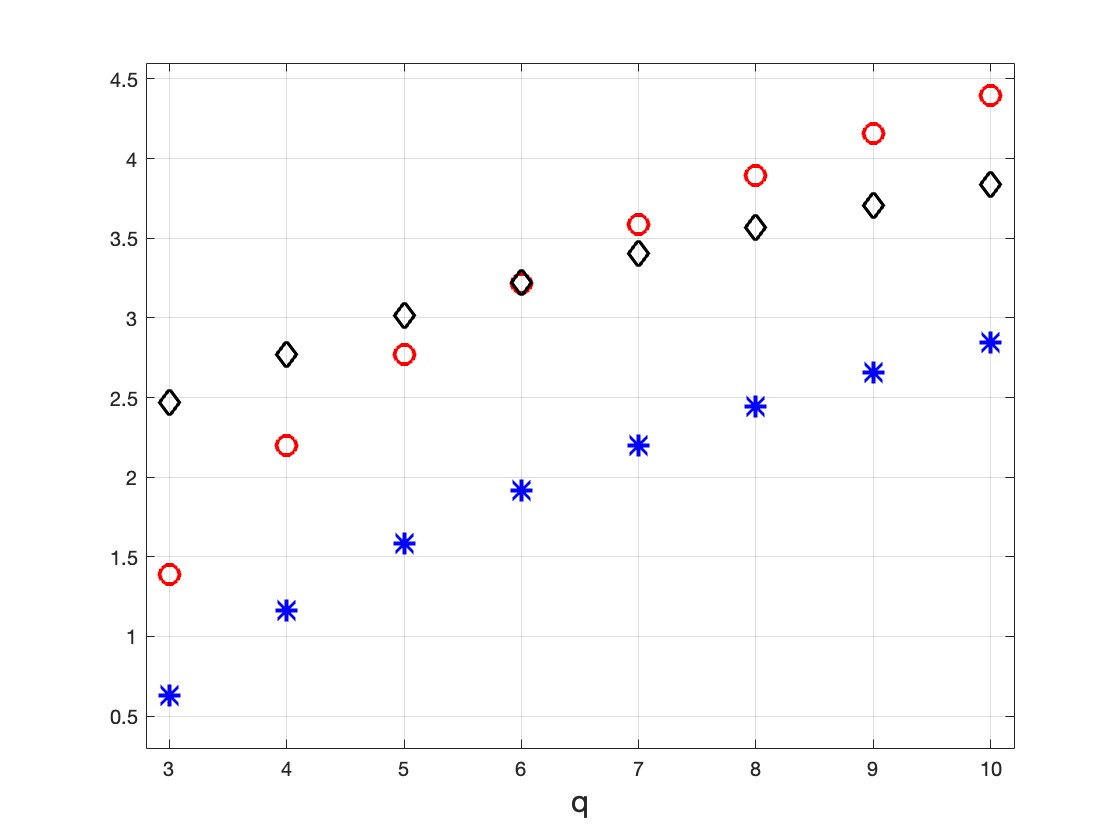}
\caption{The critical $t_c=2 \log(q-1)$ ($\circ$) in the homogeneous case for $q=3,4,\ldots, 10$ and, for the same $q$-values, the critical $t_c$ ($*$) for the Pareto case with $\tau=5$, as computed by the Newton iteration \eqref{newton}, initialized by the upper bound  $\frac{2q\log q}{(q-1)\mathbb{E}[W]}$
($\diamond$)  of Proposition  \ref{upper-bound-t_c}.}
  \label{fig:crit-pareto}
\end{figure}

The remaining integral on the r.h.s.\ of \eqref{K-Par} does not admit a closed-form expression, and so it can be expected that there is no closed-form solution $t=t_c$ of the equation $\mathscr{K}^{\text{Par}}(t_c)=0$. For $\tau \ge 4$, $\calF(t)$ is first positive and then negative and so the unique positive zero $t_c$ of $\mathscr{K}(t)$ can be computed by Newton iteration in \eqref{newton} with initialization $t^{(0)}= \frac{\tau-2}{\tau-1}\frac{2q \log q}{q-1}$. We compute
\eqan{
\mathscr{K}'(t)
&=
    \frac{\tau-2}{\tau-1}\frac{\e^t}{\e^{t }+q-1}
    -\left(\frac{q+1}{2q}-\frac{1}{\tau-1}\right)  + \frac{(\tau-3)(q-1)}{2(\tau-1)} \int_{1}^{\infty}\frac{(\tau-2)w^{-\tau+1}}{\e^{t w}+q-1}dw
    \nonumber\\
&\qquad-\frac{(\tau-3)(q-1)t}{2(\tau-1)} \int_{1}^{\infty}\frac{(\tau-2)w^{-\tau+2}\e^{tw}}{(\e^{t w}+q-1)^2}dw.}
The computation of $\mathscr{K}(t)$ and $\mathscr{K}'(t)$ requires the numerical evaluation of the two integrals appearing on the r.h.s.\ of the above equation. These integrals can be evaluated using numerical routines. When we do so for the case $\tau=5,$ and $q=7$, we find that Newton iteration in \eqref{newton} gives $t^{(j)} = 2.20111...$ from $j=5$ onward; see Figure \ref{fig:crit-pareto}.    

\end{document}